\documentclass[11pt]{amsart}
\usepackage[T1]{fontenc}
\usepackage[cp852]{inputenc}
\usepackage{amsmath}
\usepackage{amssymb}
\usepackage{moreverb}
\newtheorem{thm}{Theorem}[section]
\newtheorem{cor}[thm]{Corollary}

\newtheorem{prop}[thm]{Proposition}
\newtheorem{rem}[thm]{Remark}
\newtheorem{example}[thm]{Example}
\usepackage{graphicx}

\usepackage{amsfonts}
\usepackage{amssymb}
\usepackage{eucal}

\theoremstyle{definition}
\numberwithin{equation}{section}
\begin{document}

\begin{center}
{\bf{ON SEMI-RIEMANNIAN MANIFOLDS SATISFYING 
SOME GENERALIZED EINSTEIN METRIC CONDITIONS}}
\end{center}

\vspace{2mm}

\begin{center}
Ryszard Deszcz, Ma\l gorzata G\l ogowska, Marian Hotlo\'{s},
Miroslava Petrovi\'{c}-Torga\v{s}ev, and Georges Zafindratafa
\end{center}

\vspace{2mm}

\begin{center}
{\sl{Dedicated to the memory of Professor Old\v{r}ich Kowalski (1936--2021)}}
\end{center}

%\thispagestyle{empty}
%\documentclass[12pt]{amsart}
%\documentstyle[12pt,leqno]{article}
%\usepackage[T1]{fontenc}
%\usepackage[cp852]{inputenc}
%\topmargin -10mm
%\oddsidemargin -10mm
%\textheight 220mm

%\textwidth 135mm

%\oddsidemargin  0pt
%\evensidemargin 0pt
%\marginparwidth 1in
%\marginparwidth 0in
%\marginparsep 0pt

%\topmargin 0pt
%\headheight 0pt
%\headsep 0pt
%\topskip 0pt

%\footheight 5mm
%\footskip 10mm

%\textheight 248mm
%\textwidth 160mm

%\begin{document}

%\renewcommand{\baselinestretch}{1.1}

\vspace{5mm}

\noindent
{\bf{Abstract.}} The difference tensor $R \cdot C - C \cdot R$ of a
semi-Riemannian manifold $(M,g)$, $\dim M \geq 4$, formed by its 
Riemann-Christoffel curvature tensor
$R$ and the Weyl conformal curvature tensor $C$,
under some assumptions, 
can be expressed 
as a linear combination of $(0,6)$-Tachibana tensors $Q(A,T)$, 
where $A$ is a symmetric $(0,2)$-tensor and $T$ 
a generalized curvature tensor. These conditions
form a family of generalized Einstein metric conditions.
In this survey paper we present recent results
on manifolds and submanifolds, and in particular hypersurfaces, 
satisfying such conditions.

%\section{Curvature properties of some warped product manifolds}

\section{Introduction}

Let $(M,g)$
be a semi-Riemannian manifold.
We denote by
$g$, $\nabla$, $R$, $S$, $\kappa$ and $C$, 
the metric tensor, the Levi-Civita connection, 
the Riemann-Christoffel curvature tensor, 
the Ricci tensor,
the scalar curvature and the Weyl conformal curvature tensor 
of $(M,g)$, respectively. 
Now we can define the $(0,2)$-tensors 
$S^{2}$ and $S^{3}$,
the $(0,4)$-tensors
$R \cdot S$, $C \cdot S$ and
$Q(A,B)$, 
and 
the $(0,6)$-tensors
$R \cdot R$, 
$R \cdot C$, 
$C \cdot R$, 
$C \cdot C$ 
and 
$Q(A,T)$, 
where $A$ and $B$ are symmetric $(0,2)$-tensors
and $T$ a generalized curvature tensor.
Thus, in particular, we have the 
difference tensor 
$R \cdot C - C \cdot R$,
or $C \cdot R  - R \cdot C$,
{\cite[Section 1] {DP-TVZ}}.
Furthermore for $A$ and $B$
we define their Kulkarni-Nomizu product $A \wedge B$. 
For precise definitions of the symbols used, 
we refer to Section 2
of this paper, as well as to
{\cite[Section 1] {2016_DGJZ}}, 
{\cite[Section 1] {2020_DGZ}},
{\cite[Chapter 6] {DHV2008}} and
{\cite[Sections 1 and 2] {DP-TVZ}}.

A semi-Riemannian manifold $(M,g)$, $\dim M = n \geq 2$, is said to be 
an {\sl Einstein manifold} \cite{Besse-1987},
or an {\sl Einstein space}, if at every point of $M$ 
its Ricci tensor $S$ is proportional to $g$, 
i.e., 
\begin{eqnarray}
S = \frac{\kappa}{n}\, g 
\label{2020.10.3.c}
\end{eqnarray}
on $M$,
assuming that the scalar curvature $\kappa$ is constant when $n = 2$. 
According to {\cite[p. 432] {Besse-1987}}
this condition is called the {\sl Einstein metric condition}.
Einstein manifolds form a natural subclass
of several classes of semi-Riemannian manifolds which are determined by
curvature conditions imposed on their Ricci tensor 
{\cite[Table, pp. 432-433] {Besse-1987}}.
These conditions are called 
{\sl generalized Einstein metric conditions}
{\cite[Chapter XVI] {Besse-1987}}.

Semi-Riemannian manifolds of dimension $\geq 4$ and in particular, 
hypersurfaces in spaces of constant curvature
or Chen ideal submanifolds in Euclidean spaces, 
satisfying curvature conditions of the form 

\vspace{2mm}

\noindent
$(\ast)$
the difference tensor 
$R \cdot C - C \cdot R$ 
and a finite sum of the Tachibana tensors of the form
$Q(A,T)$ are linearly dependent, 

\vspace{2mm}

\noindent
where $A = g$, $A = S$,  or $A = S^{2}$,
and 
$T = R$, $T = C$, 
$T = g \wedge g$, $T = g \wedge S$, $T = g \wedge S^{2}$,   
$T = S \wedge S$, $T = S \wedge S^{2}$, 
or $T = S^{2} \wedge S^{2}$,
were studied in several papers, see, e.g.,
\cite{{P119}, {Ch-DDGP}, {DGHZ01}, {2016_DGHZhyper}, {2016_DGJZ}, 
{DGPSS}, {2020_DGZ}, {DH-Tsukuba}, {DHS01}, {DP-TVZ}}.
Conditions of this form are also
generalized Einstein metric conditions,
see, e.g., {\cite[Section 6] {DGHSaw}} or {\cite[Section 1] {DP-TVZ}}.  
We refer to \cite{DGHSaw} for a survey of results 
on manifolds satisfying conditions of the form $(\ast)$. 

Precisely, in that paper results obtained to 2010 are presented.
In the next years studies on manifolds 
and submanifolds, and in particular hypersurfaces, 
satisfying such conditions are continued.
As it was already mentioned in Abstract,
in this paper we present a survey of recent results on that subject.

On every semi-Riemannian Einstein manifold $(M,g)$, $\dim M = n \geq 4$,
the following identity is satisfied
{\cite[Theorem 3.1] {DHS01}}  (see also {\cite[p. 107] {DP-TVZ}})
\begin{eqnarray*}
R \cdot C - C \cdot R = \frac{\kappa}{ (n-1)n }\, Q(g,C) .
%\label{Einstein011}
\end{eqnarray*}
This, by  
(\ref{2020.10.3.c}) and 
$\frac{\kappa}{ (n-1)n } = \frac{\kappa}{ n-1 } - \frac{\kappa}{ n }$, 
turns into 
\begin{eqnarray}
C \cdot R - R \cdot C = Q(S,C) - \frac{\kappa}{ n-1 }\, Q(g,C) .
\label{Roterformula}
\end{eqnarray}
We mention that there are non-Einstein and non-conformally flat 
semi-Riemannian manifolds satisfying  (\ref{Roterformula}).
Namely, every Roter space satisfies (\ref{Roterformula}) (see Section 4).

Let $(M,g)$ be a semi-Riemannian manifold of dimension $n \geq 3$.
We define the subsets ${\mathcal{U}}_{R}$ and ${\mathcal U}_{S}$ of $M$ by 
${\mathcal{U}}_{R}  = 
\{x \in M\, | \, 
R - \frac{\kappa }{2 (n-1) n}\, g \wedge g \neq 0\ \mbox {at}\ x \}$
and 
${\mathcal U}_{S} =  \{x \in M\, | \, 
S - \frac{\kappa }{n}\, g \neq 0\ \mbox {at}\ x \}$, respectively. 
If $n \geq 4$ then
we define the set ${\mathcal U}_{C} \subset M$ as the set of all points 
of $(M,g)$ at which $C \neq 0$.
We note that on any semi-Riemannian manifold $(M,g)$, 
$\dim M = n \geq 4$, we have (see, e.g., \cite{DGHHY})
\begin{eqnarray}
{\mathcal{U}}_{S} \cup {\mathcal{U}}_{C} &=& {\mathcal{U}}_{R} . 
\label{dghhy}
\end{eqnarray}

In Section 3 we present a survey on semi-Riemannian manifolds 
satisfying curvature conditions known 
as pseudosymmetry type curvature conditions. In particular, 
curvature conditions of the form $(\ast)$ are conditions of this kind.

An extension of the class of Einstein manifolds also form
quasi-Einstein, $2$-quasi-Einstein and partially Einstein manifolds.
These manifolds satisfy some
pseudosymmetry type curvature conditions (see Section 4).

In Section 5 we consider warped product manifolds
$\overline{M} \times _{F} \widetilde{N}$ 
with a $2$-dimensional semi-Riemannian manifold 
$(\overline{M},\overline{g})$,
a warping function $F$
and an $(n-2)$-dimensional semi-Riemannian manifold 
$(\widetilde{N},\widetilde{g})$, $n \geq 4$, 
assuming that 
$(\widetilde{N},\widetilde{g})$
is a space of constant curvature 
when $n \geq 5$. 
Such manifolds satisfy some
pseudosymmetry type curvature conditions 
(see Theorems 5.1 and 5.2).
We mention that
the warped product manifolds
$\overline{M} \times _{F} \widetilde{N}$,
with a $2$-dimensional Riemannian manifold 
$(\overline{M},\overline{g})$ 
and an $(n-2)$-dimensional unit sphere ${\mathbb S}^{n-2}$, $n \geq 4$,
where the warping function $F$
is a solution of some second order
quasilinear elliptic partial differential equation in the plane
are related to Chen ideal submanifolds 
(see Section 11 for details).

Section 6 contains results on semi-Riemannian manifolds 
$(M,g)$ of dimension $\geq 4$
satisfying 
the following
generalized Einstein metric condition
on
${\mathcal U}_{S} \cap {\mathcal U}_{C} \subset M$
\begin{eqnarray}
R\cdot C - C \cdot  R &=&  L\, Q(S,C),
\label{advances4}
\end{eqnarray}
where $L$ is some function on this set.

Let $N_{s}^{n+1}(c)$, $n  \geq 4$, 
be a semi-Riemannian space of constant curvature 
with signature $(s,n+1-s)$, where
$c = \frac{\widetilde{\kappa}}{n (n+1)}$ and
$\widetilde{\kappa}$ is its scalar curvature. 
Let $M$ be a connected hypersurface
isometrically immersed in $N_{s}^{n+1}(c)$.
We denote by ${\mathcal U}_{H} \subset M$ the set of all points 
at which the tensor 
$H^{2}$ is not a linear combination of the metric tensor $g$ 
and the second fundamental tensor $H$ of $M$. We can verify that 
${\mathcal U}_{H} \subset {\mathcal U}_{S} \cap {\mathcal U}_{C} \subset M$ 
(see, e.g., {\cite[p. 137] {G108}}). 
We note that 
\begin{eqnarray}
H^{2} &=& \alpha \, H + \beta \, g  
\label{h2277aa}
\end{eqnarray}
on $({\mathcal U}_{S} \cap {\mathcal U}_{C}) \setminus {\mathcal U}_{H}$, 
where $\alpha $ and $\beta $ are some functions defined on this set.

As it was stated in {\cite[Corollary 4.1] {R99}},
for a hypersurface $M$ in $N^{n+1}_{s}(c)$, $n \geq 4$, 
if at every point of ${\mathcal U}_{H} \subset M$ one of the tensors
$R \cdot C$, $C \cdot R$ 
or $R \cdot C - C \cdot R$ is a linear combination of the tensor 
$R \cdot R$ 
and a finite sum of the Tachibana tensors of the form $Q(A,T)$, 
where $A$ is a symmetric $(0,2)$-tensor 
and $T$ a generalized curvature tensor, 
then the following equation is satisfied on ${\mathcal U}_{H}$
\begin{eqnarray}
H^{3} &=& \mathrm{tr}(H)\, H^{2} + \psi \, H + \rho \, g ,
\label{dgss01}
\end{eqnarray}
where $\psi$ and $\rho$ are some functions on ${\mathcal U}_{H}$. 
We also mention that if the condition
\begin{eqnarray*} 
R \cdot C - C \cdot R &=& Q(g, T), 
\end{eqnarray*}
where $T$ is a generalized curvature tensor,
holds on the set ${\mathcal U}_{H}$ of a hypersurface 
$M$ in $N^{n+1}_{s}(c)$, $n \geq 4$, 
then the tensor $T$ is a linear combination of the tensors:
$R$, $g \wedge g$, $g \wedge S$, $S \wedge S$ and
$g \wedge S^{2}$ {\cite[Theorem 5.2 (iii)] {DGPSS}}.

In Section 7 we present results on $2$-quasi-umbilical
hypersurfaces in $N_{s}^{n+1}(c)$, $n  \geq 4$,
contained in \cite{DGP-TV02}. 
In particular
(see Theorem 7.3),
we present curvature properties of a minimal 
$2$-quasi-umbilical hypersurface $M$
in an Euclidean space $\mathbb{E}^{n+1}$, $n \geq 4$, 
having exactly three distinct principal curvatures 
$\lambda _{1}$, $\lambda _{2}$ and $\lambda _{3}$ 
satisfying at every point of $M$: 
$\lambda _{1} = 0$, $\lambda _{2} = - (n-2) \lambda $ 
and  
$\lambda _{3} = \lambda _{4} = \ldots = \lambda _{n} = \lambda \neq 0$,
where $\lambda $ is a non-zero function on $M$.
It is easy to check that $\mathrm{tr}(H) = 0$ and 
\begin{eqnarray}
H^{3} &=& \phi \, H^{2} + \psi \, H 
\label{dgss01ab}
\end{eqnarray}
on ${\mathcal U}_{H} \subset M$,
where $\phi = - (n-3) \lambda $ and $\psi = (n-2) \lambda ^{2}$.
Because $\mathrm{tr} (H) \neq \phi$ at every point of ${\mathcal U}_{H}$, 
(\ref{dgss01}) is not satisfied on ${\mathcal U}_{H}$.
Now, in view of the presented above result 
(i.e., {\cite[Corollary 4.1] {R99}}), 
we conclude that the difference tensor $R \cdot C - C \cdot R$ 
of the considered hypersurface $M$
cannot be expressed on ${\mathcal U}_{H}$ as a linear combination of
certain Tachibana tensors  of the form $Q(A,T)$, 
where $A$ is a symmetric $(0,2)$-tensor and $T$ a generalized curvature tensor.
The tensor $R \cdot C - C \cdot R$ of that hypersurface we can express, 
for instance, by (\ref{2023.05.20.a}).  
In this section we also present results on type number two hypersurfaces in $N_{s}^{n+1}(c)$, $n  \geq 3$.

Results on hypersurfaces in $N_{s}^{n+1}(c)$, $n  \geq 4$,
satisfying some special generalized Einstein metric conditions 
of the form $(\ast)$
are given in Sections  8, 9 and 10. 
In Section 8 we present results on hypersurfaces $M$ 
in a semi-Riemannian space of constant curvature 
$N_{s}^{n+1}(c)$, $n \geq 4$,
satisfying (\ref{dgss01}) on ${\mathcal{U}}_{H} \subset M$,
for some functions 
$\psi $ and $\rho $ on ${\mathcal{U}}_{H}$.
Among other things  it was stated
that some  
generalized Einstein metric condition of the form $(\ast)$
is satisfied on ${\mathcal{U}}_{H}$ (see Theorem 8.3).
Sections 9 and 10 contain results on 
quasi-Einstein and non-quasi-Einstein 
hypersurfaces in semi-Riemannian spaces of constant curvature 
$N_{s}^{n+1}(c)$, $n \geq 4$, satisfying 
on ${\mathcal U}_{S} \cap {\mathcal U}_{C} \subset M$
the conditions:
(\ref{advances4}) and
\begin{eqnarray}
R \cdot C - C \cdot R &=& L_{1}\, Q(S,C) + L_{2}\, Q(g,C) ,
\label{cond01}
\end{eqnarray}
where $L$, $L_{1}$ and $L_{2}$ are some functions defined on this set.

In Section 11 we present results on non-Einstein and non-conformally flat
Chen ideal submanifolds satisfying some
generalized Einstein metric conditions of the form $(\ast)$.
For instance, Theorems 11.5, 11.8 and 11.11 
(see Theorems 5, 6 and 7 of \cite{DP-TVZ}) 
contain results on Chen ideal submanifolds 
$M$ in $\mathbb{E}^{n+m}$, $n \geq 4$, $m \geq 1$,
satisfying the following conditions 
of the form $(\ast)$:
(\ref{cond01}),
\begin{eqnarray}
R \cdot C - C \cdot R &=& L_3\, Q(g  , R) +  L_4\, Q(S  , R) ,  
\label{2023.05.20.b}\\
%\end{eqnarray}
%\begin{eqnarray}
R \cdot C - C \cdot R 
&=& L_5\, Q(g  , g \wedge S) +  L_6\, Q(S  , g \wedge S) ,
\label{2023.05.20.d}
\end{eqnarray}
respectively,
where $L_{3}$, $L_{4}$, $L_{5}$ and $L_{6}$ are some functions defined on $M$.

\section{Preliminaries.}

Throughout this paper, all manifolds are assumed 
to be connected paracompact
mani\-folds of class $C^{\infty }$. Let $(M,g)$, $\dim M = n \geq 3$,
be a semi-Riemannian manifold, and let $\nabla$ 
be its Levi-Civita connection and $\Xi (M)$ the Lie
algebra of vector fields on $M$. We define on $M$ the endomorphisms 
$X \wedge _{A} Y$ and ${\mathcal{R}}(X,Y)$ of $\Xi (M)$ by
\begin{eqnarray*}
(X \wedge _{A} Y)Z 
&=& 
A(Y,Z)X - A(X,Z)Y, \\ 
{\mathcal R}(X,Y)Z 
&=& 
\nabla _X \nabla _Y Z - \nabla _Y \nabla _X Z - \nabla _{[X,Y]}Z ,
\end{eqnarray*}
respectively, where 
$X, Y, Z \in \Xi (M)$
and 
$A$ is a symmetric $(0,2)$-tensor on $M$. 
The Ricci tensor $S$, the Ricci operator ${\mathcal{S}}$ 
and the scalar curvature
$\kappa $ of $(M,g)$ are defined by 
\begin{eqnarray*}
S(X,Y)\ =\ \mathrm{tr} \{ Z \rightarrow {\mathcal{R}}(Z,X)Y \} ,\ \  
g({\mathcal S}X,Y)\ =\ S(X,Y) ,\ \  
\kappa \ =\ \mathrm{tr}\, {\mathcal{S}},
\end{eqnarray*}
respectively. 
The endomorphism ${\mathcal{C}}(X,Y)$ is defined by
\begin{eqnarray*}
{\mathcal C}(X,Y)Z  &=& {\mathcal R}(X,Y)Z 
- \frac{1}{n-2}(X \wedge _{g} {\mathcal S}Y + {\mathcal S}X \wedge _{g} Y
- \frac{\kappa}{n-1}X \wedge _{g} Y)Z .
\end{eqnarray*}
Now the $(0,4)$-tensor $G$, 
the Riemann-Christoffel curvature tensor $R$ and
the Weyl conformal curvature tensor $C$ of $(M,g)$ are defined by
\begin{eqnarray*}
G(X_1,X_2,X_3,X_4) &=& g((X_1 \wedge _{g} X_2)X_3,X_4) ,\\
R(X_1,X_2,X_3,X_4) &=& g({\mathcal R}(X_1,X_2)X_3,X_4) ,\\
C(X_1,X_2,X_3,X_4) &=& g({\mathcal C}(X_1,X_2)X_3,X_4) ,
\end{eqnarray*}
respectively, where $X_1,X_2,\ldots \in \Xi (M)$.  
For a symmetric $(0,2)$-tensor $A$ we denote by 
${\mathcal{A}}$ the endomorphism related to $A$ by 
$g({\mathcal{A}}X,Y) = A(X,Y)$.
The $(0,2)$-tensors
$A^{p}$, $p = 2, 3, \ldots $, are defined by 
$A^{p}(X,Y) = A^{p-1} ({\mathcal{A}}X, Y)$,
assuming that $A^{1} = A$. In this way, for $A = S$ and 
${\mathcal{A}} = {\mathcal S}$
we get the tensors $S^{p}$, 
$p = 2, 3, \ldots $, assuming that $S^{1} = S$.

Let ${\mathcal B}$ be a tensor field sending any $X, Y \in \Xi (M)$ 
to a skew-symmetric endomorphism 
${\mathcal B}(X,Y)$, 
and let $B$ be
the $(0,4)$-tensor associated with ${\mathcal B}$ by
\begin{eqnarray}
B(X_1,X_2,X_3,X_4) &=& 
g({\mathcal B}(X_1,X_2)X_3,X_4) .
\label{DS5}
\end{eqnarray}
The tensor $B$ is said to be a {\sl{generalized curvature tensor}}  if the
following two conditions are fulfilled:
\begin{eqnarray*}
& &
B(X_1,X_2,X_3,X_4)\ =\ B(X_3,X_4,X_1,X_2),\\
& &  
B(X_1,X_2,X_3,X_4) + B(X_2,X_3,X_1,X_4) + B(X_3,X_1,X_2,X_4) \ =\ 0 . 
\end{eqnarray*}
For ${\mathcal B}$ as above, let $B$ be again defined by (\ref{DS5}). 
We extend the endomorphism ${\mathcal B}(X,Y)$
to a derivation ${\mathcal B}(X,Y) \cdot \, $ of the algebra 
of tensor fields on $M$,
assuming that it commutes with contractions and 
${\mathcal B}(X,Y) \cdot \, f  = 0$ for any smooth function $f$ on $M$. 
Now for a $(0,k)$-tensor field $T$,
$k \geq 1$, we can define the $(0,k+2)$-tensor $B \cdot T$ by
\begin{eqnarray*}
& & (B \cdot T)(X_1,\ldots ,X_k,X,Y) \ =\ 
({\mathcal B}(X,Y) \cdot T)(X_1,\ldots ,X_k)\\  
&=& - T({\mathcal{B}}(X,Y)X_1,X_2,\ldots ,X_k)
- \cdots - T(X_1,\ldots ,X_{k-1},{\mathcal{B}}(X,Y)X_k) .
\end{eqnarray*}
If $A$ is a symmetric $(0,2)$-tensor then we define the
$(0,k+2)$-tensor $Q(A,T)$ by
\begin{eqnarray*}
& & Q(A,T)(X_1, \ldots , X_k, X,Y) \ =\
(X \wedge _{A} Y \cdot T)(X_1,\ldots ,X_k)\\  
&=&- T((X \wedge _A Y)X_1,X_2,\ldots ,X_k) 
- \cdots - T(X_1,\ldots ,X_{k-1},(X \wedge _A Y)X_k) .
\end{eqnarray*}
In this manner we obtain the $(0,6)$-tensors 
$B \cdot B$ and $Q(A,B)$.

Substituting in the above formulas 
${\mathcal{B}} = {\mathcal{R}}$ or ${\mathcal{B}} = {\mathcal{C}}$, 
$T=R$ or $T=C$ or $T=S$, $A=g$ or $A=S$ 
we get the tensors $R\cdot R$, $R\cdot C$, $C\cdot R$, 
$C\cdot C$, $R\cdot S$, 
$Q(g,R)$, $Q(S,R)$, $Q(g,C)$,  $Q(S,C)$, 
and $Q(g,S)$, $Q(g,S^{2})$,  $Q(S,S^{2})$.

For a symmetric $(0,2)$-tensor $E$ and a $(0,k)$-tensor $T$, $k \geq 2$, we
define their {\sl{Kulkarni-Nomizu tensor}} $E \wedge T$ by 
(see, e.g., {\cite[Section 2] {DGHHY}}) 
\begin{eqnarray*}
& &(E \wedge T )(X_{1}, X_{2}, X_{3}, X_{4}; Y_{3}, \ldots , Y_{k})\\
&=&
E(X_{1},X_{4}) T(X_{2},X_{3}, Y_{3}, \ldots , Y_{k})
+ E(X_{2},X_{3}) T(X_{1},X_{4}, Y_{3}, \ldots , Y_{k} )\\
& &
- E(X_{1},X_{3}) T(X_{2},X_{4}, Y_{3}, \ldots , Y_{k})
- E(X_{2},X_{4}) T(X_{1},X_{3}, Y_{3}, \ldots , Y_{k}) .
\end{eqnarray*}
It is obvious that the following tensors 
are generalized curvature tensors: $R$, $C$ and 
$E \wedge F$, where $E$ and $F$ are symmetric $(0,2)$-tensors. 
We have 
\begin{eqnarray}
C &=& R - \frac{1}{n-2}\, g \wedge S + \frac{\kappa }{(n-2) (n-1)} \, G , 
\label{Weyl}\\
G &=& \frac{1}{2}\, g \wedge g ,
\label{2020.12.08.a}
\end{eqnarray}
and (see, e.g., {\cite[Lemma 2.2(i)] {DGHHY}})
\begin{eqnarray}
&(a)&
\ \
Q(E, E \wedge F) \ =\ - \frac{1}{2}\, Q(F, E \wedge E ) ,\nonumber\\
&(b)&
\ \
E \wedge Q(E,F) \ =\ - \frac{1}{2}\, Q(F, E \wedge E ).
\label{DS7}
\end{eqnarray}
By an application of (\ref{DS7})(a) we obtain on $M$ the identities
\begin{eqnarray}
Q(g, g \wedge S) &=& - Q(S,G) ,
\ \ Q(S, g \wedge S)\ = \ - \frac{1}{2}\, Q(g, S \wedge S) .  
\label{dghz01}
\end{eqnarray}
Further, by making use of (\ref{Weyl}), (\ref{2020.12.08.a})
and (\ref{dghz01}), we get immediately
\begin{eqnarray*}
Q(g,C) &=& Q(g,R) - \frac{1}{n-2}\, Q(g,g \wedge S) 
+ 
 \frac{\kappa }{(n-2) (n-1)} \, Q(g,G) \nonumber\\
&=& Q(g,R) - \frac{1}{n-2}\, Q(g,g \wedge S) ,\nonumber\\ 
Q(S,C) &=& Q(S,R) - \frac{1}{n-2}\, Q(S,g \wedge S) 
+ 
 \frac{\kappa }{(n-2) (n-1)} \, Q(S,G) \nonumber\\ 
&=&  
Q(S,R) + \frac{1}{2 (n-2)}\, Q(g,S \wedge S) 
- \frac{\kappa }{(n-2) (n-1)} \, Q(g,g \wedge S) .
\end{eqnarray*}

Let $E_{1}$, $E_{2}$ and $F$ be
symmetric $(0,2)$-tensors.
We have (see, e.g., 
{\cite[Lemma 2.1(i)] {Ch-DDGP}} and references therein) 
\begin{eqnarray}
E_{1} \wedge Q(E_{2},F) + E_{2} \wedge Q(E_{1},F) 
+ Q(F, E_{1} \wedge E_{2}) &=& 0 .
\label{2021.11.30.a4}
\end{eqnarray} 
From (\ref{2021.11.30.a4}) we get easily
(see also {\cite[Lemma 2.2(iii)] {DGHHY}}
and references therein)
\begin{eqnarray*}
Q(F, E \wedge E_{1}) + Q(E, E_{1} \wedge F)
+ Q(E_{1}, F \wedge E) &=& 0 .
%\label{2021.11.30.a456}
\end{eqnarray*}

We denote by $g_{ij}$, $g^{ij}$, $R_{hijk}$, $S_{ij}$, $C_{hijk}$,
$(R \cdot S)_{hklm}$, $Q(g,S)_{hklm}$, 
$Q(g,R)_{hijklm}$, $Q(S,R)_{hijklm}$,
$(R \cdot C)_{hijklm}$ and $(C \cdot R)_{hijklm}$,
the local components of the tensors $g$, $g^{-1}$, $R$, $S$, $C$,
$R \cdot S$, $Q(g,S)$, $Q(g,R)$, $Q(S,R)$, $R \cdot C$ and $C \cdot R$,
respectively. 
On every semi-Riemannian manifold $(M,g)$, $n \geq 4$, the following identity 
is satisfied (see, e.g., {\cite[Section 2] {DGHSaw}}) 
\begin{eqnarray}
(n-2)(R \cdot C-C \cdot R) 
&=& Q(S,R) - \frac{\kappa}{n-1}\, Q(g,R) - g \wedge (R \cdot S) + P,
\label{1abcd}
\end{eqnarray}
where the $(0,6)$-tensor $P$ is defined by
\begin{eqnarray*}
P(X_1,X_2,X_3,X_4;X,Y)
&=&
g(X, X_1)\, R( {\mathcal{S}}Y, X_2,X_3,X_4)
-
g(Y, X_1)\, R( {\mathcal{S}}X, X_2,X_3,X_4)\nonumber\\
& &
+ g(X, X_2)R( X_{1}, {\mathcal{S}}Y, X_3,X_4)
-
g(Y, X_2)\, R( X_{1}, {\mathcal{S}}X, X_3,X_4)\nonumber\\
& &
+ g(X, X_3)\, R( X_{1}, X_{2}, {\mathcal{S}}Y, X_4)
-
g(Y, X_3)\, R( X_{1}, X_{2}, {\mathcal{S}}X, X_4)\nonumber\\
& &
+ g(X, X_4)\, R( X_{1}, X_{2}, X_3, {\mathcal{S}}Y)
-
g(Y, X_4)\, R( X_{1}, X_{2}, X_3, {\mathcal{S}}X ) ,
\end{eqnarray*}
and $X_1, X_2, X_3, X_4, X, Y$ are vector fields on $M$. 
The local expression of the identity (\ref{1abcd}) reads
\begin{eqnarray*}
(n-2)(R \cdot C - C \cdot R)_{hijklm}
&=& 
Q(S,R)_{hijklm} - \frac{\kappa}{n-1}\, Q(g,R)_{hijklm} \\
& &
+ g_{hl}A_{mijk} - g_{hm}A_{lijk} - g_{il}A_{mhjk} + g_{im}A_{lhjk}\\
& &
+ g_{jl}A_{mkhi} - g_{jm}A_{lkhi} - g_{kl}A_{mjhi} + g_{km}A_{ljhi} \\
& &
- g_{ij}(R\cdot S)_{hklm} - g_{hk}(R \cdot S)_{ijlm}
+ g_{ik}(R \cdot S)_{hjlm} + g_{hj} (R \cdot S)_{iklm} ,
\end{eqnarray*}
where $A_{mijk} = g^{rs}S_{mr}R_{sijk}$ and
\begin{eqnarray*}
(R \cdot S)_{hklm} 
&=& g^{rs} ( S_{hr}R_{sklm} + S_{kr}R_{shlm}) ,\\
Q(g,S)_{hklm} 
&=& 
g_{hl}S_{km} + g_{kl}S_{hm} - g_{hm}S_{kl} - g_{km}S_{hl} ,\\
\\
Q(g,R)_{hijklm} 
&=& 
g_{hl}R_{mijk} + g_{il}R_{hmjk} + g_{jl}R_{himk} + g_{kl}R_{hijm}\\
& &
- g_{hm}R_{lijk} - g_{im}R_{hljk} - g_{jm}R_{hilk} - g_{km}R_{hijl} ,\\
Q(S,R)_{hijklm} 
&=& 
S_{hl}R_{mijk} + S_{il}R_{hmjk} + S_{jl}R_{himk} + S_{kl}R_{hijm}\\
& &
- S_{hm}R_{lijk} - S_{im}R_{hljk} - S_{jm}R_{hilk} - S_{km}R_{hijl} ,\\
\\
(R \cdot C)_{hijklm} &=& g^{rs} ( C_{rijk} R_{shlm}
+ C_{hrjk} R_{silm} + C_{hirk} R_{sjlm} + C_{hijr} R_{sklm} ) ,\\
(C \cdot R)_{hijklm}
&=& g^{rs} ( R_{rijk} C_{shlm} + R_{hrjk} C_{silm}
+ R_{hirk} C_{sjlm} + R_{hijr} C_{sklm} ) .
\end{eqnarray*}

Let $A$ be a symmetric $(0,2)$-tensor on a semi-Riemannian manifold
$(M,g)$, $\dim M = n \geq 3$. 
Let $A_{ij}$ be the local components of the tensor $A$. 
Further, let $A^{2}$ and $A^{3}$ be the $(0,2)$-tensors 
with the local components  
$A^{2}_{ij} = A_{ir}g^{rs}A_{sj}$
and
$A^{3}_{ij} = A^{2}_{ir}g^{rs}A_{sj}$, respectively.
We have 
$\mathrm{tr}_{g} (A) = \mathrm{tr} (A) = g^{rs}A_{rs}$,
$\mathrm{tr}_{g} (A^{2}) = \mathrm{tr} (A^{2}) = g^{rs}A^{2}_{rs}$
and
$\mathrm{tr}_{g} (A^{3}) = \mathrm{tr} (A^{3}) = g^{rs}A^{3}_{rs}$.
We denote by ${\mathcal U}_{A}$ the set
of points of $M$ at which $A \neq \frac{\mathrm{tr} (A)}{n}\, g$.

\begin{prop} 
Let $A$ be a symmetric $(0,2)$-tensor on a semi-Riemannian manifold
$(M,g)$, $\dim M = n \geq 3$, such that 
$\mathrm{rank} A = 2$ on ${\mathcal U}_{A} \subset M$.
Then
on this set we have {\cite[Lemma 2.1; eqs. (2.6) and (2.10)] {P106}}
\begin{eqnarray}
A^{3} &=& \mathrm{tr} (A) \, A^{2} 
+ \frac{\mathrm{tr} (A^{2}) - (\mathrm{tr} (A))^{2} }{2}\, A
\label{2023.05.05.aa} \\
A \wedge A^{2}
&=&
\frac{\mathrm{tr} (A)}{2}\, A \wedge A .
\label{2023.05.05.bb}
\end{eqnarray}
Moreover, the following identity is satisfied on ${\mathcal U}_{A}$
\begin{eqnarray}
Q(A^{2},\frac{1}{2}\,  A \wedge A) \ =\
- Q(A, A \wedge A^{2} )\ =\ 
- \frac{\mathrm{tr} (A)}{2}\, Q(A, A \wedge A ) \ =\ 0 .
\label{2023.05.05.cc}
\end{eqnarray}
\end{prop}

\begin{prop} 
Let $A$ be a symmetric $(0,2)$-tensor on a semi-Riemannian manifold
$(M,g)$, $\dim M = n \geq 4$.
\newline
(i) ( {\cite[Lemma 2.1] {DGHSaw-2022}}, 
see also {\cite[Proposition 2.1 (i)] {2023_DGHP-TZ}}) 
If the following condition is satisfied on ${\mathcal U}_{A}  \subset M$
\begin{eqnarray*}
\mathrm{rank} ( A - \alpha \, g ) = 1 
%\label{2022.12.22.aa} 
\end{eqnarray*}
then
\begin{eqnarray}
g \wedge A^{2} + \frac{n-2}{2}\, A \wedge A 
- \mathrm{tr} (A) \, g \wedge A 
+ \frac{(\mathrm{tr} (A))^{2} - \mathrm{tr} (A^{2})}{2 (n-1)} \, g \wedge g = 0 
\label{2022.12.22.bb} 
\end{eqnarray}
and
\begin{eqnarray*}
A^{2} - \frac{ \mathrm{tr} (A^{2})}{n} =
( \mathrm{tr} (A) - (n - 2) \alpha ) 
\left( A - \frac{\mathrm{tr} (A) }{n} \, g \right)  
%\label{2022.12.20.kk} 
\end{eqnarray*}
on ${\mathcal U}_{A}$,
where $\alpha$ is some function on ${\mathcal U}_{A}$. 
\newline
(ii) {\cite[Proposition 2.1 (ii)] {2023_DGHP-TZ}} 
If (\ref{2022.12.22.bb}) is satisfied 
on ${\mathcal U}_{A} \subset M$ then 
\begin{eqnarray*}
A^{2} - \frac{\mathrm{tr} (A^{2})}{n}\, g 
= 
\rho \left( A - \frac{\mathrm{tr} (A)}{n}\, g \right) 
%\label{2022.12.22.ii} 
\end{eqnarray*}
and
\begin{eqnarray*}
\left(
A - \frac{ \mathrm{tr} (A) - \rho }{n-2}\, g
\right)
\wedge
\left(
A - \frac{ \mathrm{tr} (A) - \rho }{n-2}\, g
\right) = 0
%\label{2022.12.22.jj} 
\end{eqnarray*}
on ${\mathcal U}_{A}$,
where $\rho$ is some function on ${\mathcal U}_{A}$.
\end{prop}

Let $(M,g)$ be a semi-Riemannian manifold of dimension $\dim M = n \geq 3$.
We set
\begin{eqnarray}
E 
= g \wedge S^{2} + \frac{n-2}{2} \, S \wedge S - \kappa \, g \wedge S
+ \frac{\kappa ^{2} - \mathrm{tr} (S^{2})}{2(n-1)}\, g \wedge g.
\label{2022.11.10.aaa}
\end{eqnarray}
It is easy to check
that the tensor $E$ is a generalized curvature tensor.
In the same way, we define the $(0,4)$-tensor $E(A)$ corresponding 
to a symmetric  $(0,2)$-tensor $A$ {\cite[eq. (1.7] {2023_DGHP-TZ})
\begin{eqnarray}
E(A) = g \wedge A^{2} + \frac{n-2}{2} \, A \wedge A 
- \mathrm{tr}_{g} ( A)  \, g \wedge A
+ \frac{(\mathrm{tr} ( A) )^{2} - \mathrm{tr} (A^{2})}{2(n-1)}
\, g \wedge g .
\label{2022.12.33.aaa}
\end{eqnarray}

Let $T$ be a generalized curvature tensor on a semi-Riemannian manifold 
$(M,g)$, $\dim M = n \geq 4$. We denote by
$\mathrm{Ric} (T)$, $\kappa (T)$ and $\mathrm{Weyl} (T)$ the Ricci tensor, 
the scalar curvature and the Weyl tensor of the tensor $T$, respectively.
We refer to {\cite[Section 2] {DGHHY}},  {\cite[Section 3] {2021-DGH}} 
or {\cite[Section 3] {DGHZ01}} 
for definitions of the considered tensors.
In particular, we have 
\begin{eqnarray*}
\mathrm{Weyl} (T) = T - \frac{1}{n-2} \, g \wedge \mathrm{Ric} (T) 
+  \frac{\kappa (T) }{2 (n-2)(n-1)}\, g \wedge g .
%\label{2023.03.06a}
\end{eqnarray*}

Let $A$ be a symmetric $(0,2)$-tensor on a semi-Riemannian manifold
$(M,g)$, $\dim M = n \geq 3$. Let $E(A)$ be the tensor 
defined by (\ref{2022.12.33.aaa}). It is easy to check that 
$\mathrm{Ric} (E(A))$ is a zero tensor. Therefore, we also have
$\kappa (E(A)) = 0$. Any generalized curvature tensor $T$
defined on a $3$-dimensional semi-Riemannian manifold $(M,g)$ 
is expressed by
$T = g \wedge \mathrm{Ric} (T) - (\kappa (T)/4) g \wedge g$
(see {\cite[Section 2, p. 383] {2023_DGHP-TZ}} and references therein).
Thus we see that the tensor $T = E(A)$ on any $3$-dimensional 
semi-Riemannian manifold $(M,g)$ is a zero tensor. In particular, 
on any $3$-dimensional semi-Riemannian manifold $(M,g)$ we have
$E = 0$.

\begin{prop}
{\cite[Proposition 2.2] {2023_DGHP-TZ}}
Let $T$ be a generalized curvature tensor 
on a semi-Riemannian manifold $(M,g)$, 
$\dim M = n \geq 4$.
If the following condition is satisfied at a point $x \in M$ 
\begin{eqnarray*}
T = \alpha _{1} \, R + \frac{\alpha _{2}}{2} \, S \wedge S 
+ \alpha _{3}\, g \wedge S  +  \alpha _{4}\, g \wedge S^{2} 
+  \frac{ \alpha _{5} }{2} \, g \wedge g
%\label{2023.03.06b}
\end{eqnarray*}
then 
\begin{eqnarray*}
\mathrm{Weyl} (T) =  \alpha _{1} \, C +   \frac{\alpha _{2}}{n-2} \, E 
%\label{2023.03.06c}
\end{eqnarray*}
at this point, where the tensor $E$ is defined by (\ref{2022.11.10.aaa}) 
and $\alpha_{1}, \ldots , \alpha_{5} \in \mathbb{R}$.
\end{prop}

According to \cite{P106}, 
a generalized curvature tensor $T$ on a semi-Riemannian manifold 
$(M,g)$, $\dim M = n \geq 4$,
is called a {\sl{Roter type tensor}} if 
\begin{eqnarray}
T &=& \frac{\phi}{2}\, Ric(T) \wedge Ric(T) 
+ \mu\, g \wedge Ric(T) + \eta\, G 
\label{eq:h7}
\end{eqnarray}
on ${\mathcal{U}}_{Ric(T)} \cap {\mathcal{U}}_{Weyl (T)}$,
where $\phi$, $\mu $ and $\eta $ are some functions on this set.
Manifolds admitting 
Roter type tensors were investigated (e.g.) in \cite{Kow02}. 
We have 
\begin{prop}
Let $T$ be a generalized curvature tensor 
on a semi-Riemannian manifold $(M,g)$, $\dim M = n \geq 4$,
satisfying (\ref{eq:h7}) on 
${\mathcal{U}}_{Ric(T)} \cap {\mathcal{U}}_{Weyl (T)} \subset M$. 
Then the following relations hold on this set
\newline
(i) {\cite[Proposition 3.2] {DGHZ01}}
\begin{eqnarray*}
& &\ \ \ (Ric(T))^{2} 
\ =\ \alpha _{1}\, Ric(T) + \alpha _{2} \, g ,\nonumber\\
& &\ \ \
\alpha _{1} \ =\ \kappa(T) + \phi^{-1} ( (n-2)\mu -1 ),\ \ \ \alpha _{2}
\ =\
\phi^{-1} (  \mu \kappa(T) + (n-1) \eta ) ,
%\label{roter71}
\end{eqnarray*}
\begin{eqnarray*}
&(a)&\ \ \ T \cdot T \ =\ L_{T}\, Q(g,T),
\ \ \
L_{T}
\ =\
\phi^{-1}  \left( (n-2) (\mu ^{2} - \phi \eta) - \mu \right),\nonumber\\
&(b)&\ \ \ T \cdot Weyl(T) \ =\ L_{T}\, Q(g,Weyl(T)),\nonumber\\
&(c)&\ \ \ 
T \cdot T \ =\ Q(Ric(T),T) + (L_{T} + \phi^{-1}  \mu ) \, Q(g,Weyl(T)) ,
%\label{roter72}
\end{eqnarray*}
\begin{eqnarray}
Weyl(T) \cdot T &=& L_{Weyl(T)}\, Q(g,T) ,\ \ \
L_{Weyl(T)} 
\ =\ L_{T} + \frac{1}{n-2}  \left( \frac{\kappa(T) }{n-1} 
- \alpha _{1} \right) ,
\label{roter73NN}\\
%\end{eqnarray}
%\begin{eqnarray}
Weyl(T) \cdot Weyl(T) &=& L_{Weyl(T)}\, Q(g,Weyl(T)) ,
\label{roter73}\\
Weyl(T) \cdot T &=& 
Q(Ric(T), Weyl(T)) + \left( L_{T} - \frac{\kappa(T)}{n-1} \right)  
Q(g,Weyl(T)) .
%\label{roter73abc}
\nonumber
\end{eqnarray}
We also have
\begin{eqnarray*}
& &
T \cdot Weyl(T) - Weyl(T) \cdot T \ =\ 
\left( \frac{(n-1) \mu - 1}{(n-2) \phi} 
+ \frac{\kappa(T) }{n-1} \right) Q(g,T)\nonumber\\
& &
+ \frac{1}{n-2}\, Q(Ric(T),T) + \frac{\mu ((n-1) \mu - 1) 
- (n-1) \phi \eta }{(n-2) \phi} \, Q(Ric(T),G) ,
%\label{dh2}
\end{eqnarray*}
and, equivalently,
\begin{eqnarray*}
T \cdot Weyl(T) - Weyl(T) \cdot T &=& 
\left( \phi^{-1}  
\left( \mu - \frac{1}{n-2} \right) 
+ \frac{\kappa(T) }{n-1} \right) Q(g,T)\nonumber\\
& &+ \left(  \phi^{-1} \mu \left( \mu - \frac{1}{n-2} \right) 
- \eta \right) Q(Ric(T),G) .
%\label{dh3}
\end{eqnarray*}
(ii) {\cite[Sections 1 and 4] {Kow02}}
\begin{eqnarray*}
Q(Ric(T),Weyl(T)) &=&   
\phi ^{- 1}  \left( \frac{1}{n-2} - \mu \right) Q(g,T)\nonumber\\
& &
+ \frac{1}{n-2}\left(  L_{T} - \frac{\kappa(T) }{n-1} \right)  
Q(g, g \wedge Ric(T)) .
%\label{roter74}
\end{eqnarray*}
Moreover, if $L_{T} \, =\, \frac{ \kappa(T) }{n-1}$, 
resp., $\kappa(T) \, =\, 0$, then we have
\begin{eqnarray*}
Q(Ric(T),Weyl(T)) &=& L_{Weyl(T)}\, Q(g,T) ,
%\label{roter75}
\end{eqnarray*}
and
\begin{eqnarray*}
T \cdot Weyl(T) - Weyl(T) \cdot T &=& - Q(Ric(T),Weyl(T)) ,
%\label{dh3aa}
\end{eqnarray*}
respectively. Moreover, if
$L_{Weyl(T)}\, =\, 0$ at a point of
${\mathcal{U}}_{Ric(T)} \cap {\mathcal{U}}_{Weyl (T)}$ 
then at this point we have
\begin{eqnarray*}
& &
\phi ^{-1} \left( \frac{1}{n-2} - \mu \right) \,
( T \cdot Weyl(T) - Weyl(T) \cdot T)\nonumber\\
&=& 
(n-2)
\left( \phi^{-1} \mu  
\left( \mu - \frac{1}{n-2} \right) - \eta \right) Q(Ric(T), Weyl(T)) .
%\label{dh3abab}
\end{eqnarray*}
\end{prop}

For further results on generalizing curvature tensors satisfying 
some conditions we refer to {\cite[Section 3] {2021-DGH}}.

Let $A$ be a symmetric $(0,2)$-tensor and 
$T$ a $(0,k)$-tensor, $k = 2, 3, \ldots $.  
The tensor $Q(A,T)$ is called the {\sl{Tachibana tensor}} of $A$ and $T$, 
or the Tachibana tensor for short (see, e.g., \cite{DGPSS}). 
Using the tensors $g$, $R$ and $S$ we can define the following 
$(0,6)$-Tachibana tensors: 
$Q(S,R)$, $Q(g,R)$, $Q(g,g \wedge S)$ and $Q(S, g \wedge S)$. 
We can check, by making use of (\ref{DS7})(a) and (\ref{dghz01}), 
that other $(0,6)$-Tachibana tensors  
constructed from $g$, $R$ and $S$ may be expressed 
by the four Tachibana tensors mentioned above 
or vanish identically on $M$.

According to {\cite[eq. (1.13), Theorem 3.4 (i)] {2016_DGJZ}}
the following identity is satisfied on any semi-Riemannian manifold 
$(M,g)$ of dimension $n \geq 4$ 
\begin{eqnarray}
C \cdot R  + R \cdot C 
&=& R \cdot R + C \cdot C 
- \frac{1}{(n-2)^{2}}\, Q\left( g, g \wedge S^{2}  
- \frac{\kappa }{ n-1}\, g\wedge S \right) .
\label{identity01}
\end{eqnarray}
From (\ref{identity01}), by a suitable contraction, we get 
(cf.  {\cite[Lemma 2.3] {DGP-TV02}}, {\cite[p. 217] {DHS-1999}})
\begin{eqnarray*}
C \cdot S 
&=& R \cdot S - \frac{1}{n-2} 
Q\left( g, S^{2} - \frac{\kappa }{n-1} S \right) .
%\label{uuWeyl}
\end{eqnarray*}

Let  $(\overline{M},\overline{g})$ and $(\widetilde{N},\widetilde{g})$,
$\dim \overline{M} = p$, $\dim \widetilde{N} = n-p$, $1 \leq p < n$, 
be semi-Riemannian manifolds
covered by systems of charts $\{ U;x^{a} \}$ 
and 
$\{ V;y^{\alpha } \} $,
respectively.
Let $F$ be a positive smooth function on $\overline{M}$.
It is well known that the {\sl warped product}
$\overline{M} \times _F \widetilde{N}$ of $(\overline{M},\overline{g})$ 
and $(\widetilde{N}, \widetilde{g})$ 
is the product manifold $\overline{M} \times \widetilde{N}$ 
with the metric
$g = \overline{g} \times _F \widetilde{g} $ defined by 
(see, e.g., \cite{{Kru01}, {ON}}) 
\begin{eqnarray*}
\overline{g} \times _F \widetilde{g} & = &
{\pi}_1^{*} \overline{g} + (F \circ {\pi}_1)\, {\pi}_2^{*} \widetilde{g},
\end{eqnarray*}
where 
${\pi}_1 : 
\overline{M} \times \widetilde{N} \longrightarrow \overline{M}$ and 
${\pi}_2 : 
\overline{M} \times \widetilde{N} \longrightarrow \widetilde{N}$ 
are the natural projections on $\overline{M}$ 
and $\widetilde{N}$, respectively.
Let 
$ \{ U \times V ; x^{1}, \ldots ,x^{p},x^{p+1} =  y^{1}, \ldots , 
x^{n} = y^{n-p} \} $ 
be a product chart for $\overline{M} \times \widetilde{N}$. The local
components $g_{ij}$ of the metric 
$g = \overline{g} \times _F \widetilde{g}$ with respect
to this chart are the following
$g_{ij} = \overline{g}_{ab}$ if $i = a$ and $j = b$,
$g_{ij} = F\, \widetilde{g}_{\alpha \beta }$ 
if $i = \alpha $ and $j = \beta $, and
$g_{ij} = 0$ otherwise,
where $a,b,c, d, f \in \{ 1, \ldots ,p \}$,
$\alpha , \beta , \gamma , \delta \in \{ p+1, \ldots ,n \}$
and $h, i, j, k, r, s \in \{ 1,2, \ldots ,n \}$.
We will denote by bars (resp., by tildes) tensors formed from 
$\overline{g}$ (resp., $\widetilde{g}$).
The local components 
\begin{eqnarray*}
\Gamma ^{h} _{ij} 
\ =\ \frac{1}{2}\, g^{hs} ( \partial_{i} g_{js} + \partial_{j} g_{is} 
- \partial_{s} g_{ij}),
\ \ 
\partial _j \ =\ \frac{\partial }{\partial x^{j}} ,
\end{eqnarray*}
of the Levi-Civita connection $\nabla $
of $\overline{M} \times _F \widetilde{N}$ are the following 
(see, e.g., \cite{Kru01}):
\begin{eqnarray*}
& &
\Gamma ^{a} _{bc}\ =\ \overline{\Gamma } ^{a} _{bc} , \ \
\Gamma ^{\alpha } _{\beta \gamma } 
\ =\ \widetilde{\Gamma } ^{\alpha } _{\beta \gamma } ,  \ \
\Gamma ^{a} _{\alpha \beta } 
\ =\ - \frac{1}{2} \bar{g} ^{ab} F_b \widetilde{g} _{\alpha \beta } , \ \
\Gamma ^{\alpha } _{a \beta } 
\ =\ \frac{1}{2F} F_a \delta ^{\alpha } _{\beta } ,\\
%\end{eqnarray*}
%\begin{eqnarray*}
& &
\Gamma ^{a} _{\alpha b}\ =\ \Gamma ^{\alpha } _{ab}\ =\ 0 , \ \
F_a\ =\ \partial _a F \ =\ \frac{\partial F}{\partial x^{a}} , \ \
\partial _a \ =\ \frac{\partial }{\partial x^{a}} .
\end{eqnarray*}
The local components
\begin{eqnarray*}
R_{hijk}\ =\ g_{hs}R^{s}_{\, ijk}
\ =\ g_{hs} (\partial _k \Gamma ^{s} _{ij} 
- \partial _j \Gamma ^{s} _{ik} + \Gamma ^{r} _{ij} \Gamma ^{s} _{rk}
- \Gamma ^{r} _{ik} \Gamma ^{s} _{rj} ) ,\ \ 
\partial _k \ =\ \frac{\partial }{\partial x^{k}} ,
\end{eqnarray*}
of the Riemann-Christoffel curvature tensor $R$
and the local components $S_{ij}$ of the Ricci tensor $S$
of the warped product $\overline{M} \times _F N$ which may not vanish
identically are the following:
\begin{eqnarray*}
& &
R_{abcd} \ =\ \overline{R}_{abcd} ,\ \ \ 
R_{\alpha ab \beta} \ =\
- \frac{1}{2}\, T_{ab} \widetilde{g}_{\alpha \beta} ,\ \ \ 
R_{\alpha \beta \gamma \delta}  \ =\ 
F \widetilde{R}_{\alpha \beta \gamma \delta} 
- \frac{1}{4}\, \Delta_1 F\, 
\widetilde{G}_{\alpha \beta \gamma \delta}\, ,\\
%\label{L3}
%\end{eqnarray*}
%\begin{eqnarray*}
& &
S_{ab} \ =\ \overline{S}_{ab} 
- \frac{n-p}{2}\, \frac{1}{F}\, T_{ab} ,\ \ \  
S_{\alpha \beta } \ =\ 
\tilde{S}_{\alpha \beta } 
- \frac{1}{2}\, ( tr(T)
+ \frac{n-p-1}{2F} \Delta _1 F ) \widetilde{g}_{\alpha \beta } ,
%\label{AL2}
\end{eqnarray*}
where
\begin{eqnarray*}
T_{ab}&=& 
\overline{\nabla }_b F_a - \frac{1}{2F} F_a F_b,\ \
tr(T) \ =\ \overline{g}^{ab} T_{ab} 
\ =\ \Delta F - \frac{1}{2 F}\, {\Delta}_1 F ,\ \ 
\nonumber\\ 
\Delta F &=& {\Delta}_{\overline{g}} F\ =\ 
\overline{g}^{ab} \nabla _a F_b,\ \ 
{\Delta}_1 F \ = \ {\Delta}_{1 \overline{g}} F
\ =\ \overline{g}^{ab} F_a F_b\, ,
%\label{AL3}
\end{eqnarray*}
and $T$ is the $(0,2)$-tensor with the local components $T_{ab}$.
We can express the scalar curvature $\kappa $ of
$\overline{M} \times _F \widetilde{N}$ by 
\begin{eqnarray*}
\kappa &=& \overline{\kappa } + \frac{1}{F}\, \widetilde{\kappa }
- \frac{n - p}{F}\, ( tr(T) + \frac{n - p -1 }{4F} \Delta _1 F)
\ =\ 
\overline{\kappa } + \frac{1}{F}\, \widetilde{\kappa }
- \frac{n - p}{F}\, ( \Delta F + \frac{n - p - 3 }{4F} \Delta _1 F)  .
%\label{AL4}
\end{eqnarray*}

\section{Pseudosymmetry type curvature conditions}

It is well-known that if a semi-Riemannian manifold $(M,g)$, 
$\dim M = n \geq 3$, 
is locally symmetric then $\nabla R = 0$
on $M$
(see, e.g., {\cite[Chapter 1.5] {Lumiste}}). 
This  
implies  the following integrability condition 
${\mathcal{R}}(X,Y ) \cdot R = 0$,
in short 
\begin{eqnarray}
R \cdot R &=& 0 .
\label{semisymmetry}
\end{eqnarray}
Semi-Riemannian manifold satisfying (\ref{semisymmetry})
is called {\sl semisymmetric} \cite{Sz 1} (see also 
{\cite[Chapter 8.5.3] {TEC_PJR_2015}}, {\cite[Chapter 20.7] {Chen-2011}}, 
{\cite[Chapter 1.6] {Lumiste}}, \cite{{Sz 2}, {Sz 3}, {LV3-Foreword}}).

Semisymmetric manifolds form a subclass of the class of 
pseudosymmetric manifolds.
A semi-Rieman\-nian manifold $(M,g)$, $\dim M = n \geq 3$,
 is said to be {\sl pseudosymmetric} 
if the tensors $R \cdot R$ and $Q(g,R)$ 
are linearly dependent at every point of $M$
(see, e.g., 
{\cite[Chapter 8.5.3] {TEC_PJR_2015}}, 
{\cite[Chapter 20.7] {Chen-2011}},
{\cite[Section 15.1] {Chen-2021}},
{\cite[Chapter 6] {DHV2008}},
{\cite[Chapter 12.4] {Lumiste}}, 
\cite{{D-1992}, 
{DGHHY}, {DGHSaw}, {DHV2008}, {DVV1991}, 
{HV_2007}, {HaVerSigma}, 
{SDHJK}, {LV1}, {LV2}, 
{LV3-Foreword}, {LV4}} and references therein). 
This is equivalent to
\begin{eqnarray}
R \cdot R &=& L_{R}\, Q(g,R) 
\label{pseudo}
\end{eqnarray}
on   ${\mathcal{U}}_{R}  \subset M$,
where $L_{R}$ is some function on ${\mathcal{U}}_{R}$. 
Every semisymmetric manifold is pseudosymmetric.
The converse statement is not true (see, e.g., \cite{{DVV1991}}).
A non-semisymmetric pseudosymmetric manifold is called 
{\sl properly pseudosymmetric}.
A pseudosymmetric manifold $(M,g)$, $\dim M = n \geq 3$, 
is called a {\sl pseudo\-sym\-metric manifold of constant type} 
if the function $L_{R}$ is constant on ${\mathcal{U}}_{R}$
\cite{KowSek_1997} 
(see also \cite{{HajKowSek}, {HaSek}, {OS0}, {OS2}, {OS3}, {KoSe5}}).

A semi-Riemannian manifold $(M,g)$, $\dim M = n \geq 3$, 
is called {\sl Ricci-pseudosymmetric} 
if the tensors $R \cdot S$ and $Q(g,S)$ 
are linearly dependent at every point of $M$
(see, e.g., {\cite[Chapter 8.5.3] {TEC_PJR_2015}}, 
{\cite[Section 15.1] {Chen-2021}}, \cite{DGHSaw}).
This is equivalent on ${\mathcal{U}}_{S} \subset M$ to 
\begin{eqnarray}
R \cdot S &=& L_{S}\, Q(g,S) , 
\label{Riccipseudo07}
\end{eqnarray}
where $L_{S}$ is some function on ${\mathcal{U}}_{S}$. 
Every warped product manifold $\overline{M} \times _{F} \widetilde{N}$
with a $1$-dimensional $(\overline{M}, \overline{g})$ manifold and
an $(n-1)$-dimensional Einstein semi-Riemannian manifold 
$(\widetilde{N}, \widetilde{g})$, $n \geq 3$, 
and a warping function $F$, 
is a Ricci-pseudosymmetric manifold,
see, e.g., 
{\cite[Section 1] {Ch-DDGP}}
and
{\cite[Example 4.1] {2016_DGJZ}}.
According to \cite{G6}, 
a Ricci-pseudosymmetric manifold
$(M,g)$, $\dim M = n \geq 3$, 
is called a {\sl Ricci-pseudo\-sym\-metric manifold of constant type} 
if the function $L_{S}$ is constant on ${\mathcal{U}}_{S}$.

A semi-Riemannian manifold $(M,g)$, $\dim M = n \geq 4$, 
is said to be {\sl Weyl-pseudo\-sym\-met\-ric} 
if the tensors $R \cdot C$ and $Q(g,C)$ are linearly dependent 
at every point of $M$
\cite{{DGHHY}, {DGHSaw}}. 
This is equivalent on ${\mathcal{U}}_{C} \subset M$ to 
\begin{eqnarray}
R \cdot C &=& L_{1}\, Q(g,C) ,  
\label{Weyl-pseudo}
\end{eqnarray}
where $L_{1}$ is some function on ${\mathcal{U}}_{C}$. 
We can easily check that 
on every Einstein manifold $(M,g)$, $\dim M \geq 4$,
(\ref{Weyl-pseudo}) turns into
\begin{eqnarray*}
R \cdot R &=& L_{1}\, Q(g,R) .
\end{eqnarray*}
For a presentation of results on the problem 
of the equivalence of pseudosymmetry, 
Ricci-pseudo\-sym\-met\-ry and Weyl-pseudosymmetry
we refer to {\cite[Section 4] {DGHSaw}}.  
Inclusions between mentioned above semi-Riemannian 
manifolds $(M,g)$, i.e., pseudosymmetric, 
Ricci-pseudo\-sym\-met\-ric and Weyl-pseudosymmetric manifolds, 
can be presented in the following diagram
{\cite[Section 4] {DGHSaw}}

\[
\newlength{\BW}\settowidth{\BW}{$\displaystyle 
R \cdot R = L_{R} Q(g,R)$}\addtolength{\BW}{0.6em}
\newlength{\BH}\setlength{\BH}{0.9cm}
\def\R{\rule[-0.4\BH]{0.0cm}{\BH}}
%% albo ramki prostokatne
\def\FR#1{\fboxsep0pt\fbox{\makebox[\BW]{\R $\displaystyle #1$\R}}}
%% albo z zaokraglonymi rogami
%\def\FR#1{\fboxsep0pt\ovalbox{\makebox[\BW]{\R $\displaystyle #1$\R}}}
\def\Upset{\rotatebox{90}{$\displaystyle{}\;\subset\;\;\;{}$}}
\def\Eq{ = }
\begin{array}{ccccc}
\FR{R \cdot S \Eq L_{S} Q(g,S)} & \supset & 
\FR{R \cdot R \Eq L_{R} Q(g,R)} & \subset & 
\FR{R \cdot C \Eq L_{C} Q(g,C)} \\

\Upset & & \Upset & & \Upset \\

\FR{R \cdot S \Eq 0} & \supset & 
\FR{R \cdot R \Eq 0} & \subset & 
\FR{R \cdot C \Eq 0}\\

\Upset & & \Upset & & \Upset \\

\FR{\nabla S \Eq 0} & \supset & 
\FR{\nabla R \Eq 0} & \subset & 
\FR{\nabla C \Eq 0}\\

\Upset & & \Upset & & \Upset \\

\FR{S \Eq \frac{\kappa }{n} g} & \supset & 
\FR{R \Eq \frac{\kappa }{(n-1) n} G} & \subset & 
\FR{C \Eq 0}
\end{array}\]

\vspace{3mm}

\noindent
All inclusions in the above presentation are strict, provided that 
$\dim M = n \geq 4$.

A semi-Riemannian manifold $(M,g)$, $\dim M = n \geq 4$, is said to be 
a {\sl manifold with pseudosymmetric Weyl tensor}
({\sl to have a pseudosymmetric conformal Weyl tensor})
if the tensors $C \cdot C$ and $Q(g,C)$ 
are linearly dependent at every point of $M$ 
(see, e.g., {\cite[Section 15.1] {Chen-2021}}, 
\cite{{DGHHY}, {DGHSaw}, {2016_DGJZ}}).
This is equivalent on ${\mathcal U}_{C} \subset M$ to 
\begin{eqnarray}
C \cdot C &=& L_{C}\, Q(g,C) ,  
\label{4.3.012}
\end{eqnarray}
where $L_{C}$ is some function on ${\mathcal{U}}_{C}$. 
Every warped product manifold $\overline{M} \times_{F} \widetilde{N}$, 
with $\dim\, \overline{M} = \dim \widetilde{N} = 2$, satisfies (\ref{4.3.012})
(see, e.g., \cite{{DGHHY}, {DGHSaw}, {2016_DGJZ}} and references therein).
Thus in particular,
the Schwarzschild spacetime, the Kottler spacetime
and the Reissner-Nordstr\"{o}m spacetime satisfy (\ref{4.3.012}).
Semi-Riemannian manifolds with pseudosymmetric Weyl tensor 
were investigated among others in \cite{{DGHHY}, {DeHoJJKunSh}, {43}, {DY-1994}}.

We can show that 
$C \cdot R  = L\, Q(g,R)$ implies $C \cdot C = L\, Q(g,C)$
(see {\cite[Proposition 2.1] {MADEO}}).
Let $(M,g)$, $\dim M = n \geq 4$, be a semi-Riemannian manifold satisfying 
$C \cdot R  = L\, Q(g,R)$ on ${\mathcal U}_{C} \subset M$.
From this we get easily 
$C \cdot S = L\, Q(g,S)$ on ${\mathcal U}_{C}$.
Further, we have 
\begin{eqnarray*}
C \cdot C 
&=& C \cdot \left( R 
- \frac{1}{n-2}\, g \wedge S + \frac{\kappa}{(n-2)(n-1)}\, G \right) \\
&=& C \cdot R - \frac{1}{n-2}\, g \wedge ( C \cdot S) 
+ \frac{\kappa}{(n-2)(n-1)}\, C \cdot G \\
&=& L\, Q(g,R) - \frac{L}{n-2}\, g \wedge Q(g,S) \\
&=& L\, Q(g,R) - \frac{L}{n-2}\, Q( g,   g \wedge S  )
\ =\ L\, Q(g, R -  \frac{1}{n-2}\, g \wedge S )\\
&=& L\, Q(g, R -  \frac{1}{n-2}\, g \wedge S 
+ \frac{\kappa}{(n-2)(n-1)}\, G  ) \ =\
L\, Q(g,C) .
\end{eqnarray*}

We also have the following diagram

\[
\newlength{\BWA}\settowidth{\BW}{$\displaystyle 
R \cdot R = L_{R} Q(g,R)$}\addtolength{\BW}{0.6em}
\newlength{\BHA}\setlength{\BH}{0.9cm}
\def\R{\rule[-0.4\BH]{0.0cm}{\BH}}
%% albo ramki prostokatne
\def\FR#1{\fboxsep0pt\fbox{\makebox[\BW]{\R $\displaystyle #1$\R}}}
%% albo z zaokraglonymi rogami
%\def\FR#1{\fboxsep0pt\ovalbox{\makebox[\BW]{\R $\displaystyle #1$\R}}}
\def\Upset{\rotatebox{90}{$\displaystyle{}\;\subset\;\;\;{}$}}
\def\Eq{ = }
\begin{array}{ccccc}
\FR{C \cdot S \Eq L_{S} Q(g,S)} & \supset & 
\FR{C \cdot R \Eq L_{R} Q(g,R)} & \subset & 
\FR{C \cdot C \Eq L_{C} Q(g,C)} \\

\Upset & & \Upset & & \Upset \\

\FR{C \cdot S \Eq 0} & \supset & 
\FR{C \cdot R \Eq 0} & \subset & 
\FR{C \cdot C \Eq 0}\\

\Upset & & \Upset & & \Upset \\

\FR{S \Eq \frac{\kappa }{n} g} & \supset & 
\FR{R \Eq \frac{\kappa }{(n-1) n} G} & \subset & 
\FR{C \Eq 0}
\end{array}\]

\vspace{3mm}

\noindent
All inclusions in the above presentation are strict, provided that 
$\dim M = n \geq 4$.

Warped product manifolds $\overline{M} \times _{F} \widetilde{N}$, 
of dimension $\geq 4$,
satisfying on 
${\mathcal U}_{C} \subset \overline{M} \times _{F} \widetilde{N}$,
the condition 
\begin{eqnarray}
R \cdot R - Q(S,R) &=& L\, Q(g,C) ,  
\label{genpseudo01}
\end{eqnarray}
where $L$ is some function on ${\mathcal{U}}_{C}$,
were studied among others in \cite{49}. In that paper
necessary and sufficient conditions for  
$\overline{M} \times _{F} \widetilde{N}$ 
to be a manifold satisfying (\ref{genpseudo01}) are given.
Moreover, in that paper it was proved that 
any $4$-dimensional warped product manifold 
$\overline{M} \times _{F} \widetilde{N}$, 
with a $1$-dimensional base $(\overline{M},\overline{g})$, 
satisfies (\ref{genpseudo01}) {\cite[Theorem 4.1] {49}}.

We refer to
\cite{{Ch-DDGP}, {DGHHY}, {DGHSaw}, {DGHZ01}, {2016_DGJZ}, 
{DHV2008}, {DeHoJJKunSh}, {SDHJK}} 
for details on semi-Riemannian manifolds satisfying 
(\ref{pseudo}) and (\ref{Riccipseudo07})-(\ref{genpseudo01}), 
as well other conditions of this kind, named pseudosymmetry 
type curvature conditions. 
We also refer to {\cite[Section 3] {DeHoJJKunSh}} 
for a recent survey on manifolds 
satisfying such curvature conditions.
It seems that the condition (\ref{pseudo}) 
is the most important condition of that family of curvature conditions
(see, e.g., \cite{2016_DGJZ}).
The Schwarzschild spacetime, the Kottler spacetime, 
the Reissner-Nordstr\"{o}m spacetime, 
as well as the Friedmann-Lema{\^{\i}}tre-Robertson-Walker spacetimes 
are the ``oldest'' examples 
of pseudosymmetric warped product manifolds 
(see, e.g., \cite{{2016_DGJZ}, {DHV2008}, {DVV1991}, {SDHJK}}).

\section{Quasi-Einstein, 2-quasi-Einstein and partially Einstein manifolds}

A semi-Riemannian manifold $(M,g)$, 
$\dim M = n \geq 3$, 
is said to be a {\sl quasi-Einstein manifold} if 
\begin{eqnarray}
\mathrm{rank}\, (S - \alpha\, g) &=& 1
\label{quasi02}
\end{eqnarray}
on ${\mathcal U}_{S} \subset M$, where $\alpha $ is some function 
on ${\mathcal U}_{S}$.
It is known that every non-Einstein warped product manifold 
$\overline{M} \times _{F} \widetilde{N}$
with a $1$-dimensional $(\overline{M}, \overline{g})$ base manifold and
a $2$-dimensional manifold $(\widetilde{N}, \widetilde{g})$
or an $(n-1)$-dimensional Einstein manifold
$(\widetilde{N}, \widetilde{g})$
and a warping function $F$, $n \geq 4$, 
is a  quasi-Einstein manifold (see, e.g., \cite{{Ch-DDGP}, {2016_DGJZ}}). 
A Riemannian manifold 
$(M,g)$, $\dim M = n \geq 3$, whose Ricci tensor has an eigenvalue
of multiplicity $n-1$ is a quasi-Einstein manifold 
(cf. {\cite[Introduction] {P47}}). Evidently, 
the converse statement is also true.
We mention that quasi-Einstein manifolds arose during the study 
of exact solutions
of the Einstein field equations and the investigation 
on quasi-umbilical hypersurfaces 
of conformally flat spaces 
(see, e.g., \cite{{DGHSaw}, {2016_DGJZ}} and references therein). 
Quasi-Einstein hypersurfaces 
in semi-Riemannian spaces of constant curvature
were studied among others in
\cite{{DGHS}, {R102}, {P104}, {G6}} (see also \cite{DGHSaw} 
and references therein).
Quasi-Einstein manifolds satisfying some pseudosymmetry 
type curvature conditions  
were investigated recently in
\cite{{P119}, {Ch-DDGP}, {DGHHY}, {DGHZ01}, {DeHoJJKunSh}}. 
We mention that an example 
of a non-semisymmetric ($R \cdot R \neq 0$)
Ricci-semisymmetric ($R \cdot S = 0$)
quasi-Einstein hypersurface
$M$ in an Euclidean space $\mathbb{E}^{n+1}$, $n\geq 5$,
was constructed in \cite{AbDi}.

A semi-Riemannian manifold $(M,g)$, $\dim M = n \geq 3$, 
is called a $2$-{\sl{quasi-Einstein manifold}} if 
\begin{eqnarray}
\mathrm{rank}\, (S - \alpha \, g ) &\leq & 2 
\label{quasi0202weak}
\end{eqnarray}
on ${\mathcal U}_{S} \subset M$ and $\mathrm{rank}\, (S - \alpha \, g ) = 2$
on some open non-empty subset of ${\mathcal U}_{S}$, 
where $\alpha $ is some function on ${\mathcal U}_{S}$ 
(see, e.g., \cite{{DGP-TV02}}).

Every non-Einstein and non-quasi-Einstein warped product manifold 
$\overline{M} \times _{F} \widetilde{N}$
with a $2$-dimensional base manifold $(\overline{M}, \overline{g})$ 
and a $2$-dimensional manifold $(\widetilde{N}, \widetilde{g})$
or an $(n-2)$-dimensional Einstein semi-Riemannian manifold
$(\widetilde{N}, \widetilde{g})$, when $n \geq 5$, 
and a warping function $F$ satisfies (\ref{quasi0202weak})
(see, e.g., {\cite[Theorem 6.1] {2016_DGJZ}}). 
Thus some exact solutions of the Einstein field equations are 
non-conformally flat $2$-quasi-Einstein manifolds.
For instance, the Reissner-Nordstr\"{o}m spacetime, as well as
the Reissner-Nordstr\"{o}m-de Sitter type spacetimes 
are such manifolds (see, e.g.,
\cite{Kow02}). It seems that the Reissner-Nordstr\"{o}m spacetime
is the "oldest" example of 
a non-conformally flat 
$2$-quasi-Einstein warped product manifold.
It is easy to see that every $2$-quasi-umbilical hypersurface
in a semi-Riemannian space of constant curvature 
is a $2$-quasi-Einstein manifold
(see, e.g., \cite{DGP-TV02}).

We mention that  Einstein warped product manifolds $\overline{M} \times _{F} \widetilde{N}$,
with $\dim \overline{M} = 1, 2$, were studied in {\cite[Chapter 9.J] {Besse-1987}}
(see also {\cite[Chapter 3.F] {Besse-1987}}).

The semi-Riemannian manifold $(M,g)$, $\dim M = n \geq 3$, will be called
a {\sl{partially Einstein manifold}}, or
a {\sl{partially Einstein space}}
(cf. {\cite[Foreword] {CHEN-2017}}, {\cite[p. 20] {V2}}),  
if at every point $x \in {\mathcal{U}}_{S} \subset M$ 
its Ricci operator ${\mathcal{S}}$ satisfies
${\mathcal{S}}^{2} = \lambda {\mathcal{S}} + \mu  I\!d_{x}$,
or equivalently, 
\begin{eqnarray}
S^{2} &=&  \lambda \, S + \mu \, g ,
\label{partiallyEinstein}
\end{eqnarray}
where $\lambda , \mu \in {\mathbb{R}}$  
and $I\!d_{x}$ is the identity transformation of $T_{x} M$.
Thus every quasi-Einstein manifold is a partially Einstein manifold.
The converse statement is not true. 
Contracting (\ref{partiallyEinstein}) we get
$ \mathrm{tr} (S^{2}) = \lambda \, \kappa + n\, \mu$. 
This together with 
(\ref{partiallyEinstein}) yields (cf. {\cite[Section 5] {2021-DGH}})
\begin{eqnarray} 
S^{2} - \frac{ \mathrm{tr} (S^{2})}{n} \, g 
&=&  \lambda \left( S - \frac{\kappa }{n} \, g \right) .
\label{partiallyEinstein.11} 
\end{eqnarray}
In particular, a Riemannian manifold $(M,g)$, $\dim M = n \geq 3$, 
is a partially Einstein space if at every point 
$x \in {\mathcal{U}}_{S} \subset M$ 
its Ricci operator ${\mathcal{S}}$
has exactly two distinct eigenvalues
(principal Ricci curvatures)
$\rho _{1} = \rho _{2} = \ldots = \rho _{p}$ 
and 
$\rho _{p+1} = \rho _{p+2} = \ldots = \rho _{n}$ 
with multiplicities $p$ and $n-p$, respectively, where
$1 \leq p \leq n-1$.
Now (\ref{partiallyEinstein}) and (\ref{partiallyEinstein.11})  yield
\begin{eqnarray}
S^{2} \ =\  (\rho_{1} + \rho_{p+1})\, S - \rho_{1} \rho_{p+1} \, g 
\ \ \ \mbox{and} \  \  \
S^{2} - \frac{ \mathrm{tr} (S^{2})}{n} \, g 
\ =\  (\rho_{1} + \rho_{p+1}) \left( S - \frac{\kappa }{n} \, g \right) ,
\label{partiallyEinstein.117} 
\end{eqnarray}
respectively. 
Evidently, if $p = 1$, or $p = n-1$, 
then $(M,g)$ is a quasi-Einstein manifold.

\begin{rem}
(i)
Let $(M,g)$, $n = \dim M \geq 3$, be a semi-Riemannian manifold.
In addition, let $(M,g)$ be a conformally flat manifold when $n \geq 4$.
Thus the set ${\mathcal{U}}_{C} \subset M$ is an empty set.
Now, by (\ref{dghhy}), 
the subsets ${\mathcal{U}}_{S}$ and ${\mathcal{U}}_{R}$ of $M$ satisfy
${\mathcal{U}}_{S} = {\mathcal{U}}_{R}$.
In view of {\cite[Lemma 1.2] {DDV-1989}}
(see also  {\cite[Lemma 2.1] {DK-1999}}) 
we can state that on ${\mathcal{U}}_{S}$ any of the following
three conditions is equivalent to each other: 
$R \cdot R = \rho \, Q(g,R)$, $R \cdot S = \rho \, Q(g,S)$ and
\begin{eqnarray} 
S^{2} - \frac{ \mathrm{tr} (S^{2})}{n} \, g 
&=& \left( \frac{\kappa}{n-1} + (n-2) \rho \right) 
\left( S - \frac{\kappa }{n} \, g \right) .
\label{partiallyEinstein.1133} 
\end{eqnarray}
where $\rho$ is some function on ${\mathcal{U}}_{S}$.
\newline
(ii)  
Let $(M,g)$, $n = \dim M \geq 3$, be a Riemannian manifold with 
vanishing Weyl conformal curvature tensor $C$ 
such that its Ricci operator $\mathcal{S}$ 
has at every point of $M$ exactly two Ricci principal curvatures
$\rho _{1}$ and $\rho _{p+1}$ with multiplicities $p$ and $n-p$, 
respectively, where $1 \leq p < n$. 
Using 
(\ref{partiallyEinstein.117}) and
(\ref{partiallyEinstein.1133}) we can easily check 
that $(M,g)$ is a pseudosymmetric manifold satisfying (\ref{pseudo}) with 
\begin{eqnarray} 
L_{R} \ =\  \rho 
&=&  \frac{1}{(n-2)(n-1)} ( (n-1) (\rho_{1} + \rho_{p+1}) - \kappa )  .
\label{partiallyEinstein.1166} 
\end{eqnarray}
Evidently, if $p = 1$ then (\ref{partiallyEinstein.1166}) turns into 
\begin{eqnarray*} 
L_{R} \ =\  \rho 
&=&
 \frac{1}{n-1}  \rho_{1}   .
%\label{partiallyEinstein.116688} 
\end{eqnarray*}
(iii) 
As an immediate consequence of the above presented results 
we obtain {\cite[Foreword, p. xviii] {LV3-Foreword}}:
\newline
\sl{Riemannian spaces of dimensions $\geq 3$ with wanishing Weyl 
conformal curvature tensor are Deszcz symmetric if and only if they are
Einstein or "partially Einstein", - partially Einstein spaces being 
defined by the condition that their Ricci tensor has precisely 
two distinct eigenvalues -}.
\newline
(iv)
(cf. {\cite[Proposition 2.1] {G6}})
Let $(M,g)$, $\dim M = n \geq 3$, be a semi-Rieman\-nian manifold.
We note that (\ref{quasi02}) 
holds at a point $x \in {\mathcal U}_{S} \subset M$ 
if and only if 
$(S - \alpha\, g) \wedge (S - \alpha\, g) = 0$
at $x$, i.e.,  
\begin{eqnarray}
\frac{1}{2}\, S \wedge S   - \alpha\, g \wedge S 
+ \frac{\alpha ^{2}}{2}\, g \wedge g &=& 0 .   
\label{quasi03}
\end{eqnarray}
From (\ref{quasi03}), by a suitable contraction,
we get immediately 
\begin{eqnarray}
S^{2} 
&=& (\kappa - (n-2) \alpha )\, S + \alpha ( (n-1) \alpha - \kappa )\, g.
\label{quasi03quasi03} 
\end{eqnarray}
Evidently, 
(\ref{quasi03quasi03})
is a special case of (\ref{partiallyEinstein}).
\newline
(v) Let $\rho_{1}$, $\rho_{2}$ and $\rho_{3}$ be principal Ricci curvatures 
of a $3$-dimensional Riemannian manifold $(M,g)$. Moreover, let the condition 
\begin{eqnarray}
\rho_{1} = \rho_{2} \neq \rho_{3}
\label{3quasi0604a} 
\end{eqnarray} 
be satisfied at every point of the set 
${\mathcal{U}}_{S} = {\mathcal{U}}_{R} \subset M$.
Thus
\begin{eqnarray*}
R \cdot R &=& \frac{\rho_{3}}{2} \, Q(g,R)
\end{eqnarray*}
on ${\mathcal{U}}_{R}$.
Evidently, if  $\rho_{3} = 0$ on  ${\mathcal{U}}_{R}$ then 
$(M,g)$ is semisymmetric.
\newline
(vi) $3$-dimensional Riemannian manifolds satisfying (\ref{3quasi0604a})
were studied, among others, in the following papers:
\cite{{K4}, {KowSek_1996bb}, {KowSek_1996}, {KowSek_1997}, {OS0}, 
{OS2}, {OS3}}.  
We mention that an explicit classification of $3$-dimensional 
semisymmetric Riemannian manifolds 
is given in \cite{K_1996}. 
\end{rem}

\begin{example}
The warped product manifold 
$\overline{M} \times _{F} \widetilde{N}$ 
with a $1$-dimensional manifold 
$(\overline{M},\overline{g})$, $\overline{g}_{11} = \pm 1$,
and an $(n-1)$-dimensional semi-Riemannian Einstein manifold
$(\widetilde{N}, \widetilde{g})$, $n \geq 5$, 
assumed that it is not of constant curvature, 
and a warping function $F$,
satisfies on 
${\mathcal U}_{S} \cap {\mathcal U}_{C} 
\subset \overline{M} \times _{F} \widetilde{N}$
{\cite[Theorem 4.1] {Ch-DDGP}}
\begin{eqnarray*}
& &
\mathrm{rank} (S - \alpha \, g ) \ =\ 1,\ \ \
R \cdot S \ =\ L_{S}\, Q(g,S),\nonumber\\ 
& &
(n-2)\, (R \cdot C - C \cdot R) \ =\ Q(S,R) - L_{S}\, Q(g,R),
%\label{quasi10}
\\
& &
\alpha \ =\ \frac{\kappa }{n-1} - L_{S} ,\ \ \
L_{S} \ =\ - \frac{\mathrm{tr} T}{2 F} .\nonumber
\end{eqnarray*}
Furthermore, using
\begin{eqnarray*}
Q(g,R) &=& Q(g,C) - \frac{1}{n-2}\, Q(S,G) ,\\
Q(S,R) &=& Q(S,C) + \frac{1}{n-2} 
\left( \alpha - \frac{\kappa}{n-1} \right) Q(S,G) ,
\end{eqnarray*}
we obtain {\cite[Example 4.1] {2016_DGJZ}}
\begin{eqnarray*}
(n-2)\, (R \cdot C - C \cdot R) &=& Q(S,C) - L_{S}\, Q(g,C) .
\end{eqnarray*}
\end{example}

\begin{example} 
Let $M$ be an open connected non-empty subset of $\mathbb{R}^5$ endowed with 
the metric $g$ of the form \cite{{K1}, {K2}, {K3}} 
\begin{eqnarray*}
ds^{2} 
&=& g_{ij} dx^{i} dx^{j} 
\ =\ dx^{2} + dy^{2} + du^{2} + dv^{2} 
+ \rho ^{2}\, (x du - y dv + dz )^{2}\, ,\ \ \ 
\rho \ =\ \mbox{const.} \neq 0 \, .
\end{eqnarray*}
The manifold $(M,g)$ is a non-conformally flat manifold 
satisfying the following conditions \cite{SDHJK}: 
\begin{eqnarray*}
& &
\nabla _{X} S(Y,Z) + \nabla _{Y} S(Z,X) + \nabla _{Z} S(X,Y) \ = \ 0 ,\\
& &
S =  \frac{\kappa }{2}\, g -\frac{3 \kappa }{2}\, \eta \otimes \eta ,
\ \ \ \eta \ =\ ( 0, 0, -\rho, -x \rho, y \rho ) ,
\ \ \ \kappa \ =\ \rho ^{2} ,\ \ \
S^{2} \ =\ - \frac{\kappa }{2}\, S + \frac{\kappa^{2} }{2}\, g ,\\
%\end{eqnarray*}
%\begin{eqnarray*} 
& &
S \cdot R \ =\ 
2 \kappa\, R - \frac{\kappa}{2} g \wedge S 
+ \frac{\kappa^{2}}{4}\, g \wedge g ,\ \ \
C \cdot S \ = \ 0 , \ \ \
R \cdot S \ =\  - \frac{\kappa }{4}\, Q(g,S) ,\\ 
& & 
R \cdot R \ =\  - \frac{\kappa }{4}\, Q(g,R) ,\ \ \ 
R \cdot C \ = \  - \frac{\kappa }{4}\, Q(g,C) ,\\
& &
C \cdot R 
\ = \  - \frac{1}{3}\, Q(S,C) - \frac{\kappa }{3}\, Q(g,C) ,\ \ \ 
C \cdot C \ = \ C \cdot R ,\\
& &
R \cdot C + C \cdot R 
\ =\ - \frac{1}{3}\, Q(S,C) - \frac{7 \kappa }{12}\, Q(g,C) 
\end{eqnarray*}
and the condition of the form $(\ast)$
\begin{eqnarray*}
R \cdot C - C \cdot R 
&=& \frac{1}{3}\, Q(S,C) + \frac{\kappa }{12}\, Q(g,C) .
\end{eqnarray*}
The $(0,4)$-tensor $S \cdot R$ is defined by
\begin{eqnarray*}
(S \cdot R)(X,Y,W,Z) 
\ =\ 
R({\mathcal{S}}X,Y,W,Z) + R(X,{\mathcal{S}}Y,W,Z)
+ R(X,Y,{\mathcal{S}}W,Z) + R(X,Y,W,{\mathcal{S}}Z) .
\end{eqnarray*}
\end{example}

An important subclass of the class of partially Einstein manifolds,
of dimension $\geq 4$,
form some non-conformally flat and non-quasi-Einstein manifolds 
called Roter spaces. 
A semi-Riemannian manifold $(M,g)$, $\dim M = n \geq 4$, 
satisfying on 
${\mathcal U}_{S} \cap {\mathcal U}_{C} \subset M$ 
the following equation
\begin{eqnarray}
R &=& \frac{\phi}{2}\, S\wedge S 
+ \mu\, g\wedge S + \frac{\eta}{2}\, g \wedge g ,
\label{eq:h7a}
\end{eqnarray}
where 
$\phi$, $\mu $ and $\eta $ are some functions on this set,
is called a {\sl Roter type manifold}, or a {\sl Roter manifold}, 
or a {\sl Roter space} 
(see, e.g., {\cite[Section 15.5] {Chen-2021}},
\cite{{P106}, {2016_DGJZ}, {DGP-TV02}, {DHV2008}}).
Equation (\ref{eq:h7a}) is called a {\sl{Roter equation}}
(see, e.g., {\cite[Section 1] {DGHSaw-2022}}). 
Roter spaces and in particular Roter hypersurfaces 
in semi-Riemannian spaces of constant curvature were studied in:
\cite{{DDH-2021}, {P106}, {DGHHY},  {DGHZ01}, 
{DH-1998}, {R102}, {DeKow}, {DePlaScher}, {DeScher}, 
{G2}, {G5}, {Kow01}, {Kow02}}. 
We only mention that 
every Roter space $(M,g)$, $\dim M = n \geq 4$, 
satisfies on
${\mathcal{U}}_{S} \cap {\mathcal{U}}_{C} \subset M$  among others
the following conditions: 
(\ref{Roterformula}), (\ref{pseudo}) 
and (\ref{Riccipseudo07})-(\ref{genpseudo01})
(see, e.g., \cite{ {DGHSaw}, {2016_DGJZ}} and references therein).
We also refer to {\cite[Section 4] {2023_DGHP-TZ}} for a 
more detailed presentation of results on Roter spaces.

Let $(M,g)$, $n \geq 4$, be a non-partially-Einstein 
and non-conformally flat semi-Rieman\-nian manifold.
If its Riemann-Christoffel curvature $R$
is at every point of  
${\mathcal U}_{S} \cap {\mathcal U}_{C} \subset M$ 
a linear combination of the Kulkarni-Nomizu products 
formed by the tensors
$S^{0} = g$ and $S^{1} = S, S^{2}, \ldots , S^{p-1}, S^{p}$, 
where $p$ is some natural 
number $\geq 2$, then $(M,g)$
is called 
a {\sl generalized Roter type manifold},
or a {\sl generalized Roter manifold}, 
or a {\sl generalized Roter type space}, or 
a {\sl generalized Roter space}. For instance, when $p = 2$,
we have
\begin{eqnarray} 
R &=&  \frac{\phi _{2} }{2}\, S^{2} \wedge S^{2} 
+ \phi_{1}\, S \wedge S^{2}  
+ \frac{\phi}{2}\, S \wedge S
+ \mu _{1}\, g \wedge S^{2}
+ \mu \, g \wedge S
+ \frac{ \eta }{2} \, g \wedge g ,
\label{B001simply}
\end{eqnarray}
where 
$\phi $, 
$\phi _{1}$,
$\phi _{2}$,
$\mu_{1}$, $\mu$ and 
$\eta$ 
are functions on
${\mathcal U}_{S} \cap {\mathcal U}_{C}$.
Because $(M,g)$ is a non-partially Einstein manifold,
at least one of the functions 
$\mu _{1}$, $\phi _{1}$ and $\phi _{2}$ is a non-zero function.
Equation (\ref{B001simply}) is called a 
{\sl{Roter type equation}}  
(see, e.g., {\cite[Section 1] {DGHSaw-2022}}).
We refer to 
\cite{{DGHSaw-2022}, {2013_DGJP-TZ}, {2016_DGJZ}, 
{DGP-TV02}, {DeHoJJKunSh}, 
{S3}, {Saw114}, {SDHJK}, {2016_SK}, {2019_SK}}
for results on manifolds (hypersurfaces) satisfying (\ref{B001simply}).

If $(M,g)$, $\dim M = n \geq 4$, is a semi-Riemannian 
manifold satisfying  (\ref{4.3.012}) on ${\mathcal U}_{C} \subset M$, 
i.e., $C \cdot C = L_{C}\, Q(g,C)$, then (\ref{identity01}) yields
\begin{eqnarray}
C \cdot R  + R \cdot C 
&=& R \cdot R + L\, Q(g,C)  - \frac{1}{(n-2)^{2}}\, 
Q\left( g, 
g \wedge S^{2} - \kappa \, g\wedge S 
\right) 
\label{02identity01}
\end{eqnarray}
on ${\mathcal U}_{C}$.
In addition, 
if (\ref{genpseudo01}) is satisfied on ${\mathcal U}_{C}$, i.e.,
$R \cdot R - Q(S,R) = L\, Q(g,C)$,
then (\ref{02identity01}) turns into 
{\cite[Theorem 3.4 (ii)-(iii)] {2016_DGJZ}}
\begin{eqnarray}
C \cdot R + R \cdot C 
&=& 
Q(S,C) + ( L + L_{C} )\, Q(g,C)\nonumber\\
& & 
- \frac{1}{(n-2)^{2}}\, Q\left( g, 
g \wedge S^{2} + \frac{n-2}{2}\, S \wedge S - \kappa \, g \wedge S \right) .
\label{identity05}
\end{eqnarray}
Evidently, (\ref{identity05}),
by making use of (\ref{2022.11.10.aaa}) and the identity 
$Q(g, \frac{1}{2}\, g \wedge g) = Q(g,G) = 0$ takes the form
\begin{eqnarray*}
C \cdot R + R \cdot C 
&=& Q(S,C) + ( L + L_{C} )\, Q(g,C) 
- \frac{1}{(n-2)^{2}}\, Q( g, E) .
%\label{identity05}
\end{eqnarray*}

We also have (cf. {\cite[Theorem 3.4 (iv)-(v)] {2016_DGJZ}}):
if $(M,g)$, $\dim M =  n \geq 4$, 
is a quasi-Einstein semi-Riemannian manifold satisfying 
$\mathrm{rank}\, (S - \alpha\, g) = 1$ (i.e. (\ref{quasi02})), 
(\ref{4.3.012})
and
(\ref{genpseudo01})
then (\ref{identity05}),  
by making use of (\ref{quasi03}) and (\ref{quasi03quasi03}),  
yields 
\begin{eqnarray}
C \cdot R + R \cdot C  &=& Q(S,C) + ( L + L_{C} )\, Q(g,C) .
\label{identity05quasi}
\end{eqnarray}
In particular,
if $(M,g)$ is the G\"{o}del spacetime 
then (\ref{identity05quasi}) turns into
\begin{eqnarray*}
C \cdot R + R \cdot C &=& Q(S,C) + \frac{\kappa }{6}\, Q(g,C) .
%\label{Goedelidentity05}
\end{eqnarray*}

\begin{thm}
({cf. \cite[Proposition 3.2, Theorem 3.3, Theorem 4.4] {DGHHY}})
If $(M,g)$, $\dim M = n \geq 4$, is a semi-Riemannian manifold
satisfying on 
${\mathcal U}_{S} \cap {\mathcal U}_{C} \subset M$ conditions
(\ref{4.3.012}), (\ref{genpseudo01}) and 
\begin{eqnarray}
R \cdot S  &=& Q(g,D) ,
\label{eqD}
\end{eqnarray}
where $D$ is a symmetric $(0,2)$-tensor,
then $R \cdot R = L_{R}\, Q(g,R)$ 
on this set, where $L_{R}$ is some function on this set.
Moreover, at every point of ${\mathcal U}_{S} \cap {\mathcal U}_{C}$
we have
\newline
(i)
$\mathrm{rank}\, (S - \alpha _{1}\, g) = 1$ 
and $\alpha _{1} 
= \frac{1}{2} \left( \frac{\kappa }{n-1} - L + L_{C} \right)$, or
\newline
(ii)
$\mathrm{rank}\, (S - \alpha _{1}\, g) \geq 2$ and  
$\alpha _{1} = \frac{1}{2}  \left( \frac{\kappa }{n-1} - L + L_{C} \right)$
and (\ref{eq:h7a}).
\end{thm}

\begin{cor}
{\cite[Corollary 3.6] {2016_DGJZ}}
If $(M,g)$, $\dim M = n \geq 4$, is a semi-Riemannian manifold
satisfying on 
${\mathcal U}_{S} \cap {\mathcal U}_{C} \subset M$ conditions
(\ref{4.3.012}), (\ref{genpseudo01}) and (\ref{eqD}) then 
\begin{eqnarray*}
C \cdot R  + R \cdot C &=& Q(S,C) + L_{2} \, Q(g,C)  
\end{eqnarray*}
on ${\mathcal U}_{S} \cap {\mathcal U}_{C}$, 
where $L_{2}$ is some function on this set.
\end{cor}

\begin{thm} {\cite[Theorem 3.1] {2016_DGJZ}}
If $(M,g)$, $\dim M = n \geq 4$, is a pseudosymmetric 
Einstein semi-Riemannian manifold satisfying 
$R \cdot R = L_{R}\, Q(g,R)$
(i.e. (\ref{pseudo})) 
on ${\mathcal{U}}_{R} \subset M$, then on this set 
\begin{eqnarray*}
R \cdot R - Q(S,R) &=& \left(  L_{R} - \frac{\kappa}{n} \right) Q(g,C),\\
C \cdot C &=& \left(  L_{R} - \frac{\kappa}{(n-1) n} \right) Q(g,C),\\
C \cdot R + R \cdot C 
&=& Q(S,C) + \left( 2 L_{R} - \frac{\kappa }{n-1} \right) Q(g,C).
\end{eqnarray*}
\end{thm}

\begin{thm} {\cite[Theorem 4.1] {2016_DGJZ}}
Let $\overline{M} \times _{F} \widetilde{N}$ 
be the warped product manifold 
with a $1$-dimensional manifold $(\overline{M},\overline{g})$,
$\overline{g}_{11} = \pm 1$,
and a $3$-dimensional semi-Riemannian manifold 
$(\widetilde{N}, \widetilde{g})$, and a warping function $F$.
If $(\widetilde{N}, \widetilde{g})$ 
is not a space of constant curvature then 
(\ref{genpseudo01}) and (\ref{02identity01}) 
hold 
on ${\mathcal U}_{C} \subset \overline{M} \times _{F} \widetilde{N}$.
Moreover, if $(\widetilde{N}, \widetilde{g})$ 
is a quasi-Einstein manifold then  
(\ref{4.3.012}) and (\ref{identity05})
hold 
on ${\mathcal U}_{S} \cap {\mathcal U}_{C} 
\subset \overline{M} \times _{F} \widetilde{N}$.
\end{thm}

The last theorem leads to the following results.

\begin{thm} {\cite[Theorem 4.3] {2016_DGJZ}}
Let $\overline{M} \times _{F} \widetilde{N}$ 
be the warped product manifold 
with a $1$-dimensional manifold 
$(\overline{M},\overline{g})$, $\overline{g}_{11} = \pm 1$,
and an $(n-1)$-dimensional quasi-Einstein semi-Riemannian manifold
$(\widetilde{N}, \widetilde{g})$, $n \geq 4$, 
and a warping function $F$, 
and let $(\widetilde{N}, \widetilde{g})$ be 
a conformally flat manifold, when $n \geq 5$.
Then conditions (\ref{4.3.012}), (\ref{genpseudo01}) 
and (\ref{identity05}) 
are satisfied on ${\mathcal U}_{S} \cap {\mathcal U}_{C} 
\subset \overline{M} \times _{F} \widetilde{N}$.
\end{thm}

\begin{thm} {\cite[Theorem 4.4] {2016_DGJZ}}
Let $\overline{M} \times \widetilde{N}$ be the product manifold 
with a $1$-dimensional manifold $(\overline{M},\overline{g})$, 
$\overline{g}_{11} = \pm 1$,
and an $(n-1)$-dimensional quasi-Einstein 
semi-Riemannian manifold $(\widetilde{N}, \widetilde{g})$, $n \geq 4$,
satisfying 
$\mathrm{rank}\, (\widetilde{S} - \rho\, \widetilde{g}) = 1$
on ${\mathcal U}_{\widetilde{S}} \subset \widetilde{M}$, where $\rho $ 
is some function on ${\mathcal U}_{\widetilde{S}}$,
and let $(\widetilde{N}, \widetilde{g})$ be 
a conformally flat manifold, when $n \geq 5$.
Then on ${\mathcal U}_{S} \cap {\mathcal U}_{C} 
\subset \overline{M} \times \widetilde{N}$ we have
\begin{eqnarray*}
(n-3)(n-2)\rho \, C 
&=& 
g \wedge S^{2} + \frac{n-2}{2}\, S \wedge S - \kappa \, g\wedge S 
+ (n-2) \rho \left( \frac{2 \kappa }{n-1} - \rho \right) G  ,
%\label{new-quasi11}
\\
C \cdot R + R \cdot C &=& Q(S,C) 
+ \left( \frac{\kappa }{(n-2)(n-1)} - \rho \right) Q(g,C) .
%\label{new-quasi10}
\end{eqnarray*}
\end{thm}

We present now the following results 
on Einstein and quasi-Einstein manifolds.
Using (\ref{2020.10.3.c}) we can easily check that equation
\begin{eqnarray}
g \wedge S^{2} + \frac{n-2}{2} \, S \wedge S - \kappa \, g \wedge S
+ \frac{\kappa ^{2} - \mathrm{tr} (S^{2})}{2(n-1)} \, g \wedge g &=& 0 
\label{2022.11.20.aa}
\end{eqnarray}
is satisfied on any Einstein manifold $(M,g)$. Moreover, we also have

\begin{thm} {\cite[Lemma 2.1] {DGHSaw-2022}}
If $(M,g)$, $n \geq 4$, is a quasi-Einstein manifold satisfying
(\ref{quasi02}) 
on ${\mathcal U}_{S} \cap  {\mathcal U}_{C} \subset M$ 
then (\ref{2022.11.20.aa}) holds on this set.
\end{thm}

We finish this section with the following results.

\begin{thm} {\cite[Lemma 2.2] {DGHSaw-2022}}
If $(M,g)$, $n \geq 4$, is a Roter space satisfying
(\ref{eq:h7a}) 
on ${\mathcal U}_{S} \cap  {\mathcal U}_{C} \subset M$ then 
\begin{eqnarray*}
C &=& \frac{\phi}{n-2} 
\left(g \wedge S^{2} + \frac{n-2}{2} \, S \wedge S - \kappa \, g \wedge S
+ \frac{\kappa ^{2} - \mathrm{tr} (S^{2})}{2(n-1)}
\, g \wedge g \right) 
%\label{2022.11.22.aa}
\end{eqnarray*}
on this set.
\end{thm}

\begin{rem} (i) In view of {\cite[Lemma 3.2 (ii)] {D-1992}},
we can state that the following identity is satisfied 
on every semi-Riemannian manifold 
$(M,g)$, $n = \dim M \geq 3$,
with vanishing Weyl conformal curvature tensor $C$
\begin{eqnarray} 
R \cdot R - Q(S,R) &=&
\frac{1}{(n-2)^{2}}\,
Q(g, g \wedge S^{2} + \frac{n-2}{2} \, S \wedge S - \kappa \, g \wedge S) .
\label{2022.11.20.zz}
\end{eqnarray}
(ii) On every $3$-dimensional semi-Riemannian manifold $(M,g)$
the identity $R \cdot R = Q(S,R)$ is satisfied 
{\cite[Theorem 3.1] {D-1992}}. Thus (\ref{2022.11.20.zz}) reduces to 
\begin{eqnarray*} 
Q(g, g \wedge S^{2} + \frac{n-2}{2} \, S \wedge S - \kappa \, g \wedge S)
&=& 0 .
\end{eqnarray*}
From this, by suitable contractions we easily get (\ref{2022.11.20.aa}).
\newline
(iii) From (i) we easily deduce that
on every semi-Riemannian conformally flat manifold 
$(M,g)$, $n = \dim M \geq 4$, conditions: $R \cdot R = Q(S,R)$
and (\ref{2022.11.20.aa}) are equivalent.
\end{rem}

\begin{rem}
Warped product manifolds  $\overline{M} \times _{F} \widetilde{N}$, 
with a $1$-dimensional base manifold 
$(\overline{M}, \overline{g})$, $\overline{g}_{11} = \pm 1$,  
and an $(n-1)$-dimensional Einsteinian or non-Einsteinian fiber 
$(\widetilde{N},\widetilde{g})$, $n \geq 4$, 
satisfying the condition 
\begin{eqnarray*} 
R \cdot C - C \cdot R &=& L\, Q(S,R) ,
\end{eqnarray*}
were studied in \cite{P119}.
\end{rem}

\section{Warped product manifolds with 2-dimensional base manifold}

Let $\overline{M} \times _{F} \widetilde{N}$ be the warped product manifold 
with a $2$-dimensional semi-Riemannian manifold 
$(\overline{M},\overline{g})$
and an $(n-2)$-dimensional semi-Riemannian manifold 
$(\widetilde{N},\widetilde{g})$, 
$n \geq 4$, and a warping function $F$, 
and let $(\widetilde{N},\widetilde{g})$ 
be a space of constant curvature when $n \geq 5$. 

Let $S_{hk}$ and $C_{hijk}$ be the local components of the Ricci tensor $S$ 
and the Weyl conformal curvature tensor $C$ 
of $\overline{M} \times _{F} \widetilde{N}$, respectively. 
We have
\begin{eqnarray}
S_{ad} &=& \frac{\overline{\kappa } }{2}\, g_{ab} 
- \frac{n-2}{2 F}\, T_{ab} ,\ \  
S_{\alpha \beta } \ =\ \tau _{1} \, g_{\alpha \beta } ,\ \ 
S_{a \alpha } \ =\ 0,
\label{RcciRicci01}\\
\tau _{1} &=&  
\frac{\widetilde{\kappa} }{(n-2) F} - \frac{\mathrm{tr}(T) }{2F }  
- (n-3)\, \frac{\Delta _1 F }{4 F^{2}} ,\ \ \
\Delta _{1} F \ =\  {\Delta}_{1\, \overline{g}} F\ =\ 
\overline{g}^{ab} F_a F_b ,
\label{2022.11.21.aa}\\
T_{ab}&=& 
\overline{\nabla }_a F_b - \frac{1}{2F} F_a F_b ,\ \ 
\mathrm{tr}(T) \ =\ \overline{g}^{ab} T_{ab} ,\nonumber
%\label{2022.11.21.bb}
\end{eqnarray}
where $T$ is the $(0,2)$-tensor with the local components $T_{ab}$.
We also have {\cite[eqs. (5.10)-(5.14)] {2016_DGJZ}}
\begin{eqnarray}
C_{abcd} 
&=& \frac{n-3}{n-1}\, \rho_{1}\, G_{abcd} 
\ =\ \frac{n-3}{n-1}\, \rho_{1}\, ( g_{ad}g_{bc} - g_{ac}g_{bd} ),
\label{WeylWeyl01.a}\\
%\end{eqnarray}
%\begin{eqnarray}
C_{\alpha bc \beta} 
&=& - \frac{n-3}{(n-2) (n-1)}\, \rho_{1}\, G_{\alpha bc \beta } 
\ =\ - \frac{n-3}{(n-2) (n-1) \rho_{1} }\, \rho_{1}\, 
g_{bc} g_{\alpha \beta} ,
\label{WeylWeyl01.b}\\
%\end{eqnarray}
%\begin{eqnarray}
C_{\alpha \beta \gamma \delta} 
&=&  \frac{2 \rho_{1} }{(n-2) (n-1) }\, G_{\alpha \beta \gamma \delta }
\ =\ 
\frac{2 \rho_{1} }{(n-2) (n-1) }\, 
( g_{\alpha \delta} g_{\beta \gamma } - g_{\alpha \gamma } 
g_{ \beta  \delta} ) ,
\label{WeylWeyl01.c}\\
%\end{eqnarray}
%\begin{eqnarray}
C_{abc\delta } 
&=& C_{ab \alpha \beta } \ =\ C_{a \alpha \beta \gamma } \ =\ 0 ,
\label{WeylWeyl01}
\end{eqnarray}
where $G_{hijk} = g_{hk}g_{ij} - g_{hj}g_{ik}$ and
\begin{eqnarray*}
\rho_{1} &=& 
\frac{ \overline{\kappa } }{2} 
+ \frac{ \widetilde{\kappa } }{ (n-3)(n-2) F } 
+ \frac{1}{2 F} \left( \Delta F - \frac{\Delta _1 F}{ F} \right) ,\ \ \
\Delta F \ =\  \overline{g}^{ab} \overline{\nabla }_a F_b\, .
%\label{WeylWeyl02}
\end{eqnarray*}
If we set 
{\cite[eqs. (5.13)] {2016_DGJZ}} 
\begin{eqnarray}
\rho &=&  \frac{2 (n-3) }{ n-1} \, \rho_{1} \ = \
\frac{2 (n-3) }{ n-1}
\left(  
\frac{ \overline{\kappa } }{2} 
+ \frac{ \widetilde{\kappa } }{ (n-3)(n-2) F } 
+ \frac{1}{2 F} \left( \Delta F - \frac{\Delta _1 F}{ F} \right)
 \right)
\label{WeylWeyl06}
\end{eqnarray}
then (\ref{WeylWeyl01.a})-(\ref{WeylWeyl01}) turn into 
{\cite[eqs. (5.14)] {2016_DGJZ}} 
\begin{eqnarray*}
C_{abcd} &=& \frac{\rho }{2}\, G_{abcd} ,\ \ \
C_{\alpha bc \beta} 
\ =\ - \frac{\rho}{2 (n-2)}\,  G_{\alpha bc \beta } ,\nonumber\\
C_{\alpha \beta \gamma \delta} &=& \frac{\rho }{(n-3) (n-2)}\,   
G_{\alpha \beta \gamma \delta } ,\ \ \
C_{abc\delta } \ =\ C_{ab \alpha \beta } 
\ =\ C_{a \alpha \beta \gamma } \ =\ 0 .
%\label{WeylWeyl03}
\end{eqnarray*}
Further, 
by making use of the formulas
for 
the local components $(C \cdot C)_{hijklm}$ and $Q(g,C)_{hijklm}$ 
of the tensors $C \cdot C$ and $Q(g,C)$, i.e.
\begin{eqnarray*}
(C \cdot C)_{hijklm} 
&=&
g^{rs}(
C_{rijk}C_{ shlm}
+ C_{hrjk}C_{ silm}
+ C_{hirk}C_{sjlm}
+ C_{hijr}C_{sklm}) ,\\
Q(g,C)_{hijklm}
&=&
g_{hl}C_{ mijk} 
+ g_{il}C_{ hmjk} 
+ g_{jl}C_{ himk} 
+ g_{kl}C_{ hijm}\\ 
& &-
g_{hm}C_{ lijk} 
- g_{im}C_{ hljk} 
- g_{jm}C_{ hilk} 
- g_{km}C_{ hijl}  ,
\end{eqnarray*}
we obtain {\cite[eqs. (7.7)-(7.8)] {2016_DGJZ}} 
\begin{eqnarray*}
& &
(C \cdot C)_{\alpha abcd \beta } 
\ =\
 - \frac{(n-1)\rho^{2}}{4 (n-2)^{2}}\, g_{\alpha \beta} G_{dabc},\ \ \
(C \cdot C)_{a \alpha \beta \gamma d \delta } 
\ =\
\frac{(n-1)\rho^{2}}{4 (n-2)^{2}(n-3)}\, 
g_{ad} G_{\delta \alpha \beta \gamma },\nonumber\\
& &
Q(g,C)_{\alpha abcd \beta } 
\ =\ 
\frac{(n-1)\rho }{2 (n-2)}\, g_{\alpha \beta} G_{dabc}, \ \  \
Q(g,C)_{a \alpha \beta \gamma d \delta } 
\ =\
- \frac{(n-1)\rho }{2 (n-2)(n-3)}\,  g_{ad} 
G_{\delta \alpha \beta \gamma }.
\end{eqnarray*}

\begin{thm} {\cite[Theorem 7.1 (i)] {2016_DGJZ}}
Let $\overline{M} \times _{F} \widetilde{N}$ be the warped product manifold 
with a $2$-dimensional semi-Riemannian manifold 
$(\overline{M},\overline{g})$
and an $(n-2)$-dimensional semi-Riemannian manifold 
$(\widetilde{N},\widetilde{g})$, 
$n \geq 4$, and a warping function $F$, 
and let $(\widetilde{N},\widetilde{g})$ 
be a space of constant curvature, when $n \geq 5$. 
\newline
(i) The following equation is satisfied on 
${\mathcal U}_{C} \subset \overline{M} \times _{F} \widetilde{N}$ 
\begin{eqnarray}
C \cdot C &=& L_{C} \, Q(g,C)\, ,\nonumber\\
L_{C} &=& - \frac{\rho }{2(n-2)}
\ = \
- \frac{ n-3 }{(n-2)(n-1)}
\left(  
\frac{ \overline{\kappa } }{2} 
+ \frac{ \widetilde{\kappa } }{ (n-3)(n-2) F } 
+ \frac{1}{2 F} \left( \Delta F - \frac{\Delta _1 F}{ F} \right)
 \right),
\label{function1} 
\end{eqnarray}
where the function $\rho$ is defined by (\ref{WeylWeyl06}).
\newline
(ii) Equation (\ref{genpseudo01}) is satisfied on  
${\mathcal U}_{C} \subset \overline{M} \times _{F} \widetilde{N}$, 
where the functions $\tau _{1}$ and $L$ are defined by
(\ref{2022.11.21.aa}) and
\begin{eqnarray}
L &=&
- \frac{n-2}{ (n-1) \rho }   
\left(
\overline{\kappa} \left( \tau_{1} + \frac{ \mathrm{tr}(T) }{2F } \right)
+ \frac{n-3}{4 F^{2} } \left( \mathrm{tr}(T^{2}) 
- (\mathrm{tr}(T))^{2} \right) \right) ,
\label{function2}
\end{eqnarray}
$T$ is the $(0,2)$-tensor with the local components 
$T_{ab} = \overline{\nabla }_a F_b - \frac{1}{2F} F_a F_b $, 
$\mathrm{tr}(T) \ =\ \overline{g}^{ab} T_{ab} $,
$T^{2}_{ad} = T_{ac } \overline{g}^{cd} T_{db}$
and 
$\mathrm{tr}(T^{2}) =  \overline{g}^{ab} T^{2}_{ab}$. 
\newline
(iii) The following equation is satisfied on 
${\mathcal U}_{C} \subset \overline{M} \times _{F} \widetilde{N}$ 
\begin{eqnarray*}
C \cdot R + R \cdot C &=& Q(S,C) + ( L_{C} + L )\, Q(g,C)
 - \frac{1}{(n-2)^{2}}\, 
Q\left( g, 
g \wedge S^{2}
+ \frac{n-2}{2}\, S \wedge S - \kappa \, g \wedge S \right) ,
%\label{identity05}
\end{eqnarray*} 
where $L_{C}$ and $L$ are functions defined by 
(\ref{function1}) and
(\ref{function2}).
\end{thm}

%\vspace{3mm}

We have (see, eq. (\ref{RcciRicci01}))
\begin{eqnarray*}
S_{ad} &=& \frac{\overline{\kappa } }{2}\, g_{ab} 
- \frac{n-2}{2 F}\, T_{ab} ,\ \  
S_{\alpha \beta } 
\ =\ \tau _{1} \, g_{\alpha \beta } ,\ \ S_{a \alpha } \ =\ 0,
\end{eqnarray*}
where $\tau _{1}$ is defined by (\ref{2022.11.21.aa}).

We define now on 
${\mathcal U}_{S} \subset \overline{M} \times _{F} \widetilde{N}$ 
the $(0,2)$-tensor $A$ by $A = S - \tau_{1}\, g$.
We can check that 
$\mathrm{rank}(A) = 2$ at a point of ${\mathcal U}_{S}$ 
if and only if $\mathrm{tr}(A^{2}) - (\mathrm{tr}(A))^{2} \neq 0$ 
at this point {\cite[Section 6] {2016_DGJZ}}. At all points of 
${\mathcal U}_{S}$, at which $\mathrm{rank}(A) = 2$,
we set
\begin{eqnarray}
\tau _{2} &=& (\mathrm{tr}(A^{2}) - (\mathrm{tr}(A))^{2})^{-1}\, .
\label{2022.11.23.ff}
\end{eqnarray}
Let $V$ be the set of all points of 
${\mathcal U}_{S} \cap {\mathcal U}_{C}$ at which:
$\mathrm{rank}(A) = 2$ and $S_{ad}$ is not proportional to $g_{ad}$.

\begin{thm} {\cite[Theorem 7.1 (ii)] {2016_DGJZ}}
Let $\overline{M} \times _{F} \widetilde{N}$ be the warped product manifold 
with a $2$-dimensional semi-Riemannian manifold 
$(\overline{M},\overline{g})$
and an $(n-2)$-dimensional semi-Riemannian manifold 
$(\widetilde{N},\widetilde{g})$, 
$n \geq 4$, and a warping function $F$, 
and let $(\widetilde{N},\widetilde{g})$ 
be a space of constant curvature when $n \geq 5$. 
Then on the set $V \subset {\mathcal U}_{S} \cap {\mathcal U}_{C}$ 
(defined above) we have: 
\begin{eqnarray*}
C &=& - \frac{ (n-1) \rho \tau _{2} }{(n-3)(n-2)}  
\left( 
g \wedge S^{2} + \frac{ n-2 }{2}  S \wedge S - \kappa  g \wedge S 
+ \frac{ \kappa^{2} - \mathrm{tr}(S^{2})  }{n-1} G \right) ,
%\label{mainTT032main}
\end{eqnarray*}
\begin{eqnarray*}
R \cdot C + C \cdot R 
&=& Q(S,C)
+ \left( L - \frac{\rho }{2 (n-2)} 
+ \frac{n-3 }{ (n-2) (n-1) \rho \tau_{2} } \right) Q(g,C) ,
%\label{identity06dede}
\end{eqnarray*}
where $\rho$, $\tau _{1}$ and $\tau _{2}$ are defined by
(\ref{WeylWeyl06}), (\ref{2022.11.21.aa}) and (\ref{2022.11.23.ff}),
respectively.
\end{thm}

\begin{thm} {\cite[Theorem 6.2] {2016_DGJZ}}
Let $\overline{M} \times _{F} \widetilde{N}$ be the warped product manifold 
with a $2$-dimensional semi-Riemannian manifold 
$(\overline{M},\overline{g})$
and an $(n-2)$-dimensional semi-Riemannian manifold 
$(\widetilde{N},\widetilde{g})$, 
$n \geq 4$, and a warping function $F$, 
and let $(\widetilde{N},\widetilde{g})$ be an Einstein space, when 
$n \geq 5$. 
On the set $V \subset {\mathcal U}_{S} \cap {\mathcal U}_{C}$ we have 
\begin{eqnarray*}
R \cdot S &=&
(\phi _{1} - 2 \tau_{1} \phi_{2} + \tau _{1}^{2} \phi_{3}) \, Q(g,S)
+ (\phi_{2} - \tau_{1} \phi_{3})\, Q(g, S^{2}) + \phi _{3}\, Q(S,S^{2}) ,
%\label{2022.11.25aa}
\\
\phi _{1} &=& \frac{2 \tau _{1} - \overline{\kappa} }{2 (n-2) } ,\ \ \ 
\phi _{2} = \frac{1}{n-2}, \ \ \
\phi _{3} = \frac{ \tau_{2} (2 \kappa - \overline{\kappa} 
- 2 (n-1) \tau_{1}) }{n-2} ,
%\label{2022.11.25bb}
\end{eqnarray*}
where $\tau _{1}$ is defined by (\ref{2022.11.21.aa}).
\end{thm}

Theorems 4.10 and 4.11 imply  
\begin{thm} {\cite[Proposition 2.3] {DGHSaw-2022}}
If $(M,g)$, $n \geq 4$, is a semi-Riemannian manifold satisfying
(\ref{quasi02}) or (\ref{eq:h7a}) at every point 
of ${\mathcal U}_{S} \cap  {\mathcal U}_{C} \subset M$ 
then the following equation is satisfied on this set  
\begin{eqnarray}
\tau \, C &=&  
g \wedge S^{2} + \frac{n-2}{2} \, S \wedge S - \kappa \, g \wedge S
+ \frac{\kappa ^{2} - \mathrm{tr} (S^{2})}{2(n-1)}
\, g \wedge g ,
\label{2022.11.22.jj}
\end{eqnarray}
where $\tau$ is some function 
on ${\mathcal U}_{S} \cap {\mathcal U}_{C}$.
\end{thm}

Theorem 5.4, {\cite[Theorem 7.1 (ii)] {2016_DGJZ}} and
{\cite[Theorem 4.1] {DeKow}} imply 
\begin{thm} {\cite[Theorem 2.4] {DGHSaw-2022}}
Let $\overline{M} \times _{F} \widetilde{N}$ 
be the warped product manifold 
with a $2$-dimensional semi-Riemannian manifold 
$(\overline{M},\overline{g})$,
an $(n-2)$-dimensional semi-Riemannian manifold 
$(\widetilde{N},\widetilde{g})$, 
$n \geq 4$, a warping function $F$, 
and let $(\widetilde{N},\widetilde{g})$ 
be a space of constant curvature when $n \geq 5$. 
Then (\ref{2022.11.22.jj}) holds on
${\mathcal U}_{S} \cap {\mathcal U}_{C} 
\subset \overline{M} \times _{F} \widetilde{N}$.
\end{thm}

\begin{rem}Let $\overline{M} \times _{F} \widetilde{N}$ 
be the warped product manifold 
with a $2$-dimensional semi-Riemannian manifold 
$(\overline{M},\overline{g})$
and an $(n-2)$-dimensional semi-Riemannian manifold 
$(\widetilde{N},\widetilde{g})$, 
$n \geq 4$, and a warping function $F$, 
and let $(\widetilde{N},\widetilde{g})$ 
be a space of constant curvature when $n \geq 5$. 
\newline
(i) 
From {\cite[Theorem 4.1] {DeKow}} it follows that
at all points of the set ${\mathcal U}_{S} \cap {\mathcal U}_{C}$, 
at which 
$S_{ad}$ is proportional to $g_{ad}$ and $\mathrm{rank}(A) = 2$,
the Riemann-Christoffel curvature tensor $R$
is a linear combination 
of the Kulkarni-Nomizu products 
$S \wedge S$, $g \wedge S$ and $g \wedge g$,
i.e., (\ref{eq:h7a}) is satisfied. 
Thus, in view of 
{\cite[Theorem 6.7] {DGHSaw}} 
(see also {\cite[Theorem 3.2] {2016_DGJZ}}),
we have $R \cdot R = L_{R}\, Q(g,R)$
and in a consequence $R \cdot S = L_{R}\, Q(g,S)$, for some function $L_{R}$.
\newline
(ii) On the set $V \subset {\mathcal U}_{S} \cap {\mathcal U}_{C}$ 
we also have {\cite[Section 7] {2016_DGJZ}}
\begin{eqnarray*}
R \cdot C 
&=& 
Q(S,C) + \left( L + \frac{n-3 }{ (n-2) (n-1) \rho \tau_{2} } \right) Q(g,C)
+ \frac{ (n-1)\, \rho \tau _{2} }{(n-2)^{2}}\, g \wedge  Q(S , S^{2}) 
\nonumber\\
& & 
+ \frac{ 1 }{(n-2)^{2}}\, Q\left( \left(  \frac{ \rho  }{2} 
+ (n-1)\, \rho  \tau _{1}^{2} \tau _{2} \right) S
- (n-1)\, \rho \tau _{1} \tau _{2}\, S^{2},G \right) ,\\
%\label{identity07}
C \cdot R &=& 
- \frac{\rho }{2 (n-2)}\, Q(g,C) 
- \frac{ (n-1)\, \rho \tau _{2} }{(n-2)^{2}}\, 
g \wedge  Q(S , S^{2})\nonumber\\
& & 
- \frac{ 1 }{(n-2)^{2}}\, Q\left(  \left( \frac{ \rho  }{2} 
+ (n-1)\, \rho  \tau _{1}^{2} \tau _{2} \right) S
- (n-1)\, \rho \tau _{1} \tau _{2}\, S^{2} , G \right) ,
%\label{geneinst01}
\end{eqnarray*}
and in a consequence 
\begin{eqnarray*}
R \cdot C - C \cdot R
&=&
Q(S,C)
+ \left( L + \frac{n-3 }{ (n-2) (n-1) \rho \tau_{2} }
+ \frac{\rho }{2 (n-2)}  \right) Q(g,C)\nonumber\\
& & 
+ \frac{ 2 }{(n-2)^{2}}\, Q\left( \left(  \frac{ \rho  }{2} 
+ (n-1)\, \rho  \tau _{1}^{2} \tau _{2} \right) S
- (n-1)\, \rho \tau _{1} \tau _{2}\, S^{2},G \right) \nonumber\\
& & 
+ \frac{2 (n-1)\, \rho \tau _{2} }{(n-2)^{2}}\, g \wedge  Q(S , S^{2}), 
%\label{identity07aa}
\end{eqnarray*}
where $\rho$, $\tau _{1}$ and $\tau _{2}$ are defined by
(\ref{WeylWeyl06}), (\ref{2022.11.21.aa}) and (\ref{2022.11.23.ff}),
respectively.
\newline
(iii)
A short presentation on quasi-Einstein and $2$-quasi-Einstein 
warped product manifolds satisfying conditions of the form $(\ast)$ 
is given in \cite{DGHS-2022} (see also \cite{DG-2023}).
\end{rem}

\section{Semi-Riemannian manifolds satisfying the condition $R \cdot C - C \cdot R = L\, Q(S,C)$}

We can check that on any Einstein manifold $(M,g)$, $\dim M = n \geq 4$, 
the  tensors 
$Q(g,R)$, $Q(S,R)$, $Q(g,C)$ and $Q(S,C)$ satisfy
\begin{eqnarray}
\frac{\kappa}{n}\, Q(g,R) \ =\ Q(S,R) &=& Q(S,C) 
\ =\  \frac{\kappa }{n}\, Q(g,C) .
\label{tachib01}
\end{eqnarray}
Further, in 
{\cite[Theorem 3.1] {DHS01}} 
it was stated 
that on every Einstein manifold $(M,g)$, $\dim M = n \geq 4$, 
the following identity is satisfied
\begin{eqnarray}
R\cdot C - C \cdot R &=&  \frac{\kappa }{(n-1)n}\, Q(g,R) .
\label{tachib01a}
\end{eqnarray}
The remarks above lead to the problem of investigation 
of curvature properties of non-Einstein and non-conformally flat 
semi-Riemannian manifolds $(M,g)$, $\dim M = n \geq 4$,
satisfying at every point of $M$
the curvature condition of the following form: 
the difference tensor $R\cdot C - C \cdot R$ is proportional 
to $Q(g,R)$, $Q(S,R)$, $Q(g,C)$ and $Q(S,C)$.
Such conditions are strongly related to some pseudosymmetry 
type curvature conditions, 
see, e.g., \cite{DGHSaw} and references therein.

Semi-Riemannian manifolds $(M,g)$, $\dim M = n \geq 4$, 
satisfying at every point of $M$ the following condition
\begin{eqnarray}
\mbox{the tensors}\ \ R\cdot C - C \cdot R\ \ 
\mbox{and}\ \ Q(S,C)\ \ \mbox{are linearly dependent},
\label{2022.08.08aa}
\end{eqnarray}
were investigated in \cite{DGHZ01}.
It is obvious that 
(\ref{2022.08.08aa})
is satisfied at every point of $M$ 
at which $C$ vanishes. 
It is also clear that 
(\ref{tachib01}) and (\ref{tachib01a}) imply  that 
\begin{eqnarray*}
R \cdot C - C \cdot R &=&  \frac{1}{n-1}\, Q(S,C) 
%\label{tachib02}
\end{eqnarray*}
holds on any Einstein manifold $(M,g)$, $\dim M = n \geq 4$. Therefore  
we will restrict our considerations to manifolds $(M,g)$, 
$\dim M = n \geq 4$, 
satisfying (\ref{2022.08.08aa})
on the set 
${\mathcal U}_{S} \cap {\mathcal U}_{C} \subset M$.
Thus on ${\mathcal U}_{S} \cap {\mathcal U}_{C}$ 
we have (\ref{advances4}), i.e.,
$R\cdot C - C \cdot  R =  L\, Q(S,C)$,  where $L$ is some function on this set.
We mention that if the tensor $R\cdot C - C \cdot R$ vanishes 
on ${\mathcal U}_{S} \cap {\mathcal U}_{C}$ then
on this set {\cite[Theorem 4.1] {DHS01}}: 
\begin{eqnarray}
R\cdot C \ =\ C \cdot R &=& 0 ,
\label{tachib0278}
\\
%\end{eqnarray}
%\begin{eqnarray}
Q(S,C) &=& 0 .
%\label{tachib0278add}
\nonumber
\end{eqnarray}
On the other hand, if  
$Q(S,C)$ vanishes on ${\mathcal U}_{S} \cap {\mathcal U}_{C}$ then
at every point $x \in {\mathcal U}_{S} \cap {\mathcal U}_{C}$ we have:
\newline
(i) if $\mbox{rank}\, S = 1$ at $x$ then $\kappa = 0$ 
and (\ref{tachib0278}) hold at $x$ {\cite[Section 3] {DGHZ01}}, or
\newline
(ii) if $\mbox{rank}\, S > 1$ at $x$ then $C \cdot R = 0$ 
and $R\cdot C  =\frac{\kappa }{n-1}\, Q(g,C)$ hold at $x$ 
{\cite[Section 4] {DGHZ01}}. 
\newline
Thus we see that in the case (ii), if a manifold satisfies 
(\ref{advances4}) then its scalar curvature must vanish
on ${\mathcal U}_{S} \cap {\mathcal U}_{C}$.

The main result of 
{\cite[Section 3] {DGHZ01}}  
states that pseudosymmetric manifolds satisfying 
some additional curvature conditions are quasi-Einstein manifolds 
satisfying the conditions:
$C \cdot C  = 0$, $C \cdot R = 0$, 
and (\ref{advances4}) with $L = \frac{1}{n-1}$. 
Precisely we have 
\begin{thm}
{\cite[Theorem 3.4] {DGHZ01}} 
Let $(M,g)$, $\dim M = n \geq 4$,
be a semi-Riemannian manifold. 
If the following conditions: 
\begin{eqnarray*}
R \cdot R &=& \frac{\kappa }{(n-1)n} \, Q(g, R) ,
%\label{RRDD01}
\nonumber\\
R \cdot R - Q(S,R) &=& - \frac{(n-2) \kappa }{(n-1)n} \, Q(g, C) ,
%\label{RRDD02}
\nonumber\\
R \cdot C &=& \frac{1}{n-1} \, Q(S, C)
%\label{RRDD03}
\nonumber
\end{eqnarray*}
are satisfied 
on ${\mathcal U}_{S} \cap {\mathcal U}_{C} \subset M$, 
then on this set we have: 
\begin{eqnarray*}
C \cdot C  &=& 0 ,
%\label{RRDD03ggg}\\
\nonumber\\ 
\mbox{rank}\, ( S - \frac{\kappa }{n}\, g) &=& 1 , 
%\label{RRDD03bbb}\\
\nonumber\\
C \cdot R &=& 0 ,
%\label{RRDD03fff}
\\
(n-1)\, (R \cdot C - C \cdot R) &=& Q( S,C) . 
%\label{RRDD03aa}
\nonumber
\end{eqnarray*}
\end{thm}
\begin{thm}
{\cite[Proposition 3.9] {DGHZ01}} 
Let $T$ be a generalized curvature tensor 
on a semi-Riemannian manifold $(M,g)$, $\dim M = n \geq 4$.
\newline
(i) 
If the conditions
\begin{eqnarray}
Q(Ric(T),T) &=& 0 ,
%\label{WDWD01}
\nonumber\\
\mbox{rank}\,( Ric(T) ) &=& 1 
\label{WDWD02}
\end{eqnarray}
are satisfied 
on ${\mathcal U}_{Ric(T)} \cap {\mathcal U}_{Weyl(T)} 
\subset M$, 
then on this set we have 
\begin{eqnarray}
\kappa(T) &=& 0 ,
\label{WDWD03}\\
%\end{eqnarray}
%\begin{eqnarray}
T \cdot Weyl(T) \ =\  Weyl (T) \cdot T &=&  Q(Ric(T), Weyl(T)) \ =\ 0 .
%\label{WDWD04}
\nonumber
\end{eqnarray}
(ii) 
If the conditions (\ref{WDWD02}) and
\begin{eqnarray*}
Q(Ric(T),Weyl(T) ) &=& 0 
%\label{WDWD61}
\end{eqnarray*}
are satisfied 
on ${\mathcal U}_{Ric(T)} \cap {\mathcal U}_{Weyl(T)} 
\subset M$, 
then on this set we have (\ref{WDWD03}) and 
\begin{eqnarray*}
T \cdot Weyl(T) &=&  Weyl (T) \cdot T  \ =\ 0 .
%\label{WDWD62}
\end{eqnarray*}
\end{thm}
In {\cite[Section 3] {DGHZ01}} an example of warped product manifolds 
satisfying assumptions of {\cite[Theorem 3.4] {DGHZ01}} is also given.

Let $(M,g)$, $\dim M = n \geq 4$, be a semi-Riemannian manifold 
with parallel Weyl conformal curvature tensor,
i.e. $\nabla C  = 0$ on $M$. 
It is obvious that the last condition implies $R \cdot C  = 0$.
Moreover, let the manifold $(M,g)$ be neither conformally flat nor 
locally symmetric. 
Such manifolds are called {\sl{essentially conformally 
symmetric manifolds}}, e.c.s. manifolds/metrics, or ECS manifolds/metrics, 
in short (see, e.g.,
\cite{{DerRot01}, {DerRot02}, {DerRot2007PADGE}, 
{DerTer2210a}, {DerTer2210b}}).
E.c.s. manifolds 
are semisymmetric manifolds
($R \cdot R  = 0$ {\cite[Theorem 9] {DerRot01}})
satisfying
$\kappa  = 0$ and $Q(S,C) \, =\, 0$
({\cite[Theorems 7 and 8] {DerRot01}}).
In addition, 
\begin{eqnarray}
F\, C &=& \frac{1}{2}\, S \wedge S 
\label{roter77}
\end{eqnarray}
holds on $M$, where $F$ is some function on $M$, called the 
{\sl{fundamental function}} \cite{DerRot02}.  
At every point of $M$ we also have $\mbox{rank}\, S \leq 2$ 
{\cite[Theorem 5] {DerRot02}}.
We mention that the local structure 
of e.c.s. manifolds is already determined.
We refer to \cite{{DerRot2007}, {DerRot2009}} for results 
related to this subject.
We also mention that certain e.c.s. metrics are realized 
on compact manifolds 
\cite{{DerRot2008}, {DerRot2010}, {DerTer2210a}, {DerTer2210b}, {DerTer2301a}, {DerTer2304b},  {DerTer2306a}}.

Equation (\ref{roter77}), by suitable contraction, leads immediately to
$S^{2} = \kappa \, S$, which by $\kappa = 0$, reduces to $S^{2} = 0$. 
Evidently, $\mathrm{tr}_{g} (S^{2}) = 0$.
Now using (\ref{roter77}) we get 
\begin{eqnarray*}
F\, C &=& \frac{n-2}{2 (n -2)} \, S \wedge S
\ =\ \frac{1}{n-2} 
\left(  
g \wedge S^{2} + \frac{n-2}{2} \, S \wedge S - \kappa \, g \wedge S
+ \frac{\kappa ^{2} - \mathrm{tr} (S^{2})}{2(n-1)} \, g \wedge g 
\right) .
\end{eqnarray*}
Thus we have
\begin{thm} {\cite[Theorem 6.1] {2023_DGHP-TZ}}
Condition 
(\ref{2022.11.22.jj}), with $\tau = (n-2) F$,
is satisfied on every essentially conformally symmetric manifold $(M,g)$.
\end{thm}

We assume that $F  = 0$ at $x \in M$. Now (\ref{roter77}) implies 
$\mbox{rank}\, S \leq 1$ at $x$.
It is clear that if $S$ vanishes then
\begin{eqnarray}
R \cdot C &=& C \cdot R \ = \ Q(S,C) \ = \ 0 
\label{2020.09.20.a}
\end{eqnarray}
holds at $x$.
If $\mbox{rank}\, S = 1$ then in view 
of
{\cite[Proposition 3.9(ii)] {DGHZ01}}
we also have 
(\ref{2020.09.20.a}) 
at $x$.
Next, we assume that $F$ is non-zero at $x \in M$. 
Thus $\mbox{rank}\, S  =  2$ at $x$. 
Now
(\ref{roter77}) turns into (\ref{eq:h7}) with 
$T = R$, $Ric(T) = S$, 
$\phi  = F^{-1}$, $\mu =  \frac{1}{n-2}$ and $\eta  = 0$.
Therefore (\ref{roter73NN}) and (\ref{roter73}) reduce to 
$C \cdot R  = 0$ and
$C \cdot C  = 0$, respectively. Consequently,  
(\ref{2020.09.20.a}) 
holds at $x$.
Thus we have  
\begin{thm}
{\cite[Theorem 4.1] {DGHZ01}}
Condition 
(\ref{2020.09.20.a}) 
is satisfied on every essentially conformally symmetric manifold $(M,g)$. 
\end{thm}
Thus we see that e.c.s. manifolds satisfy (\ref{advances4}). 
We also mention that the tensor $C \cdot C$ of every e.c.s. manifold
is the zero tensor ({\cite[Remark 4.2(ii)] {DGHZ01}}).  

E.c.s. warped product manifolds were investigated in \cite{Hotlos}.
In that paper examples of such manifolds 
are given {\cite[Remark 4.2(i)] {DGHZ01}}.

In 
{\cite[Section 5] {DGHZ01}}  
Roter type manifolds satisfying (\ref{advances4}) are investigated. 
In 
{\cite[Theorem 5.2] {DGHZ01}}  
it was stated
that if $(M,g)$, $\dim M = n \geq 4$, is a Roter type manifold
with vanishing scalar curvature $\kappa$ on ${\mathcal U} \subset M$
then (\ref{advances4}), with $L = - 1$, holds on this set.

\begin{thm}
{\cite[Theorem 5.2] {DGHZ01}}  
Let $(M,g)$, $\dim M = n \geq 4$, be a a semi-Riemannian 
manifold satisfying 
(\ref{eq:h7a}) 
on ${\mathcal U}_{S} \cap {\mathcal U}_{C} \subset M$. 
If $\kappa = 0$ on 
${\mathcal U}_{S} \cap {\mathcal U}_{C}$  then 
(\ref{advances4}), with $L\, =\, - 1$, holds on this set.
\end{thm}
This result is also an immediate consequence of the fact that
every Roter space $(M,g)$ satisfies 
(\ref{Roterformula})  on 
${\mathcal U}_{S} \cap {\mathcal U}_{C} \subset M$  (see Section 4).
However in the proof of the last theorem formula (\ref{Roterformula})
was not applied.

We also have 
\begin{thm}
{\cite[Theorem 5.3] {DGHZ01}}  
Let $(M,g)$, $\dim M = n \geq 4$, be a semi-Riemannian manifold
satisfying 
(\ref{eq:h7a}) and (\ref{advances4}) 
on ${\mathcal U}_{S} \cap {\mathcal U}_{C} \subset M$,
and let ${\mathcal U}_{1} \subset {\mathcal U}_{S} \cap {\mathcal U}_{C}$ 
be the set of all points at which
the functions $L$ and $L_{C}$, defined by 
(\ref{advances4}), (\ref{roter73NN}) and (\ref{roter73}) 
(for $T  = R$), 
respectively, 
are nowhere zero on this set. Then we have on ${\mathcal U}_{1}$: 
$L = - 1$ and $\kappa = 0$. 
\end{thm}

In {\cite[Example 5.4] {DGHZ01}} it was shown that 
under some conditions the Cartesian product 
of two semi-Riemannian spaces of constant curvature satisfies assumptions 
of Theorem 6.6 (i.e., {\cite[Theorem 5.3] {DGHZ01}}).

\section{2-quasi-umbilical hypersurfaces}

Let $N_{s}^{n+1}(c)$, $n \geq 3$, 
be a semi-Riemannian space of constant curvature 
$c = \frac{\widetilde{\kappa}}{n (n+1)}$ 
with signature $(s,n+1-s)$, 
where $\widetilde{\kappa}$ is its scalar curvature. 
Let $M$ be a connected hypersurface
isometrically immersed in $N_{s}^{n+1}(c)$.
We have (see, e.g., \cite{R99})
\begin{eqnarray}
R_{hijk} &=& \varepsilon \, ( H_{hk} H_{ij} - H_{hj} H_{ik})
+ \frac{\widetilde{\kappa} }{n(n+1)}\, G_{hijk}\, ,\ \ \ 
\varepsilon \ =\ \pm 1 ,
\label{C5}
\end{eqnarray}
where 
$R_{hijk}$, $G_{hijk}$ and $H_{hk}$, respectively, are
the local components of the curvature tensor $R$ of $M$, the tensor $G$ and 
the second fundamental tensor $H$, respectively.
Contracting (\ref{C5}) with $g^{ij}$ and $g^{kh}$, respectively, we obtain
\begin{eqnarray*}
S_{hk} &=& \varepsilon \, ( \mathrm{tr}(H) \, H_{hk} - H^{2}_{hk})
+ \frac{(n-1) \widetilde{\kappa} }{n(n+1)}\, g_{hk} ,
%\label{C55}
\\
\frac{\kappa }{n-1} &=& \frac{ \varepsilon }{n-1}\, 
( (\mathrm{tr}(H))^{2} - \mathrm{tr}(H^{2})) 
+ \frac{ \widetilde{\kappa}}{n+1}  ,
%\label{C55a}
\end{eqnarray*}
respectively, where
$\mathrm{tr}(H)  = g^{hk}H_{hk}$,
$\mathrm{tr}(H^{2})  =  g^{hk}H^{2}_{hk}$
and $S_{hk}$ are the local components of the Ricci tensor $S$ of $M$
and $\kappa $ is the scalar curvature of $M$.
It is known that (\ref{genpseudo01}) holds on $M$. Precisely, 
\begin{eqnarray}
R \cdot R - Q(S,R) &=&  - \frac{(n-2) \widetilde{\kappa} }{n(n+1)}\, Q(g,C) 
\label{900ab}
\end{eqnarray}
on $M$ {\cite[Lemma 2.1] {DV1993}} (see also {\cite[Lemma 2.1 (ii)] {27}}). 
Evidently, if $n=3$ then 
(\ref{900ab}) turns into
\begin{eqnarray}
R \cdot R &=&  Q(S,R) .
\label{900abeucl}
\end{eqnarray}
If the ambient space is a semi-Euclidean space $\mathbb{E}^{n+1}_{s}$, 
$n \geq 3$, then (\ref{900ab}) also reduces to (\ref{900abeucl}).
As it was proved in {\cite[Lemma 4.1] {47}},
we also have the following identity on $M$ in $N_{s}^{n+1}(c)$, $n \geq 3$,
\begin{eqnarray}
R \cdot R - \frac{ \widetilde{\kappa} }{n(n+1)}\, Q(g,R)
&=& - Q(H^{2}, \frac{1}{2}\, H \wedge H ).
\label{2023.05.05.dd}
\end{eqnarray}

Further we have 
\begin{thm} {\cite[Theorem 3.7] {2016_DGJZ}}
Let $M$ be a hypersurface isometrically immersed in $N_{s}^{n+1}(c)$, 
$n \geq 4$.
Then
\begin{eqnarray*}
C \cdot R  + R \cdot C 
&=& Q(S,C) 
- \frac{(n-2) \widetilde{\kappa} }{n(n+1)}\, Q(g,C) 
+ C \cdot C \nonumber\\
& & - \frac{1}{(n-2)^{2}}\, Q \left( 
g,  g \wedge S^{2} + \frac{n-2}{2}\, S \wedge S 
- \kappa \, g\wedge S  \right) 
%\label{02identity01hyper}
\end{eqnarray*}
holds on $M$.
Moreover, if the condition $C \cdot C = L_{C} Q(g,C)$ 
(i.e., (\ref{4.3.012})) 
is satisfied on ${\mathcal U}_{S} \cap {\mathcal U}_{C} \subset M$ 
then on this set we have
\begin{eqnarray}
C \cdot R + R \cdot C &=& Q(S,C) + \left( L_{C} 
-  \frac{(n-2) \widetilde{\kappa} }{n(n+1)} \right) Q(g,C)\nonumber\\
& & - \frac{1}{(n-2)^{2}}\, 
Q\left(  g, 
 g \wedge S^{2}
+ \frac{n-2}{2}\, S \wedge S - \kappa \, g \wedge S \right) ,
\label{identity05hyper} 
\end{eqnarray}
and in addition, if $M$ is a quasi-Einstein hypersurface satisfying 
the condition $\mathrm{rank}\, (S - \alpha\, g) = 1$
(i.e., (\ref{quasi02})),  
on ${\mathcal U}_{S} \cap {\mathcal U}_{C}$ then on this set we have 
\begin{eqnarray*}
C \cdot R + R \cdot C  &=& Q(S,C) 
+ \left( L_{C} -  \frac{(n-2) \widetilde{\kappa} }{n(n+1)} \right) Q(g,C) .
%\label{identity05quasihypera}
\end{eqnarray*}
\end{thm}

It is known that every $2$-quasi-umbilical hypersurface 
in a semi-Riemannian space of constant curvature $N_{s}^{n+1}(c)$, 
$n \geq 4$, 
satisfies (\ref{4.3.012}) (see, e.g., {\cite[Theorem 3.8] {2016_DGJZ}}). 
We have
\begin{thm} {\cite[Theorem 3.8] {2016_DGJZ}}
If $M$ is a $2$-quasi-umbilical hypersurface isometrically immersed 
in $N_{s}^{n+1}(c)$, $n \geq 4$,
then (\ref{identity05hyper}) holds on 
${\mathcal U}_{S} \cap {\mathcal U}_{C} \subset M$.
\end{thm}

\begin{thm}
{\cite[Proposition 4.2] {DGP-TV02}}  
Let $M$ be a hypersurface in an Euclidean space $\mathbb{E}^{n+1}$, 
$n \geq 4$, 
having exactly three distinct principal curvatures 
$\lambda _{1}$, $\lambda _{2}$ and $\lambda _{3}$ 
satisfying at every point of $M$: $\lambda _{1} = 0$, 
$\lambda _{2} = - (n-2) \lambda $ 
and  
$\lambda _{3} = \lambda _{4} = \ldots = \lambda _{n} = \lambda \neq 0$.
Then $M$ is a minimal, $2$-quasi-umbilical  
and $2$-quasi-Einstein hypersurface
satisfying: (\ref{dgss01ab}), (\ref{900abeucl}) and
\begin{eqnarray*}
\mathrm{rank}\, \left( S  -  \frac{\kappa }{(n-2)(n-1)}\, g \right)  
&=& 2 ,
\end{eqnarray*}
\begin{eqnarray*}
S & = & - H^{2} ,\ \ \ \kappa \ =\ - tr(H^{2})
\ =\ - (n-2)(n-1) \lambda ^{2} ,\nonumber\\
 S^{2} &=&
- ( \phi ^{2} + \psi  ) \, S + \phi \psi \, H     ,
\nonumber\\
S^{3} 
&=&
- ( \phi ^{2} + 2 \psi )\, S^{2} 
- \psi ^{2} S ,
\end{eqnarray*}
\begin{eqnarray}
R &=&
\frac{1}{2  ( \phi \psi )^{2} } 
\left( S^{2} + (  \phi ^{2} + \psi   ) S \right)
\wedge
\left( S^{2} + (  \phi ^{2} + \psi   ) S \right)  ,
\label{B001simplyhypersurf}\\
%\end{eqnarray*}
%\begin{eqnarray*}
R \cdot S 
&=& \phi \, Q(H, S) \ =\ \frac{n-1}{\kappa }\, Q(S,S^{2}) ,\nonumber\\ 
C \cdot S &=& 
 \phi \, Q(H,S)
+ \frac{\phi ^{2}}{ n-2} \, Q(g,S)
-  \frac{\phi \psi }{ n-2} \, Q(g,H) \nonumber\\
&=&
\frac{\kappa }{(n-2) (n-1) }\, Q(g,S)
- \frac{  1 }{ n-2 }\, Q(g,S^{2} )
+ \frac{n-1}{ \kappa }\, Q(S,S^{2})\nonumber ,
\end{eqnarray}
\begin{eqnarray*}
(n-2)\, R \cdot C &=& (n-2)\, Q(S,R) - \phi \, g \wedge Q(H,S) , \\
(n-2)\, C \cdot R &=&
(n-3)\, Q(S,R) - \phi \, H \wedge Q(g,S) ,\\
C \cdot C & = & 0 ,
%\label{ffnewZZ2ff}
\end{eqnarray*}
where $\phi  = - (n-3) \lambda $ and  $\psi  = (n-2) \lambda ^{2}$. 
Moreover, we have
\begin{eqnarray}
(n-2) (R \cdot C - C \cdot R) &=& Q(S,R) 
+ \phi \, ( H \wedge Q(g,S) - g \wedge Q(H,S)).
\label{2023.05.20.a}
\end{eqnarray}
\end{thm}
We note that (\ref{B001simplyhypersurf})
is a particular form of the Roter type equation (\ref{B001simply}).

Biharmonic hypersurfaces with three distinct principal curvatures 
in an Euclidean $5$-space $\mathbb{E}^{5}$ were investigated in \cite{Fu}.
The main result of {\cite[Theorem 3.2] {Fu}} states that  
every biharmonic hypersurface $M$ with three distinct 
principal curvatures in $\mathbb{E}^{5}$ is minimal.
The principal curvatures of $M$ are the following:
$\lambda _{1} = 0$, $\lambda _{2} = -2 \lambda $ 
and $\lambda _{3} = \lambda _{4} = \lambda \neq 0$, 
where $\lambda$ is some function on $M$.
Curvature properties of such hypersurfaces 
are expressed in the last theorem,
provided that $n = 4$ 
(see also {\cite[Theorem 4.3] {DGP-TV02}}). 

We refer to \cite{OU-CHEN-2020} for a survey of results on biharmonic hypersurfaces.

Let now $M$ be a type number two hypersurface in $N_{s}^{n+1}(c)$, $n \geq 3$, 
i.e.,
let $\mathrm{rank} (H) = 2$ at every point of $M$ 
(see, e.g., \cite{CHEN-YILDIRIM-2015}). 
Using  
(\ref{2023.05.05.aa}), (\ref{2023.05.05.bb}), (\ref{2023.05.05.cc})
and (\ref{2023.05.05.dd}) we obtain on $M$ (cf. {\cite[Theoem 4.2] {47}})
\begin{eqnarray}
& &
H^{3} \ =\ \mathrm{tr} (H)\, H^{2} 
+ \frac{ \mathrm{tr} ( H^{2}) - (\mathrm{tr} (H))^{2} }{2} \, H ,
\label{2023.05.05.ee}\\
& &
Q(H, H \wedge H^{2}) \ =\ 0 ,
\label{2023.05.05.ff}\\
& &
R \cdot R \ =\ \frac{ \widetilde{\kappa} }{n(n+1)}\, Q(g,R) .
\label{2023.05.05.gg}
\end{eqnarray}
Thus we have
\begin{thm} (cf. {\cite[Theorem 4.2] {47}})
Every type number two hypersurface $M$ in  $N_{s}^{n+1}(c)$, $n \geq 3$,
is a pseudosymmetric manifold of constant type satisfying 
(\ref{2023.05.05.ee}), (\ref{2023.05.05.ff}) and (\ref{2023.05.05.gg}).
\end{thm}

We mention that type number two hypersurfaces in $4$-dimensional 
Riemannian space of constant curvature were investigated in \cite{HajKowSek}. 
We also note that in {\cite[Section 5] {44}} it was stated that the Cartan
hypersurface in a $4$-dimensional sphere is a pseudosymmetric manifold of constant
type. 
Evidently, (\ref{2023.05.05.ee}) is a special form of the equation
\begin{eqnarray}
H^{3} &=& \mathrm{tr} (H)\, H^{2} + \psi \, H ,
\label{DS4aa}
\end{eqnarray}
where $\psi$ is a function on $M$.
We refer to Remark 8.7 (ii) of this paper for results on 
hypersurfaces in $N_{s}^{n+1}(c)$, $n \geq 4$, satisfying (\ref{DS4aa}).

As an immediate consequence of Proposition 2.1 and 
{\cite[Remark 7.1 (iii)] {2023_DGHP-TZ}} we obtain 
the following 
\begin{thm}
If $M$, $\dim M = n \geq 4$, is a type number two hypersurface 
isometrically immersed 
in a semi-Riemannian conformally flat manifold $N$, $\dim N = n + 1$, then
$C \cdot C = L_{C} Q(g,C)$ (i.e., (\ref{4.3.012})) holds on 
${\mathcal U}_{C} \subset M$, where $L_{C}$ is some function on this set.
\end{thm} 

Let $M$ be a type number two hypersurface in  $N_{s}^{n+1}(c)$, $n \geq 4$.
Thus, in view of the last theorem,
(\ref{4.3.012}) holds on ${\mathcal U}_{C} \subset M$. 
Moreover, in view of Theorem 7.1, (\ref{identity05hyper}) is satisfied 
on ${\mathcal U}_{S} \cap {\mathcal U}_{C} \subset M$.
In addition, we note that (\ref{2023.05.05.gg}) implies
\begin{eqnarray}
R \cdot C \ =\ \frac{ \widetilde{\kappa} }{n(n+1)}\, Q(g,C) .
\label{2023.05.35.aa}
\end{eqnarray}
Now using
Theorem 7.1, Theorem 7.5 and (\ref{2023.05.35.aa}) we get easily 
the following result.
\begin{thm}
Let $M$ be a type number two hypersurface in  $N_{s}^{n+1}(c)$, $n \geq 4$.
Then 
the following condition of the form $(\ast)$ is satisfied 
on ${\mathcal U}_{S} \cap {\mathcal U}_{C}$
\begin{eqnarray}
\ \ \ 
C \cdot R - R \cdot C &=& Q(S,C) + L_{1}\, Q(g,C) 
- \frac{1}{(n-2)^{2}}\, 
Q\left(  g, 
 g \wedge S^{2}
+ \frac{n-2}{2}\, S \wedge S - \kappa \, g \wedge S \right) ,
\label{2023.05.35.bb}
\end{eqnarray}
where $L_{1}$ is some function on this set. Moreover, 
if $M$ is a quasi-Einstein hypersurface satisfying 
the condition $\mathrm{rank}\, (S - \alpha\, g) = 1$
(i.e., (\ref{quasi02}))  
on ${\mathcal U}_{S} \cap {\mathcal U}_{C}$ then on this set we have 
\begin{eqnarray}
\ \ \ 
C \cdot R - R \cdot C &=& Q(S,C) + L_{1}\, Q(g,C) .
\label{2023.05.35.cc}
\end{eqnarray}
\end{thm}

We refer to {\cite[Section 5] {R102}} for further results 
on type number two hypersurfaces in $N_{s}^{n+1}(c)$, $n \geq 4$,
satisfying curvatures conditions of pseudosymmetry type.

\section{The condition $H^{3} = \mathrm{tr} (H)\, H^{2} + \psi \, H + \rho \, g$}

Let $M$ be a hypersurface isometrically immersed in  
$N_{s}^{n+1}(c)$, $n \geq 4$,
satisfying 
on ${\mathcal U}_{H} \subset M$
a curvature condition of the kind:
the tensor $R \cdot C$, $C \cdot R$ or 
$R \cdot C - R \cdot C$ is a linear combination
of the tensor $R \cdot R$ and of a finite sum 
of the Tachibana tensors of the form $Q(A,T)$, 
where $A$ is a symmetric $(0,2)$-tensor 
and $T$ a generalized curvature tensor.
As it was mentioned in Introduction, 
if such condition is satisfied on ${\mathcal U}_{H}$ then 
(\ref{dgss01}) holds on this set. 
We present now the following results on hypersurfaces $M$
in $N_{s}^{n+1}(c)$, $n \geq 4$,
satisfying (\ref{dgss01})
on ${\mathcal U}_{H} \subset M$. 

\begin{thm} {\cite[Proposition 5.1, eq. (29)]{Saw5}}
If $M$ is a hypersurface in a semi-Rieman\-nian space of constant curvature 
$N_{s}^{n+1}(c)$, $n \geq 4$,
satisfying (\ref{dgss01}) on ${\mathcal U}_{H} \subset M$,
for some functions 
$\psi $ and $\rho $ on ${\mathcal U}_{H}$,
then on ${\mathcal U}_{H}$ we have
\begin{eqnarray*}
R \cdot C 
&=& 
Q(S,R) 
- \frac{(n-2) \widetilde{\kappa }}{n (n+1)}\, Q(g,R)
+ \alpha _{2}\, Q(S,G)
+ \frac{\rho }{n-2}\, Q(H,G) ,
%\label{ZZ1}
\\
%\end{eqnarray}
%\begin{eqnarray}
C \cdot R &=&
\frac{n-3}{n-2}\, Q(S,R) 
+ \alpha _{1}\, Q(g,R)
+ \alpha _{2}\,  Q(S,G) ,
%\label{ZZ2}
\end{eqnarray*}
\begin{eqnarray}
(n-2)\, (R \cdot C - C \cdot R) &=& Q(S,R) + \rho \, Q(H,G)
 + \left( \frac{(n-1) \widetilde{\kappa} }{n(n+1)} - \frac{\kappa }{n-1} 
- \varepsilon \psi \right) Q(g, R),
\label{ZZ3}\\
%\end{eqnarray}
%\begin{eqnarray}
(n-2)\, C \cdot C &=& 
(n-3)\, Q(S,R) + (n-2) \alpha _{1}\, Q(g,R)\nonumber\\
& &
+ ( \alpha _{1} - \alpha _{2} )\, Q(S,G)
+ \frac{n-3}{n-2} \rho \, Q(H,G) ,
%\label{DS16A}
\nonumber
\end{eqnarray}
\begin{eqnarray}
R \cdot S &=& 
\frac{\widetilde{\kappa }}{n(n+1)}\, Q(g,S) + \rho \, Q(g,H) ,
%\label{DZ004}
\nonumber\\
%\end{eqnarray}
%\begin{eqnarray}
\alpha _{1} 
&=&
\frac{1}{n-2} \left( \frac{\kappa }{n-1} + \varepsilon \psi
- \frac{(n^{2} - 3n + 3) \widetilde{\kappa }}{n(n+1)} \right) ,
\label{AA00two}\\
\alpha _{2} 
&=&
- \frac{(n-3) \widetilde{\kappa} }{(n-2)n (n+1)} .
\label{AA00}
\end{eqnarray}
\end{thm}

\begin{thm} {\cite[Proposition 4.2]{2020_DGZ}}
Let $M$ be a hypersurface in 
a semi-Riemannian space of constant curvature 
$N_{s}^{n+1}(c)$, $n \geq 4$,
satisfying (\ref{dgss01}) on ${\mathcal U}_{H} \subset M$,
for some functions 
$\psi $ and $\rho $ on ${\mathcal U}_{H}$.
\newline
(i) The following conditions are satisfied on ${\mathcal U}_{H}$:
\begin{eqnarray*}
Q(\rho\, H - \alpha_{3}\, S - S^{2},G) &=& 0 ,
%\label{DZ005}
\end{eqnarray*}
\begin{eqnarray}
\rho\, H  &=& S^{2} + \alpha_{3}\, S  + \frac{\lambda}{n}\, g ,\ \ \ 
\lambda \ =\ 
\rho\, \mathrm{tr}(H) - \kappa \, \alpha_{3} - \mathrm{tr}(S^{2}) , 
\label{DZ008}\\
%\end{eqnarray}
%\begin{eqnarray}
\alpha _{3}
&=&
(n-2)^{2} 
\left( \frac{1}{n-2} (\alpha_{1} - \alpha_{2}) - 2\alpha_{2} 
- \frac{ \widetilde{\kappa} }{n(n+1)} \right) 
- \frac{\kappa }{n-1} 
\ =\
\varepsilon \psi - \frac{ 2 ( n-1) \widetilde{\kappa} }{n(n+1)} , 
\label{DZ006}
\end{eqnarray}
where 
$\alpha _{1}$ and $\alpha _{2}$ are defined by (\ref{AA00two}) 
and (\ref{AA00}), respectively.
Moreover,
\begin{eqnarray*}
R \cdot S &=&  
Q(g, S^{2}) + \left( \varepsilon \psi 
- \frac{( 2 n - 3) \widetilde{\kappa }}{n(n+1)} \right) Q(g,S ) ,
%\label{GGG01}
\\
R \cdot S^{2} &=& Q(S, S^{2}) 
+ \rho _{1}\, Q(g, S^{2}) + \rho _{2}\, Q(g, S ) ,
%\label{EEE01}
\\
S^{3} &=& \left( - 2 \varepsilon \psi 
+ \frac{ 3 (n - 1) \widetilde{\kappa }}{n(n+1)} \right) S^{2} 
+ \rho_{2} \, S + \rho _{3}\, g ,
%\label{EEE01new}
\end{eqnarray*}
hold on ${\mathcal U}_{H}$, 
where the functions $\rho _{1}$, $\rho _{2}$ and $\rho _{3}$ are defined by 
\begin{eqnarray*}
\rho _{1} &=& - \frac{ (n-2) \widetilde{\kappa} }{n (n+1) } 
- \alpha _{3} ,\nonumber\\ 
\rho _{2} &=& - \frac{\lambda}{n} 
- \left( \frac{ (n-1) \widetilde{\kappa} }{n (n+1) } 
+ \alpha _{3}\right)  \alpha _{3} ,\nonumber\\
\rho _{3} &=&
\frac{1}{n} 
\left( \mathrm{tr}(S^{3}) + \left( 2 \varepsilon \psi 
- \frac{ 3 (n-1) \widetilde{\kappa} }{n (n+1) } \right)
\mathrm{tr}(S^{2}) - \kappa \rho _{2} \right) .
%\label{EEE01bb}
\end{eqnarray*}
(ii) If at a point $x \in {\mathcal U}_{H}$ we have
$S^{2} = \beta_{1}\, S + \beta _{2}\, g$, 
for some $\beta_{1}, \, \beta_{2} \in {\mathbb{R}}$, 
then $\rho = 0$, $\beta _{1} = \alpha _{3}$ and $\beta _{2} 
= - ( \lambda / n)$ at this point. 
\end{thm}

Theorems 8.1 and 8.2 (precisely, (\ref{ZZ3}), (\ref{DZ008}) and (\ref{DZ006}))
lead to the following result. 
\begin{thm} 
If $M$ is a hypersurface in 
a semi-Riemannian space of constant curvature 
$N_{s}^{n+1}(c)$, $n \geq 4$,
satisfying (\ref{dgss01}) on ${\mathcal{U}}_{H} \subset M$,
for some functions 
$\psi $ and $\rho $ on ${\mathcal{U}}_{H}$
then 
\begin{eqnarray*}
(n-2)\, (R \cdot C - C \cdot R) 
&=& Q(S,R) 
 + \left( \frac{(n-1) \widetilde{\kappa} }{n(n+1)} - \frac{\kappa }{n-1} 
- \varepsilon \psi \right) Q(g, R)\nonumber\\
& &
+ Q(S^{2},G)
+ \left( \varepsilon \psi 
- \frac{ 2 ( n-1) \widetilde{\kappa} }{n(n+1)} \right) Q(S,G) 
%\label{ZZ3mod}
\end{eqnarray*}
on ${\mathcal{U}}_{H}$.
\end{thm}

\begin{thm} {\cite[Theorem 4.5]{2020_DGZ}}
Let $M$ be a hypersurface in 
a semi-Riemannian space of constant curvature 
$N_{s}^{n+1}(c)$, $n \geq 4$,
satisfying (\ref{dgss01}) on ${\mathcal{U}}_{H} \subset M$,
for some functions 
$\psi $ and $\rho $ on ${\mathcal{U}}_{H}$.
If the tensor $C \cdot C$ and a generalized curvature tensor $T$ satisfy 
$C \cdot C = Q(g,T)$
on ${\mathcal{U}}_{H}$, then on this set we have
\begin{eqnarray*}
T &=&
\left( \frac{ \kappa }{n-1} + \frac{ 2 \varepsilon \psi}{n-1} 
- \frac{  \widetilde{\kappa} }{n+1} \right) C 
 + \lambda \, G
- \frac{n-3}{(n-2)^{2} (n-1)} 
\left(
g \wedge S^{2} + \frac{n-2}{2}\, S \wedge S - \kappa\, g \wedge S 
\right) ,
%\label{dgpss777ccdd}
\end{eqnarray*}
where $\lambda$ is some function on ${\mathcal{U}}_{H}$.
\end{thm}

Further, we also have the following results obtained 
in {\cite[Section 4]{2020_DGZ}}.
\begin{thm}
Let $M$ be a hypersurface isometrically immersed in  
$N_{s}^{n+1}(c)$, $n \geq 4$,
satisfying (\ref{dgss01}) 
on ${\mathcal U}_{H} \subset M$.
\newline
(i) {\cite[Theorem 4.6]{2020_DGZ}} 
If on ${\mathcal U}_{H}$
the tensor $Q(S,R)$ is equal to the Tachibana tensor $Q(g,T_{1})$, 
where $T_{1}$ is a generalized curvature tensor, then any of the tensors:  
$R \cdot R$, $R \cdot C$, $C \cdot R$, 
$R \cdot C - C \cdot R$ and $C \cdot C$
is equal to some Tachibana tensor $Q(g,T_{2})$, 
where $T_{2}$ is a linear combination of the tensors
$R$, $g \wedge g$, $g \wedge S$, $g \wedge S^{2}$ and $S \wedge S$.
\newline
(ii) {\cite[Theorem 4.7]{2020_DGZ}} 
The following conditions 
are satisfied on ${\mathcal U}_{H}$
\begin{eqnarray*}
C \cdot C 
&=&
\frac{n-3}{n-2}\, R \cdot C
+ \frac{1}{n-2} \left( \frac{\kappa }{n-1} + \varepsilon \psi
- \frac{( 2 n - 3) \widetilde{\kappa }}{n(n+1)} \right)
Q(g,C),
%\label{DS16Anew01}
\\
\nonumber\\
(n-2)\, C \cdot R  + R \cdot C 
&=& (n-2)\, Q(S,C) 
+ \left(
\frac{\kappa}{n-1} + \varepsilon \psi 
- \frac{ (n-1)^{2} \widetilde{\kappa}}{n (n+1)} 
\right)
Q(g,C)\nonumber\\
& & - \frac{1}{(n-2)}\, Q\left( g, 
g \wedge S^{2}
+ \frac{n-2}{2}\, S \wedge S 
- \kappa \, g\wedge S  \right) .
%\label{02identity01hyper17}
\end{eqnarray*}
\end{thm}
\begin{thm}
{\cite[Theorem 4.8]{2020_DGZ}}
Let $M$ be a hypersurface isometrically immersed in  
$N_{s}^{n+1}(c)$, $n \geq 4$, such that 
(\ref{quasi02}) and (\ref{DS4aa})
are satisfied on ${\mathcal U}_{H} \subset M$, where $\psi$ is a function on 
this set. 
Then (\ref{Roterformula}) holds on ${\mathcal U}_{H}$ 
if and only if the following two conditions hold
on this set 
\begin{eqnarray}
(a)\ \  \frac{\kappa }{n-1} \ =\ \frac{\widetilde{ \kappa }}{n+1}
\ \ \ &\mbox{and}& \ \ \ 
(b)\ \  Q\left( S - \frac{\kappa }{n}\, g, C \right)\ =\ 0 .
\label{qqee08}
\end{eqnarray}
\end{thm}
For further results related to Theorem 8.5 (i) 
we refer to {\cite[Section 8] {2023_DGHP-TZ}}. 

\begin{rem}
(i) 
Let $M$ be a hypersurface in $N_{s}^{n+1}(c)$, $n \geq 4$.
As it was proved in {\cite[Proposition 3.1 (ii), Proposition 3.2] {47}}
the conditions (\ref{DS4aa}) and 
\begin{eqnarray}
R \cdot S &=& \frac{\widetilde{\kappa}}{n (n+1)}\, Q(g,S)
\label{900r22}
\end{eqnarray}
are equivalent on the set ${\mathcal U}_{H} \subset M$. 
If (\ref{DS4aa}) holds on ${\mathcal U}_{H} \subset M$ then on this set we have
{\cite[Theorem 3.1)] {DG90}} 
\begin{eqnarray}
(n-2) (R \cdot C - C \cdot R) 
&=& Q(S,R) + \left( \frac{(n-1) \widetilde{\kappa} }{n (n+1)}
- \frac{\kappa }{n-1} - \varepsilon \psi \right) Q(g,R) .
\label{900r77}
\end{eqnarray}
(ii) 
Let now $M$ be a hypersurface in a Riemannian space of constant curvature
$N^{n+1}(c)$, $n \geq 4$, 
such that at every point of $M$ there are principal curvatures
$0, \ldots , 0, \lambda, \ldots , \lambda, - \lambda, \ldots , - \lambda$,
with the same multiplicity of $\lambda $ and $- \lambda$, 
and $\lambda $ is a positive function on $M$. 
Thus we have on $M$:
$\mathrm{tr} (H) = 0$ and $H^{3}  = \lambda ^{2} \, H$,   
and, in a consequence, we also  have (\ref{900r22}) and (\ref{900r77}).
We mention that the Cartan hypersurfaces, as well as 
generalized Cartan hypersurfaces {\cite[Section 6] {BYCH1996}}
have at every point principal curvatures
$0, \ldots , 0, \lambda, \ldots , \lambda, - \lambda, \ldots , - \lambda$,
with the same multiplicity of $\lambda $ and $- \lambda$, 
and $\lambda $ is a positive function on $M$.
Curvature properties of pseudosymmetry type of Cartan hypersurfaces
(see, e.g., {\cite[Chapter 3] {TEC_PJR_2015}}, 
{\cite[Chapter 7.5] {OU-CHEN-2020}}) 
are given in \cite{{DG90}, {44}, {DY}}. 
Both classes of the considered hypersurfaces are austere 
hypersurfaces \cite{Bryant-1991}.
\end{rem}

\section{Hypersurfaces satisfying the condition $R \cdot C - C \cdot R = L Q(S,C)$}

In \cite{DGHSaw} a survey on some family of generalized Einstein 
metric conditions was given. 
Those curvature conditions are strongly related to pseudosymmetry. 
In particular, {\cite[Section 6] {DGHSaw}} contains results  
of non-Einstein and non-conformally flat semi-Riemannian manifolds 
$(M,g)$, of dimension $n \geq 4$,
satisfying conditions of the form: the tensor 
$R\cdot C - C \cdot  R$ is proportional 
to the Tachibana tensor: $Q(g,R)$, $Q(S,R)$, $Q(g,C)$ or $Q(S,C)$.
More precisely, those conditions are considered on 
${\mathcal U}_{S} \cap {\mathcal U}_{C} \subset M$. 
Among other results in that section  
it was shown that 
some hypersurfaces $M$ isometrically immersed 
in a semi-Riemannian space of constant curvature
$N^{n+1}_{s}(c)$, $n \geq 4$, satisfy 
(\ref{advances4})
on ${\mathcal U}_{S} \cap {\mathcal U}_{C} \subset M$.
We recall that an example of a hypersurface having mentioned properties
was constructed in {\cite[Section 5] {R102}}. 
We also mention that semi-Riemannian manifolds satisfying (\ref{advances4}) 
were investigated in \cite{DGHZ01}. 

Let $M$ be a hypersurface isometrically immersed 
in $N^{n+1}_{s}(c)$, $n \geq 4$,
satisfying (\ref{advances4}) on 
$({\mathcal U}_{S} \cap {\mathcal U}_{C}) \setminus {\mathcal U}_{H}$. 
We recall that (\ref{h2277aa}) holds 
on $({\mathcal U}_{S} \cap {\mathcal U}_{C}) \setminus {\mathcal U}_{H}$, 
where $\alpha $ and $\beta $ are some functions defined on this set. 

According to {\cite[Proposition 3.3] {G108}}, 
the Riemann-Christoffel curvature tensor $R$ of $M$
is expressed
on 
$({\mathcal U}_{S} \cap {\mathcal U}_{C}) \setminus {\mathcal U}_{H}$ 
by a linear combination of 
the Kulkarni-Nomizu products $S \wedge S$, $g \wedge S$ and 
$G = \frac{1}{2}\, g \wedge g$ 
formed by the Ricci tensor $S$ and the metric tensor $g$ of $M$. 
Precisely, we have  (\ref{eq:h7a}), i.e.,
\begin{eqnarray*}
R &=& \frac{\phi }{2}\, S \wedge S + \mu \, g \wedge S + \eta \, G ,
%\label{eq:h7a}
\end{eqnarray*}
where $\phi$, $\mu $ and $\eta $ are some functions on 
$({\mathcal U}_{S} \cap {\mathcal U}_{C}) \setminus {\mathcal U}_{H}$
(see
{\cite[eqs. (5.2)] {2016_DGHZhyper}}).
We also can express the tensors $C \cdot C$, $Q(g,C)$ and $Q(S,C)$
by some linear combinations of the Tachibana tensors formed 
by the tensors $g$ and $H$. 
In 
{\cite[Theorem 5.1] {2016_DGHZhyper}}
it was stated that if the scalar curvature $\kappa $ of 
a hypersurface $M$ in $N^{n+1}_{s}(c)$,  $n \geq 4$,
vanishes on  
$( {\mathcal U}_{S} \cap {\mathcal U}_{C} ) 
\setminus {\mathcal U}_{H} \subset M$ then 
\begin{eqnarray}
R \cdot C - C \cdot R &=& - Q(S,C) 
\label{h2299aa}
\end{eqnarray}
on this set. 
From that theorem it follows immediately 
{\cite[Corollary 5.3] {2016_DGHZhyper}}
that if
$M$ is a hypersurface in a Riemannian space of constant curvature
$N^{n+1}(c)$,  $n \geq 4$, having at every point exactly 
two distinct principal curvatures,
and if its scalar curvature $\kappa $ vanishes on  
$({\mathcal U}_{S} \cap {\mathcal U}_{C}) 
\setminus {\mathcal U}_{H} \subset M$ 
then (\ref{h2299aa}) holds on this set. 
In 
{\cite[Examples 5.4, 5.5 and 5.7] {2016_DGHZhyper}}
examples of non-conformally flat and non-Einstein hypersurfaces,
with $\kappa = 0$, 
having at every point exactly two distinct principal curvatures 
are presented.

As it was mentioned in Introduction,
if at every point of ${\mathcal U}_{H}$ 
of a hypersurface $M$ in $N^{n+1}_{s}(c)$, $n \geq 4$, 
one of the tensors
$R \cdot C$, $C \cdot R$ 
or $R \cdot C - C \cdot R$ is a linear combination of the tensor 
$R \cdot R$ 
and a finite sum of the Tachibana tensors of the form $Q(A,T)$, 
where $A$ is a symmetric $(0,2)$-tensor 
and $T$ a generalized curvature tensor, then 
(\ref{dgss01}) holds on ${\mathcal U}_{H}$ {\cite[Corollary 4.1] {R99}}.
Thus in particular, if (\ref{advances4}) 
is satisfied on ${\mathcal U}_{H}$ then (\ref{dgss01}) holds on this set.
Hypersurfaces in $N^{n+1}_{s}(c)$, $n \geq 4$, satisfying (\ref{dgss01}), 
or in particular
(\ref{dgss01}) with $\rho = 0$, i.e.,  (\ref{DS4aa}),
were studied in several papers:  
\cite{{47}, {DG90}, {DGHS}, {R99}, 
{DGPSS}, {R102}, {DHS02}, {P104}, {DeKow}, 
{DeScher},
{44}, {Glog}, {G6}, 
{Saw4}, {Saw5}, {S3}, {Sawicz}, {Saw114}}. 
Section 6 
of
\cite{2016_DGHZhyper}
contains some results on hypersurfaces satisfying (\ref{dgss01}).
In  Section 7 
of \cite{2016_DGHZhyper}
hypersurfaces $M$ in $N^{n+1}_{s}(c)$, $n \geq 4$,  
satisfying (\ref{advances4}) on ${\mathcal U}_{H} \subset M$
were investigated. 
The main result of this section 
(see, {\cite[Theorem 7.3] {2016_DGHZhyper}})
states that if 
$M$ is a hypersurface in $N^{n+1}_{s}(c)$, $n \geq 4$, 
satisfying on ${\mathcal U}_{H} \subset M$ the equalities: 
(\ref{DS4aa})
and
\begin{eqnarray}
\mathrm{rank}\, ( S - \alpha \, g ) &=& 1 ,
\label{dgss01dd}
\end{eqnarray}
for some function $\alpha$ on $\mathcal{U}_H$, then on this set 
\begin{eqnarray}
& &
\mathrm{rank} \left( S - \left( \frac{\kappa }{n-1} 
- \frac{\widetilde{\kappa} }{n (n+1) } 
\right) g \right) \ = \ 1 ,
\label{dgss01zz}\\
& &
(n-2)\, ( R \cdot C - C \cdot R ) \ =\  Q(S,C) 
- \frac{\widetilde{\kappa} }{n (n+1) } \, Q(g,C) .
\label{dgss01pp}
\end{eqnarray}
In particular, 
if the ambient space is a semi-Euclidean space $\mathbb{E}_{s}^{n+1}$, 
$n \geq 4$,
then
(\ref{dgss01zz}) and (\ref{dgss01pp}) turn into
\begin{eqnarray}
\mathrm{rank} \left( S - \frac{\kappa }{n-1}\, g \right) & = & 1 ,
\label{dgss01zzsimply}\\
(n-2)\, ( R \cdot C - C \cdot R ) &=&  Q(S,C) 
\label{dgss01pvpv}
\end{eqnarray}
respectively.

Let $M$ be a hypersurface in an Euclidean space $\mathbb{E}^{n+1}$, 
$n = 2 p + 1$, $p \geq 2$,
having at every point three principal curvatures
$\lambda_{1} = \lambda \neq 0$,  $\lambda_{2} 
= - \lambda $ and $\lambda_{3} = 0$,
provided that the multiplicity of $\lambda_{1}$, as well as of 
$\lambda_{2}$ is $p$. 
Clearly, $M$ is an austere hypersurface {\cite[p. 102] {HarLaw}}.
Evidently, $M = \mathcal{U}_H$ and 
(\ref{DS4aa}),
(\ref{dgss01zzsimply}) and (\ref{dgss01pvpv}) hold on $M$
(see, 
{\cite[Example 7.5(i)] {2016_DGHZhyper}}
for details). 
We recall that in \cite{AbDi} it was stated that 
$M$ is a non-semisymmetric ($R \cdot R \neq 0$) Ricci-symmetric 
($R \cdot S = 0$) hypersurface.
Further results on hypersurfaces $M$ in $N^{n+1}_{s}(c)$, $n \geq 4$,  
satisfying (\ref{advances4}) on ${\mathcal U}_{H} \subset M$, with 
$L \neq - 1$, are given in 
{\cite[Proposition 7.1] {2016_DGHZhyper}}
and {\cite[Example 7.5(ii)-(iv)] {2016_DGHZhyper}}.
We also have
\begin{thm} 
{\cite[Theorem 7.4] {2016_DGHZhyper}}
If $M$ is a hypersurface in $N^{n+1}_{s}(c)$, $n \geq 4$, 
satisfying (\ref{advances4}) with $L = - 1$
on ${\mathcal U}_{H}$ then on this set we have: $\kappa = 0$ and 
\begin{eqnarray*}
& &
(n-1)\, (R\cdot C - C \cdot  R) \nonumber\\
&=& 
Q\left( g, \left( \frac{ (n-1) \widetilde{\kappa }}{n (n+1)} 
- \varepsilon \psi \right) R
+ 
\left( \frac{ 2 (n-1) \widetilde{\kappa }}{n (n+1)} 
- \varepsilon \psi  \right) g \wedge S
- \frac{1}{2 (n-2)}\, S \wedge S
-  g \wedge S^{2} \right) ,
%\label{advances4kr02}
\end{eqnarray*}
where 
$\widetilde{\kappa }$ is the scalar curvature of $N^{n+1}_{s}(c)$,
$\kappa $ the scalar curvature of $M$ and the function $\psi$
is defined by (\ref{dgss01}).
\end{thm}

In {\cite[Example 7.6] {2016_DGHZhyper}}
it was stated that some tubular hypersurfaces with vanishing 
scalar curvature satisfy (\ref{advances4}) with $L = - 1$.

\section{The condition $R \cdot C - C \cdot R = L_{1}\, Q(S,C) + L_{2}\, Q(g,C)$}
In Section 5 of \cite{2020_DGZ} we consider hypersurfaces 
$M$ in  $N_{s}^{n+1}(c)$, $n \geq 4$, 
satisfying (\ref{cond01})
on ${\mathcal U}_{H} \subset M$.

We recall that a semi-Riemannian manifold $(M,g)$, $n \geq 3$, 
is said to be 
a quasi-Einstein manifold (see Section 4) 
if (\ref{quasi02}) holds
on ${\mathcal U}_{S} \subset M$, where $\alpha $ is some function 
on this set.
Let $M$ be a quasi-Einstein hypersurface in  $N_{s}^{n+1}(c)$, $n \geq 4$, 
satysfying (\ref{cond01}). We have

\begin{thm} {\cite[Theorem 5.1]{2020_DGZ}}
If $M$ be a hypersurface in 
a semi-Riemannian space of constant curvature 
$N_{s}^{n+1}(c)$, $n \geq 4$,
satisfying (\ref{cond01}) and (\ref{quasi02}) 
on ${\mathcal U}_{H} \subset M$,
for some functions 
$\alpha$, $L_{1}$ and $L_{2}$,
then on this set we have  (\ref{DS4aa}) and
\begin{eqnarray*}
(n-2)\, ( R \cdot C - C \cdot R ) &=& Q(S,C) 
- \frac{\widetilde{\kappa }}{ n (n+1) }\, Q(g,C) ,
%\label{quasi321}
\\
\alpha &=& \frac{\kappa}{n-1} - \frac{ \widetilde{\kappa} }{n(n+1)} ,
%\label{qqee02}
\\
R \cdot C - C \cdot R &=&  \frac{\kappa}{n-1}\, Q(g,C) - Q(S,C) + B ,
%\label{qqee03}
\end{eqnarray*}
where 
$\psi$ is some function on ${\mathcal U}_{H}$,
and the function $\alpha$ 
and the $(0,6)$-tensor $B$ are defined by
(\ref{quasi02}) and  
\begin{eqnarray*}
B &=& Q\left(
S -  \frac{1}{n-1} \left( \frac{ (n-2) \kappa }{n-1}
 +  \frac{ \widetilde{\kappa} }{n(n+1)} \right) g , C \right) ,
%\label{qqee10}
\end{eqnarray*}
respectively. Moreover, 
(\ref{Roterformula}) is satisfied on ${\mathcal U}_{H}$ 
if and only if (\ref{qqee08}) holds  on this set. 
\end{thm}

We consider non-quasi-Einstein hypersurfaces satisfying (\ref{cond01}). 
Precisely, we consider (\ref{cond01}) at all points of ${\mathcal U}_{H}$ 
at which the following condition is satisfied
\begin{eqnarray} 
\mathrm{rank} ( S - \alpha \, g ) > 1 ,\ \ 
\mbox{for any}\ \ \alpha \in  {\mathbb{R}} .
\label{dhs_quasi_Einstein}
\end{eqnarray}
We have
\begin{thm} {\cite[Theorem 5.2]{2020_DGZ}}
Let $M$ be a hypersurface in a semi-Riemannian space of constant curvature 
$N_{s}^{n+1}(c)$, $n \geq 4$, satisfying (\ref{cond01})
on ${\mathcal U}_{H} \subset M$.
If at a point $x \in {\mathcal U}_{H}$ 
(\ref{dhs_quasi_Einstein}) 
is satisfied then at this point we have (\ref{Roterformula}) and 
\begin{eqnarray*}
(n-1)\, Q(S,R) 
&=& Q\left( g,
\left( \varepsilon \psi  + \kappa 
- \frac{ (n-1) \widetilde{\kappa}}{n (n+1)} \right)  \right) R 
+
\left( \varepsilon \psi 
- \frac{ 2 (n-1) \widetilde{\kappa} }{n (n+1)} \right) g \wedge S 
+ g \wedge S^{2} - \frac{1}{2}\, S \wedge S  \,  .
%\label{cond02uuu}
\end{eqnarray*}
\end{thm} 

\begin{thm} {\cite[Theorem 5.4]{2020_DGZ}}
If at every point of
a non-quasi-Einstein hypersurface $M$ 
in a semi-Rieman\-nian space of constant curvature 
$N_{s}^{n+1}(c)$, $n \geq 4$, 
the difference tensor $R \cdot C - C \cdot R$  
is a linear combination of the Tachibana tensors 
$Q(g,C)$ and $Q(S,C)$, then
(\ref{Roterformula}) holds on $M$.
\end{thm}

\section{Chen ideal submanifolds}

Let $ M $ be a
submanifold of  dimension $ n $ in the Euclidean space 
$ \mathbb{E}^{n+m} $,
$n \geq 2$, $m \geq 1$.
Let $ g $ be the Riemannian metric induced
on $ M $ from the standard metric on $ \mathbb{E}^{n+m} $,
$ \nabla $ the corresponding Levi-Civita connection on $ M $,
and
$ R $, $ S $, $ \tau $ and $C$, 
the Riemann-Christoffel curvature tensor, 
the Ricci tensor, the scalar curvature 
and the Weyl conformal curvature tensor of $g$, respectively.
For the scalar curvature $ \tau $ of $ (M , g) $ we use the
calibration
$$
 \tau (p)  \ =\  \sum_{i < j} K\left(p , e_i(p) \wedge
e_j(p)\right)
$$
  where $  K\left(p , \pi\right) $ denotes  the
{\sl Riemannian sectional curvature} of  $ (M , g) $ at the point  $
p $ for a plane section $ \pi $ in
   the tangent  space $ T_p M $.
    For each point  $ p $ in $ M $, considering the
    number
 $$
  \left(\mbox{\rm inf}\, K\right) (p) : =
  \mbox{\rm inf} \left\{K(p , \pi)  \vert   \pi \ \text{is a  plane
    section in} \
  T_p M\right\} ,
 $$
 B.-Y. Chen %\cite{BYCH2011}
(see \cite{{BYCH1993}, {Chen-2011}})
 introduced the
  $ \delta(2) $-curvature by
  $$
 \left(\delta(2)\right)(p)
   = \delta(p)
  := \tau(p) - \left(\mbox{\rm inf}\, K\right) (p).
  $$
 This  $ \delta(2) $-, for  short,
   {\sl $ \delta  $-curvature of Chen}
   thus is a well
   defined real  function on $ M $
   which clearly is a
 {\sl Riemannian invariant} of    $ (M , g) $.
  % \noindent
  
From \cite{Chen-2011}
 (see also \cite{{BYCH1993}, {BYCH2008}, {ChenMarMol}}), we have
 the following basic result which, in
 particular, answered a question raised by
S.S. Chern
 \cite{CHERN1968} long before,
 concerning {\sl intrinsic obstructions}
 on Riemannian manifolds
 {\sl in view of minimal immersibility}
 in Euclidean   spaces.

  \begin{thm} 
{\rm \cite{BYCH1993}} (see also {\cite[Theorem 1] {DP-TVZ}}) 
\label{BYCH1} \
 %  {\bf (B. Y. Chen, 1993)} \label{BYCH1}\
 For any submanifold
   $ M $  of  dimension $ n $
   in the Euclidean space  $ \mathbb{E}^{n+m} $,
$ n \geq 2 $, $ m \geq 1 $,
\begin{eqnarray}
 \delta  \leq
     \dfrac{n^2(n-2)}{2(n-1)}   H^{2} , 
\label{2022.08.08.a}         
\end{eqnarray}
and in (\ref{2022.08.08.a}) equality
holds at a point $ p \in M $ if and
  only if,
  with respect to some suitable adapted orthonormal
 frame
  $  \left\{e_i , \xi_\alpha\right\}    $
 around $ p $ on $ M $ in $ \mathbb{E}^{n+m} $, the
  shape operators are given by
 %   $$
 %  A_\alpha := A_{\xi_\alpha} , \quad \alpha = 1 , \cdots , m ,
 % $$
  $$
  A_1 =  \begin{pmatrix}
  a & 0 & 0 & \cdots & 0 \\
  0 & b & 0 & \cdots & 0 \\
  0 & 0 & z & \cdots & 0 \\
  \vdots & \vdots & \vdots & \ddots & \vdots \\
  0 & 0 & 0 & \cdots & z \\
  \end{pmatrix}
  , \quad
 A_\beta =  \begin{pmatrix}
  c_\beta & d_\beta & 0 & \cdots & 0 \\
  d_\beta & -c_\beta & 0 & \cdots & 0 \\
  0 & 0 & 0 & \cdots & 0 \\
  \vdots & \vdots & \vdots & \ddots & \vdots \\
  0 & 0 & 0 & \cdots & 0 \\
  \end{pmatrix}
, \
  \beta >  1 ,
  $$
 where
 $ z =  a + b $ and
 $ \inf K     =
  ab - \sum_{\beta >1}  \left(c^2_\beta +  d^2_\beta\right)
  :  M \longrightarrow \mathbb{R} $.
      \end{thm}

\noindent Evidently, if $m = 1$ then $\inf K = ab$.

    With respect to the above theorem,
 one has the following    definition 
 \cite{{BYCH1993}, {BYCH2008}, {DGPV},  {2008_DP-TVZ}, {LVerstr1998}}, 
also see {\cite[Definition 1] {DP-TVZ}}.
Let $ M $ be a  submanifold of  dimension $ n $
in the Euclidean space  $ \mathbb{E}^{n+m} $,  $ n \geq 2 $,
$ m \geq 1 $.
It is called a {\sl{Chen ideal submanifold}} if, at each of its points,
the Chen's basic inequality (\ref{2022.08.08.a}) in
the Theorem  \ref{BYCH1} is actually an equality.

     Let  $ M $ be a Chen ideal submanifold
	of  dimension $ n $
    in the Euclidean space  $ \mathbb{E}^{n+m} $, $ n \geq 4 $,
    $ m \geq 1 $.
	We    use the notations
 as  in Theorem \ref{BYCH1}.
The  Riemann-Christoffel curvature tensor $ R $ satisfies:
 \begin{equation}   \label{eq001*}
 \left\{
  \begin{aligned}
  R\left(e_1 , e_2 , e_2 , e_1\right) &= \inf K
   = ab  - \sum_{\beta = 2}^m \left(c^2_\beta + d^2_\beta\right)   ; \\
   R\left(e_1 ,  e_i , e_i , e_1\right)      &
   =    az  \quad \text{\rm for} \quad i \geq 3   ;  \\
   R\left(e_2 , e_i , e_i , e_2\right)        &
   =   bz  \quad \text{\rm for} \quad i \geq 3   ;  \\
   R\left(e_i ,  e_j , e_j , e_i\right)      &
   =     z^2    \quad \text{\rm for} \quad 3 \leq i < j \leq n  .  \\
     \end{aligned}
     \right.
 \end{equation}
The other values of $ R\left(e_u , e_v , e_w , e_t\right) $ are null.
The Ricci tensor $ S $ satisfies:
   \begin{equation}\label{eq001**}
 \left\{
  \begin{aligned}
 S\left(e_1 , e_1\right)    &=  \inf K  + (n-2)az    ; \\
 S\left(e_2 , e_2\right)    &=  \inf K  + (n-2)bz        ;  \\
 S\left(e_i , e_i\right)    
 &=    (n-2)z^2   \quad \text{\rm for} \quad 3 \leq i \leq n    ; \\
 S\left(e_u , e_v\right)    
 &=      0    \quad \text{\rm for} \quad 1 \leq u < v \leq n  .  \\
     \end{aligned}
     \right.
 \end{equation}

    \noindent    The scalar curvature $ \tau $ is given by
 \begin{equation}\label{eq5}
\tau = \sum_{i=1}^n \mbox{\rm S}(e_i , e_i) = 2 \inf K +
 (n-1)(n-2)z^2 .
\end{equation}

 % \item
  \noindent    The  Weyl conformal curvature tensor $ C $ of $M$ is
  determined by the following relations:
 \begin{equation}\label{eq701}
 \left\{
  \begin{aligned}
 C\left(e_1 , e_2 , e_2 , e_1\right)
  &=  \dfrac{(n - 3) \inf K}{n - 1}     ; \\
 C\left(e_1 , e_i , e_i , e_1\right)
   &=    -\dfrac{(n - 3) \inf K }{(n - 1)(n -2)} \ \ \mbox{for}\ \  
   \geq 3   \\
 C\left(e_2 , e_i , e_i , e_2\right)
  &=   -\dfrac{(n - 3) \inf K }{(n - 1)(n -2)} \ \ \mbox{for}\ \ \geq 3 \\
  C\left(e_i , e_j , e_j , e_i\right)
  &=  \dfrac{2 \inf K }{(n - 1)(n -2)}
   \quad \text{\rm for} \quad 3 \leq i < j \leq n . \\
     \end{aligned}
     \right.
 \end{equation}

 % \end{Myitemize}

 \noindent From (\ref{eq701})
  it follows (cf. \cite{DGPV}, {\cite[Theorem F]{2008_DP-TVZ}}) 
that every Chen ideal submanifold $ M $
of  dimension $ n $
in the Euclidean space ${\mathbb E}^{n+m}$, $n \geq 4$, $m \geq 1$,
has a pseudosymmetric Weyl conformal curvature $ C $,
i.e., it satisfies the identity
\begin{eqnarray}
C \cdot C &=& L_C\, Q(g,C) ,\ \ \ 
L_C \ =\ - \dfrac{(n-3)\inf K}{(n-1)(n-2)} . 
\label{tachibconf}
\end{eqnarray}

 \noindent As it was mentioned in Section 8, 
 at every point of a hypersurface $M$
 in a space form $ \widetilde{N}^{n+1}(c)$, $n \geq 4$, the tensors 
 $R \cdot R - Q(S,R)$
 and $Q(g,C)$ are linearly dependent. 
 Precisely, (\ref{900ab}) holds on $M$.
 Thus, in particular, 
$R \cdot R - Q(S,R) = 0$ on every Chen ideal hypersurface $M$ in
$ \mathbb{E}^{n+1} $, $ n \geq 4 $.

 \noindent Now  let us compute the difference
$ R \cdot R -   Q(S,R) $
on the Chen ideal submanifold $ M $
of codimension $ m $ in $ \mathbb{E}^{n+m} $,
	$ n \geq 4 $, $ m \geq 1 $.
With respect to
 the notations
   in Theorem \ref{BYCH1} and from
   the equalities
 (\ref{eq001*}),
(\ref{eq001**})
  and  (\ref{eq701}),
    we can prove the following

 \begin{thm}\label{RRQSSR} {\cite[Theorem 2] {DP-TVZ}}
    The identity
\begin{eqnarray*}
    R \cdot R -   Q(S,R) = \dfrac{ab - \inf K}{\inf K}
    (n-2)z^2\, Q(g,C)
\end{eqnarray*}
  holds
 on the subset
  $   {\mathcal U}_{C} $ 
   of every Chen ideal submanifold $ M $
	of  dimension $ n $
	in the Euclidean space $ \mathbb{E}^{n+m} $,
	$ n \geq 4 $, $ m \geq 1 $.
  In addition,
    at each point $ p \in M $ where $ C $ vanishes
		($ \inf K = 0 $),
    the following equalities hold:
\begin{eqnarray*}
 \left\{
\begin{aligned}
 \left(R \cdot R -   Q(S,R)\right)\left(e_1 , e_2 , Z , W; e_1 , e_i\right)
        &=  (n-3)abz^2\left<\left(e_2\wedge_g e_i\right)(Z),   W\right>, \\
     \left(R \cdot R 
     -   Q(S,R)\right)\left(e_1 , e_j , Z , W; e_1 , e_i\right)
        &=  -abz^2\left<\left(e_j \wedge_g e_i\right)(Z) ,W\right>, \\
 \left(R \cdot R -   Q(S,R)\right)\left(e_1 , e_2 , Z , W; e_2 , e_i\right)
      &=  -(n-3)abz^2\left<\left(e_1\wedge_g e_i\right)(Z),   W\right>, \\
 \left(R \cdot R -   Q(S,R)\right)\left(e_2 , e_j , Z , W; e_2 , e_i\right)
        &=  -abz^2\left<\left(e_j \wedge_g e_i\right)(Z) ,W\right> ,\\
\end{aligned}
\right.
\end{eqnarray*}
 \noindent and the other values of
$\left(R \cdot R - Q(S,R)\right)\left(e_u , e_v, Z , W; e_w , e_t\right) $
  being  null.
    \end{thm}

%\noindent  

As it was proved in \cite{2016_DGJZ}, every warped product manifold
   $\overline{M} \times _{F} \widetilde{N}$
    of a $2$-dimensional base manifold
 $(\overline{M},\overline{g})$ and an $(n-2)$-dimensional
fibre, which is a space of constant curvature
   $(\widetilde{N},\widetilde{g})$, $n\geq 4$,
    with the warping  function $F$, satisfies
\begin{equation}  \label{WeylWeyl07b}
  C \cdot C =
- \frac{(n-3) \rho }{(n-2)(n-1)} \, Q(g,C)  ,\ \ \ \
\rho  =  \frac{ \overline{\tau}}{2}
 + \frac{\widetilde{\tau}}{(n-3) (n-2) F }
  + \frac{\Delta F}{2 F} - \frac{\Delta _1 F}{2 F^{2}},
\end{equation}

  \noindent where $\Delta F = g^{ab} \nabla _{b} F_{a}$, 
  $\Delta _1 F = g^{ab} F_{a} F_{b}$, 
  and $\overline{\tau }$, $\widetilde{\tau }$
are the scalar curvatures of the base and the fibre, respectively,
see to \cite{2016_DGJZ} for details.

As we noted in Section 1,
some warped product manifolds
$\overline{M} \times _{F} \widetilde{N}$ 
with a $2$-dimensional Riemannian manifold 
$(\overline{M},\overline{g})$ and 
an
$(n-2)$-dimensional unit sphere ${\mathbb S}^{n-2}$, 
are related to Chen ideal submanifolds. 
Namely, according to \cite{MDLAF}, every
non-trivial and non-minimal Chen ideal
submanifold  $M$ of  dimension $ n $ in the Euclidean
space ${\mathbb E}^{n+m}$,
$ n \geq 4 $, $ m \geq 1 $ is  isometric to an open subset of a warped
product
   $\overline{M} \times _{F} {\mathbb S}^{n-2}$ of a
 $2$-dimensional base manifold
  $(\overline{M},\overline{g})$ and an
$(n-2)$-dimensional unit sphere ${\mathbb S}^{n-2}$,
 where  the warping function $F$ is a solution of some second order
quasilinear elliptic partial differential equation in the plane.
Thus we see that (\ref{WeylWeyl07b}) holds on $M$. Furthermore, from
(\ref{tachibconf}) and (\ref{WeylWeyl07b}) it follows that $\inf K$
is expressed on $M$ by
   \begin{eqnarray*}
\inf K &=& \frac{ \overline{\tau } }{2} + \frac{ 1 }{ F } +
\frac{\Delta F}{2 F} - \frac{\Delta _1 F}{2 F^{2}} .
\label{WeylWeyl10}
   \end{eqnarray*}
Since the scalar curvature $\tau $ of $M$
 is given by (\ref{eq5})
 and satisfies (see, e.g., \cite{2016_DGJZ})
     \begin{eqnarray*}
\tau &=& \overline{\tau} + \frac{ (n-3)(n-2) }{ F } - \frac{(n-2)
\Delta F}{ F} - \frac{(n-2) (n-5) \Delta _1 F}{4 F^{2}} ,
\label{WeylWeyl12}
    \end{eqnarray*}
  we get
      \begin{eqnarray*}
(n-2) z^{2} &=& \frac{ n-4 }{ F } - \frac{\Delta F}{ F} - \frac{
(n-6) \Delta _1 F}{4 F^{2}} . \label{WeylWeyl11}
  \end{eqnarray*}

   \noindent  In \cite{2016_DGJZ},
 among other things, it was shown that the Hessian of the function
 $ f = \sqrt{F} $ is proportional to the metric $ \overline{g} $
 on $ \overline{M} $.

We mention that Chen ideal submanifolds which 
are semisymmetric were classified
in \cite{DPV}. We have
 
 \begin{thm} 
 {\cite{DPV}} (see also {\cite[Theorem 3] {DP-TVZ}}) \label{seisym}
  A Chen ideal submanifold $ M $ of  dimension $ n $
  in the Euclidean space $ {\mathbb E}^{n+m} $, $ n \geq 3 $, $ m \geq 1 $,
   is semisymmetric if and only if
   $ M $  is minimal (in which case $ M $ is $ (n-2) $-ruled) or
    $ M $  is
   a round hypercone in some totally geodesic
   subspace $ {\mathbb E}^{n+1} $ of  $ {\mathbb E}^{n+m} $.
     \end{thm}

Chen ideal pseudosymmetric submanifolds were classified
in \cite{{DGPV}, {2008_DP-TVZ}}. We have

\begin{thm} \cite{{DGPV}, {2008_DP-TVZ}} 
 (see also {\cite[Theorem 4] {DP-TVZ}})\label{pseudosemisym}  
   A Chen ideal submanifold $ M $ of codimension $ m $
  in $ {\mathbb E}^{n+m} $
   ($ n \geq 3 $, $m \geq 1 $)
   is pseudosymmetric if and only if
\newline   
(i) either  $ M $  is  semisymmetric (see Theorem \ref{seisym}),
\newline
(ii) or at every point $ p $ of $ M $ where
        $ R \cdot R \ne 0 $,
        the $ 2D $ normal section
        $ \Sigma^2_{\tilde \pi} \subset {\mathbb E}^{2+m} $
        of $ M^n $ at $ p $ in the direction of the tangent
        plane $ {\tilde \pi} \subset T_p M^n $
        for which the sectional curvature function
        $ K(p , \pi)  $ at $ p $ attains  its minimal
        value $ \left(\inf K\right)(p) $ is pseudo-umbilical at $ p $, or
    equivalently, if $ p $ is a spherical point
                of the projection
 $ \overline{\Sigma}^2_{\tilde \pi} \subset {\mathbb E}^{3} $
        of this $ 2D $ normal section $ \overline{\Sigma}^2_{\tilde \pi} $
        on the space $ {\mathbb E}^{3} $ spanned by
        $ {\tilde \pi}  $ and the mean curvature vector
        $ \stackrel{\rightarrow}{H}(p) $  of
        $ M^n $   in  $ {\mathbb E}^{n+m} $ at $ p $
        (and in this case
        $ L_R = \dfrac{n^2}{2(n-1)^2} H^2$, 
        where $ H $  is the mean curvature
        of $ M^n $   in  $ {\mathbb E}^{n+m} $).
    \end{thm}

The above presented part of this section based on 
{\cite[Section 3]{DP-TVZ}})
and the remains part of this section on 
{\cite[Section 5]{DP-TVZ}}. 
We present now results on Chen ideal submanifolds
 % $ M $ of codimension $ m $
 in Euclidean spaces
 whose difference tensor  $ R \cdot C - C \cdot R $ can be expressed
 in terms of some of the Tachibana  tensors
 $ Q(g , R) $, $  Q(S , R) $,
 $ Q(g , C) $, $  Q(S , C) $,
  $ Q(g , g\wedge S) $, $  Q(S , g\wedge S) $.

 \begin{thm} 
 \label{16decembre1C} {\cite[Theorem 6] {DP-TVZ}}\ % 
 \label{r.c-c.r=L1Q(gC) + L2Q(SC)}    % thrm 6  
 Let $ M $ be a non-conformally flat
   Chen   ideal  submanifold  of  dimension $ n $
   in the Euclidean space  $ \mathbb{E}^{n+m} $, $n \geq 4$, $m \geq 1$.
  Then 
there exist two real valued functions $ L_1 $ and $ L_2 $ on $ M $
   such that (\ref{cond01}), i.e.,
\begin{eqnarray*}
  R \cdot C - C \cdot R = L_1 Q(S  , C) + L_2 Q(g  , C)   ,
%\label{2023.05.20.c}
\end{eqnarray*}
if and only if   there exists
    an orthonormal  tangent  framefield
  $  \left\{e_1 , \cdots , e_n\right\}    $ and
     an orthonormal  normal  framefield
$  \left\{\xi_1 , \cdots , \xi_m\right\}    $
  such that the shape operators
  $A_\alpha :=  A_{\xi_\alpha} , \quad 1 \leq \alpha \leq m $,
    are given by:
    $$
 %  \begin{equation}\label{eqbli0QSgWS1}
  A_1 =  \begin{pmatrix}
  a      &     0        &        0          & \cdots & 0 \\
  0      &  \epsilon a  &        0          & \cdots & 0 \\
  0      &     0        &  (1 + \epsilon)a  & \cdots & 0 \\
  \vdots &   \vdots     &        \vdots     & \ddots & \vdots \\
  0      &     0        &        0          & \cdots & (1 + \epsilon)a \\
  \end{pmatrix}
  ,  \quad
 A_\beta =  \begin{pmatrix}
  c_\beta & d_\beta & 0 & \cdots & 0 \\
  d_\beta & -c_\beta & 0 & \cdots & 0 \\
  0 & 0 & 0 & \cdots & 0 \\
  \vdots & \vdots & \vdots & \ddots & \vdots \\
  0 & 0 & 0 & \cdots & 0 \\
  \end{pmatrix}
, \
  \beta \geq 2 ,
   $$
 where
  $ \epsilon = \pm 1 $,
  $ a $,
 $ c_\beta $ , $ d_\beta $ (for $ 2 \leq \beta \leq m $)
  are real functions on $ M $
 such that
  $$
    \sum_{\beta = 2}^m \left(c^2_\beta + d^2_\beta\right)
   = \epsilon a^2 -  \inf K ,\ \ \
   L_1  =  - 1 , \quad
    L_2 =
   - \dfrac{2\, \inf K + 2(n-1)(n-2)(1 + \epsilon)a^2}{n-1} .
   $$

  \noindent     In this  case, $ M $ is a Roter space.
  In addition  one has one of the following two cases:
\newline
(i)  either $ \epsilon = -1 $
  and
   $ M  $ is a semisymmetric
   and   minimal  submanifold
    (see   Theorem 11.3)
 such that
      $$
   R \cdot C - C \cdot R    =
    \dfrac{2\, \inf K}{n-1} Q(g , C)
    - Q(S , C) ,
   $$
(ii) or  $ \epsilon = +1 $ and
  $ M $ is  a properly  pseudosymmetric
 and    non-minimal
   submanifold
     (see   Theorem 11.4)
       such that
      $$
      R \cdot C - C \cdot R    =
     - \dfrac{2\, \inf K + 4(n-1)(n-2)a^2}{n-1} Q(g , C)
    - Q(S , C) .
    $$
\end{thm}

   \begin{cor}  \label{r.c-c.r=LQ(gC)Resum} 
   {\cite[Corollary 3] {DP-TVZ}}  % corll 3   
   Let $ M $ be a
   Chen   ideal   submanifold   of  dimension $ n $
   in the Euclidean space  $ \mathbb{E}^{n+m} $, $n \geq 4$, $m \geq 1$.
   Then the difference tensor
    $ R \cdot C - C \cdot R $ and
    the Tachibana tensor  $ Q(g , C) $ are linearly dependent
if and only if  $ M $ is conformally flat.

    \end{cor}

   \begin{cor}  \label{r.c-c.r=LQ(SC)} 
   {\cite[Corollary 4] {DP-TVZ}}  % thrm 7  
   Let $ M $ be a non-conformally flat
   Chen   ideal   submanifold   of  dimension $ n $
   in the Euclidean space  $ \mathbb{E}^{n+m} $,
$n \geq 4$, $m \geq 1$.
    Then there exists a real valued function $ L $ on $ M $
   such that
 $$
  R \cdot C - C \cdot R = L Q(S  , C)   ,
  $$
 if and only if $ M $ is not minimal, and
   there exists
    an orthonormal  tangent  framefield
  $  \left\{e_1 , \cdots , e_n\right\}    $ and
     an orthonormal  normal  framefield
$  \left\{\xi_1 , \cdots , \xi_m\right\}    $
  such that the shape operators
  $A_\alpha :=  A_{\xi_\alpha} , \quad 1 \leq \alpha \leq m$,
  are given by
   $$
  %  \begin{equation}\label{eqbli0QSC}
  A_1 =  \begin{pmatrix}
  a      &    0    &   0    & \cdots & 0 \\
  0      &    a    &   0    & \cdots & 0 \\
  0      &    0    &  2a    & \cdots & 0 \\
  \vdots & \vdots  & \vdots & \ddots & \vdots \\
  0      &    0    &   0    & \cdots & 2a \\
  \end{pmatrix}
  ,  \quad
 A_\beta =  \begin{pmatrix}
  c_\beta & d_\beta & 0 & \cdots & 0 \\
  d_\beta & -c_\beta & 0 & \cdots & 0 \\
  0 & 0 & 0 & \cdots & 0 \\
  \vdots & \vdots & \vdots & \ddots & \vdots \\
  0 & 0 & 0 & \cdots & 0 \\
  \end{pmatrix}
, \
 \beta \geq 2 ,
   $$
 where $ a $,
 $ c_\beta $ , $ d_\beta $ (for $ 2 \leq \beta \leq m $)
  are real functions on $ M $
 such that
$$ \sum_{\beta = 2}^m \left(c^2_\beta + d^2_\beta\right)   
= \left(2n^2 - 6n + 5\right) a^2 > 0 $$
and $L = - 1$. In this case,  $ M^n $ is  properly  pseudosymmetric 
(see Theorem 11.4).
     
      \end{cor}

\begin{thm} {\cite[Theorem 5] {DP-TVZ}} 
   \label{16decembre1R} % \label{r.c-c.r=L1Q(gR) + L2Q(SR)}    % thrm 5
   Let $ M $ be a non-conformally flat
   Chen   ideal   submanifold of  codimension $ m $
   in the Euclidean space  $ \mathbb{E}^{n+m} $, $n \geq 4$, $m \geq 1$.
    Then there exist
		two real valued functions $ L_3 $ and $ L_4 $ on $ M $
   such that (\ref{2023.05.20.b}) is satisfied on $M$
if and only if there exists an orthonormal  tangent  framefield
  $  \left\{e_1 , \cdots , e_n\right\}    $ and
     an orthonormal  normal  framefield
$  \left\{\xi_1 , \cdots , \xi_m\right\}    $
  such that the shape operators
  $ A_\alpha :=  A_{\xi_\alpha}$, 
  $\quad 1 \leq \alpha \leq m $ are given by:
    $$
 %\begin{equation}\label{eqbli0QSgWS1}
  A_1 =  \begin{pmatrix}
     a   &   0        &       0          & \cdots &      0       \\
     0   & \epsilon a &       0          & \cdots &      0       \\
     0   &   0        &  (1 + \epsilon)a & \cdots &      0       \\
  \vdots & \vdots     &   \vdots         & \ddots &    \vdots      \\
    0    &    0       &       0          & \cdots &  (1 + \epsilon)a \\
  \end{pmatrix}
  ,  \quad
 A_\beta =  \begin{pmatrix}
  c_\beta & d_\beta & 0 & \cdots & 0 \\
  d_\beta & -c_\beta & 0 & \cdots & 0 \\
  0 & 0 & 0 & \cdots & 0 \\
  \vdots & \vdots & \vdots & \ddots & \vdots \\
  0 & 0 & 0 & \cdots & 0 \\
  \end{pmatrix}
, \   \beta \geq 2 ,
    $$
 where
 $ \epsilon = \pm 1 $,
  $ a $,
 $ c_\beta $ , $ d_\beta $ (for $ 2 \leq \beta \leq m $)
  are real functions on $ M $
 such that
  $$
    \sum_{\beta = 2}^m \left(c^2_\beta + d^2_\beta\right)
   = \epsilon a^2 -  \inf K, \ \ \
   2\, \inf K-  (1+\epsilon)a^2 \ne 0 , \ \ \
 \inf K \ne (n-2)(1 + \epsilon)a^2
   $$
    and
      $$
  \begin{aligned}
   &L_3  =  \dfrac{(n-3)\, \inf K -  (n - 2)(1+\epsilon)a^2}{(n-1)(n-2)}
           \dfrac{2\, \inf K}{2\, \inf K -  (1+\epsilon)a^2} ,   \\
    &(1+\epsilon)a^2 \left(L_4  -
\dfrac{1}{(n-2)}\dfrac{\inf K}{2\, \inf K - (1+\epsilon)a^2}\right) = 0 .
 \end{aligned}
   $$
\noindent     
In this  case, $ M $ is a Roter space.
  In addition  one has one of the following two cases:
\newline  
(i)  either $ \epsilon = -1 $
  and
   $ M  $ is a semisymmetric
 and minimal
   submanifold
    (see Theorem 11.3)
 such that
      $$
   R \cdot C - C \cdot R  = \dfrac{(n-3)\, \inf K}{(n-1)(n-2)} Q(g , R) ,
   $$
(ii) or  $ \epsilon = +1 $ and
  $ M $ is  a
    properly  pseudosymmetric
 and  non minimal
    submanifold
     (see   Theorem 11.4)
       such that
      $$
 %   \begin{aligned}
   R \cdot C - C \cdot R    =
    \dfrac{(n-3)\, 
    \inf K - 2(n - 2)a^2}{(n-1)(n-2)}\dfrac{\inf K}{\inf K-  a^2} Q(g,R)
    % \\
  %   &\quad
   + \dfrac{1}{2(n-2)}\dfrac{\inf K}{\inf K -  a^2} Q(S , R) .
   %    \\      \end{aligned}
   $$
\end{thm}

\begin{cor} {\cite[Corollary 1] {DP-TVZ}} \label{r.c-c.r=LQ(gR)}  
   % coroll 1
 Let $ M $ be a
 non-conformally flat Chen ideal submanifold
 of  dimension $ n $  in the Euclidean space  $ \mathbb{E}^{n+m} $, 
 $n \geq 4$, $m \geq 1$.
   Then there exists a real valued function $ L $ on $ M $
   such that
 $$
  R \cdot C - C \cdot R = L Q(g, R)
  $$
%    for
 %  some function $ L : M \to {\mathbb R} $,
 if and only if
   $ M $ is minimal and
 there exists
    an orthonormal  tangent  framefield
  $  \left\{e_1 , \cdots , e_n\right\}    $ and
     an orthonormal  normal  framefield
$  \left\{\xi_1 , \cdots , \xi_m\right\}    $
  such that the shape operators
  $A_\alpha :=  A_{\xi_\alpha} , \quad 1 \leq \alpha \leq m$,
    are given by:
    $$
 % \begin{equation}\label{eqbli0QgR}
  A_1 =  \begin{pmatrix}
  a & 0 & 0 & \cdots & 0 \\
  0 & -a & 0 & \cdots & 0 \\
  0 & 0 & 0 & \cdots & 0 \\
  \vdots & \vdots & \vdots & \ddots & \vdots \\
  0 & 0 & 0 & \cdots & 0 \\
  \end{pmatrix}
  , \quad
 A_\beta =  \begin{pmatrix}
  c_\beta & d_\beta & 0 & \cdots & 0 \\
  d_\beta & -c_\beta & 0 & \cdots & 0 \\
  0 & 0 & 0 & \cdots & 0 \\
  \vdots & \vdots & \vdots & \ddots & \vdots \\
  0 & 0 & 0 & \cdots & 0 \\
  \end{pmatrix}
, \
   \beta \geq  2 ,
    $$
 where $ a $,
 $ c_\beta $ , $ d_\beta $ (for $ 2 \leq \beta \leq m $)
  are real functions on $ M $
    such that
  $$
   \sum_{\beta = 2}^m \left(c^2_\beta + d^2_\beta\right)
   = -a^2 -  \inf K,\ \ \  L = \dfrac{(n-3)\,\inf K}{(n-1)(n-2)} .$$
\noindent  
In this case,
   $ M $ is semisymmetric (see   Theorem 11.3).
 
   \end{cor}

\begin{cor}  \label{r.c-c.r=LQ(SR)Resum} 
{\cite[Corollary 2] {DP-TVZ}}  % corll 2  
  Let $ M $ be a
   Chen   ideal   submanifold  of  dimension $ n $
   in the Euclidean space  $ \mathbb{E}^{n+m} $, $n \geq 4$, $m \geq 1$.
   Then
      the difference tensor   $ R \cdot C - C \cdot R $ and
        the Tachibana tensor $ Q(S , R) $
        are linearly dependent
    if and only if
        $ M $ is conformally flat
				($ \inf K = 0 $).

\end{cor}

   \begin{thm} \label{29juillet1} 
{\cite[Theorem 7] {DP-TVZ}}\ % \label{r.c-c.r=L1Q(ggWS) + L2Q(SgWS)} 
% thrm 7
    Let $ M $ be a non-conformally flat
   Chen   ideal   submanifold  of  dimension $ n $
   in the Euclidean space  $ \mathbb{E}^{n+m} $,
   $n \geq 4$, $m \geq 1$.
    Then there exist two real valued functions 
    $ L_5 $ and $ L_6 $ on $ M $
   such that (\ref{2023.05.20.d}) is satisfied on $M$
if and only if
    there exists
    an orthonormal  tangent  framefield
  $  \left\{e_1 , \cdots , e_n\right\}    $ and
     an orthonormal  normal  framefield
$  \left\{\xi_1 , \cdots , \xi_m\right\}    $
  such that the shape operators
  $  A_\alpha :=  A_{\xi_\alpha} , \quad 1 \leq \alpha \leq m $,
   are given by:
  
  $$
 %  \begin{equation}\label{eqbli0QSgWS1}
  A_1 =  \begin{pmatrix}
     a   &   0        &       0          & \cdots &      0       \\
     0   & \epsilon a &       0          & \cdots &      0       \\
     0   &   0        &  (1 + \epsilon)a & \cdots &      0       \\
  \vdots & \vdots     &   \vdots         & \ddots &    \vdots      \\
    0    &    0       &       0          & \cdots &  (1 + \epsilon)a \\
  \end{pmatrix}
  ,  \quad
 A_\beta =  \begin{pmatrix}
  c_\beta & d_\beta & 0 & \cdots & 0 \\
  d_\beta & -c_\beta & 0 & \cdots & 0 \\
  0 & 0 & 0 & \cdots & 0 \\
  \vdots & \vdots & \vdots & \ddots & \vdots \\
  0 & 0 & 0 & \cdots & 0 \\
  \end{pmatrix}
, \
  \beta \geq 2 ,
    $$
 where
 $ \epsilon = \pm 1 $,
  $ a $,
 $ c_\beta $ , $ d_\beta $ (for $ 2 \leq \beta \leq m $)
  are real functions on $ M $
 such that
  $$
   \sum_{\beta = 2}^m \left(c^2_\beta + d^2_\beta\right)
   = \epsilon a^2 -  \inf K  , \ \ \
 \inf K \ne (n-2)(1 + \epsilon)a^2   $$
   and moreover
      $$
    \begin{aligned}
   L_5 &=
   - \dfrac{2(n-2)(1 + \epsilon)a^2}{n-1}\dfrac{\inf K \, 
   \left(\inf K + (n-2)(1  + \epsilon)a^2 \right)}{\left(\inf K 
   - (n-2)(1 + \epsilon)a^2 \right)^2} ,
     \\
   L_6 &=
   - \dfrac{1}{(n-1)(n-2)}\dfrac{\inf K\left[(n-3)\inf K
  + (n-1)(n-2)(1 + \epsilon)a^2\right]}{\left(\inf K 
  - (n-2)(1 + \epsilon)a^2 \right)^2} .
     \\
        \end{aligned}
   $$
In this  case, $ M $ is a Roter space.
  In addition  one has one of the following two cases:
\newline
(i)  either $ \epsilon = -1 $
  and
   $ M  $ is a semisymmetric
 and     minimal submanifold
    (see   Theorem 11.3)
 such that
 $$
 R \cdot C - C \cdot R =
 - \dfrac{n-3}{(n-1)(n-2)} Q(S , g \wedge S) ,
$$
(ii) or  $ \epsilon = +1 $ and
  $ M $ is  a  properly  pseudosymmetric
 and
    non-minimal
  submanifold
     (see   Theorem 11.4)
       such that
      $$
    \begin{aligned}
   R \cdot C - C \cdot R   &=
  - \dfrac{4(n-2)a^2}{n-1}\dfrac{\inf K\left(\inf K
	             + 2(n-2)a^2\right)}{\left(\inf K - 2(n-2)a^2\right)^2}
    Q(g , g \wedge S)  \\
     &\quad
  - \dfrac{1}{(n-1)(n-2)}\dfrac{\inf K \left( (n-3)\, \inf K
    	+ 2(n-1)(n-2)a^2\right)}{\left(  \inf K - 2(n-2)a^2\right) ^2}
  Q(S , g \wedge S) .
       \\
         \end{aligned}
   $$
   \end{thm}

   \begin{cor}  \label{r.c-c.r=LQ(ggWS)} 
   {\cite[Corollary 5] {DP-TVZ}} % coroll 5   
   Let $ M $ be a non-conformally
   Chen   ideal   submanifold   of  dimension $ n $
   in the Euclidean space  $ \mathbb{E}^{n+m} $,
   $n \geq 4$, $m \geq 1$.
   Then there exists a real valued function $ L $ on $ M $
   such that
 $$
  R \cdot C - C \cdot R = L Q(g, g \wedge S)   ,
  $$
 if and only if
  $ M  $ is not    minimal,
   and there exists
    an orthonormal  tangent  framefield
  $  \left\{e_1 , \cdots , e_n\right\}    $ and
     an orthonormal  normal  framefield
$  \left\{\xi_1 , \cdots , \xi_m\right\}    $
  such that the shape operators
  $A_\alpha :=  A_{\xi_\alpha} , \quad 1 \leq \alpha \leq m$,
are given by
     $$
 % \begin{equation}\label{eqbli0ggwS01}
  A_1 =  \begin{pmatrix}
  a & 0 & 0 & \cdots & 0 \\
  0 & a & 0 & \cdots & 0 \\
  0 & 0 & 2a & \cdots & 0 \\
  \vdots & \vdots & \vdots & \ddots & \vdots \\
  0 & 0 & 0 & \cdots & 2a \\
  \end{pmatrix}
  ,  \quad
 A_\beta =  \begin{pmatrix}
  c_\beta & d_\beta & 0 & \cdots & 0 \\
  d_\beta & -c_\beta & 0 & \cdots & 0 \\
  0 & 0 & 0 & \cdots & 0 \\
  \vdots & \vdots & \vdots & \ddots & \vdots \\
  0 & 0 & 0 & \cdots & 0 \\
  \end{pmatrix}
, \
  \beta \geq 2 ,
    $$
 where $ a $,
 $ c_\beta $ , $ d_\beta $ (for $ 2 \leq \beta \leq m $)
  are real functions on $ M $
 such that
  $$
    \sum_{\beta = 2}^m \left(c^2_\beta + d^2_\beta\right)
   =  \dfrac{2n^2 - 5n + 1}{n-1} a^2 > 0,\ \ \
L =    - \dfrac{2a^2}{n-2}  .
    $$
    \noindent  In this case,
  $M$ is  a properly  pseudosymmetric manifold 
  (see  Theorem 11.4).

   \end{cor}

   \begin{cor}  \label{r.c-c.r=LQ(SgWS)}  
   {\cite[Corollary 6] {DP-TVZ}}   % coroll 6     
   Let $ M $ be a non-conformally flat
   Chen   ideal   submanifold   of  dimension $ n $
   in the Euclidean space  $ \mathbb{E}^{n+m} $,
    $n \geq 4$, $m \geq 1$.
    Then there exists a real valued function $ L $ on $ M $
   such that
 $$
  R \cdot C - C \cdot R = L Q(S  , g \wedge S)
  ,
  $$
 if and only if  one has one of the two cases which follow:
\newline
(i) either   $ M $ is minimal, and
    there exists
    an orthonormal  tangent  framefield
  $  \left\{e_1 , \cdots , e_n\right\}    $ and
     an orthonormal  normal  framefield
$  \left\{\xi_1 , \cdots , \xi_m\right\}    $
  such that the shape operators
  $
  A_\alpha :=  A_{\xi_\alpha} , \quad 1 \leq \alpha \leq m $,
are given by:
    $$
 %  \begin{equation}\label{eqbli0QSgWS1}
  A_1 =  \begin{pmatrix}
  a & 0 & 0 & \cdots & 0 \\
  0 & -a & 0 & \cdots & 0 \\
  0 & 0 & 0 & \cdots & 0 \\
  \vdots & \vdots & \vdots & \ddots & \vdots \\
  0 & 0 & 0 & \cdots & 0 \\
  \end{pmatrix}
  ,  \quad
 A_\beta =  \begin{pmatrix}
  c_\beta & d_\beta & 0 & \cdots & 0 \\
  d_\beta & -c_\beta & 0 & \cdots & 0 \\
  0 & 0 & 0 & \cdots & 0 \\
  \vdots & \vdots & \vdots & \ddots & \vdots \\
  0 & 0 & 0 & \cdots & 0 \\
  \end{pmatrix}
, \
  \beta \geq 2 ,
    $$
 where $ a $,
 $ c_\beta $ , $ d_\beta $ (for $ 2 \leq \beta \leq m $)
  are real functions on $ M $
 such that
  $$
   \sum_{\beta = 2}^m \left(c^2_\beta + d^2_\beta\right)
   = - a^2 -  \inf K ,\ \ \
    L =  - \dfrac{n-3}{(n-1)(n-2)} ,$$
(ii)   or    $ M $ is not  minimal, and
  there exists
    an orthonormal  tangent  framefield
  $  \left\{e_1 , \cdots , e_n\right\}    $ and
     an orthonormal  normal  framefield
$  \left\{\xi_1 , \cdots , \xi_m\right\}    $
  such that the shape operators
  $A_\alpha :=  A_{\xi_\alpha} , \quad 1 \leq \alpha \leq m$,
are given by
     $$
  %  \begin{equation}\label{eqbli0QSgWS2}
  A_1 =  \begin{pmatrix}
    a    &    0   &   0    & \cdots & 0 \\
    0    &    a   &   0    & \cdots & 0 \\
    0    &   0    &   2a   & \cdots & 0 \\
  \vdots & \vdots & \vdots & \ddots & \vdots \\
    0    &   0    &   0    & \cdots & 2a \\
  \end{pmatrix}
  ,  \quad
 A_\beta =  \begin{pmatrix}
  c_\beta & d_\beta & 0 & \cdots & 0 \\
  d_\beta & -c_\beta & 0 & \cdots & 0 \\
  0 & 0 & 0 & \cdots & 0 \\
  \vdots & \vdots & \vdots & \ddots & \vdots \\
  0 & 0 & 0 & \cdots & 0 \\
  \end{pmatrix}
, \
 \beta \geq 2 ,
   $$
 where $ a $,
 $ c_\beta $ , $ d_\beta $ (for $ 2 \leq \beta \leq m $)
  are real functions on $ M $
 such that
 $$ \sum_{\beta = 2}^m \left(c^2_\beta + d^2_\beta\right) 
 = (2n - 3) a^2 ,\ \ \
L =
 \dfrac{1}{2\left(n-1\right)\left(n-2\right)}   .
$$
In the first  case, $ M  $ is a semisymmetric
    manifold
  (see   Theorem 11.3).
   In the second case, $ M  $ is a properly  pseudosymmetric
        manifold  (see Theorem 11.4).
 
   \end{cor}

   \begin{cor}  \label{r.c-c.r MINIM} 
   {\cite[Corollary 7] {DP-TVZ}}  % coroll 7    
   Let $ M $ be a non-conformally flat Chen ideal
     submanifold   of  dimension $ n $
   in the Euclidean  space  $ \mathbb{E}^{n+m} $,
   $n \geq 4$, $m \geq 1$.
   If $ M $ is minimal, then
$$
 R \cdot C - C \cdot R =
 - \dfrac{n-3}{(n-1)(n-2)} Q(S , g \wedge S) .
$$
     \end{cor}

  \begin{cor}  \label{r.c-c.r=0} {\cite[Corollary 8] {DP-TVZ}}  % Coroll 8
 A  Chen   ideal   submanifold $ M $ of dimension
 $ n \geq 4  $
    in the Euclidean space  $ \mathbb{E}^{n+m} $
  satisfies
 the curvature condition
 $ R \cdot C - C \cdot R = 0$
 if and only if
  $ M $ is conformally flat.
   \end{cor}

According to {\cite[eq. (28)] {DP-TVZ}}
(see also {\cite[Theorem 4.1] {2023_DGHP-TZ}}), 
  as an immediate
   consequence of Theorem 11.11,
    for Chen ideal and Roter submanifolds,
    we express the difference  tensor
$ R \cdot C - C \cdot R $  as a
linear combination of the Tachibana tensors
$ Q(g , g \wedge S) $,  $  Q(S , g \wedge S)  $,
 in terms of
 $ \inf K $  and the scalar curvature $ \tau $.

   \begin{cor}  \label{r.c-c.r PseudoSym} 
   {\cite[Corollary 9] {DP-TVZ}}  % coroll 9    
   Let $ M $ be a non-conformally flat Chen ideal
     submanifold   of  dimension $ n $
   in the Euclidean space  $ \mathbb{E}^{n+m} $,
   $n \geq 4$, $m \geq 1$.
   If $ M $ is a Roter space,  then
\begin{eqnarray*}
R \cdot C - C \cdot R &=&
 -  \dfrac{2\left(\tau  - 2\, \inf K\right) 
 \left(\tau + 2(n - 2)\, \inf K \right)}{(n-1)
 \left(\tau - 2n\, \inf K \right )^2}Q(g , g \wedge S)\\
& &
-
  \dfrac{2(n  - 1)\, \inf K
\left(\tau + 2(n - 4)\, \inf K\right)}{(n-2)
\left(\tau - 2n \, \inf K \right)^2}Q(S , g \wedge S) .
\end{eqnarray*}
     \end{cor}

\newpage

%\vspace{5mm}

\noindent
\footnotesize{Ryszard Deszcz\\
retired employee of the
Department of Applied Mathematics\\
Wroc\l aw University of Environmental and Life Sciences\\ 
Grunwaldzka 53, 50-357 Wroc\l aw, Poland}\\
{\sf E-mail: Ryszard.Deszcz@upwr.edu.pl}\\
\textbf{ORCID ID: 0000-0002-5133-5455} 
\newline

\noindent
\footnotesize{Ma\l gorzata G\l ogowska\\
Department of Applied Mathematics \\
Wroc\l aw University of Environmental and Life Sciences\\ 
Grunwaldzka 53, 50-357 Wroc\l aw, Poland}\\
{\sf E-mail: Malgorzata.Glogowska@upwr.edu.pl}\\
\textbf{ORCID ID: 0000-0002-4881-9141}
\newline

\noindent
\footnotesize{Marian Hotlo\'{s}\\
retired employee of the Department of Applied Mathematics\\ 
Wroc{\l}aw University of Science and Technology\\
Wybrze\.{z}e Wyspia\'{n}skiego 27\\ 
50-370 Wroc{\l}aw, Poland}\\ 
{\sf E-mail: Marian.Hotlos@pwr.edu.pl}\\
\textbf{ORCID ID: 0000-0002-4165-4348}
\newline

\noindent
\footnotesize{Miroslava Petrovi\'{c}-Torga\v{s}ev\\
Department of Sciences and Mathematics\\ 
State University of Novi Pazar\\ 
Vuka Karad\v{z}i\'{c}a 9\\ 
36300 Novi Pazar, Serbia}\\
{\sf E-mail: mirapt@kg.ac.rs}\\
\textbf{ORCID ID: 0000-0002-9140-833X}
\newline

\noindent
\footnotesize{Georges Zafindratafa\\
Professeur \'{E}m\'{e}rite of the 
Laboratoire 
de Math\'{e}matiques pour l'Ing\'{e}nieur (LMI)\\
Universit\'{e} Polytechnique Hauts-de-France\\
59313 Va\-len\-cien\-nes cedex 9, France}\\
{\sf E-mail: Georges.Zafindratafa@uphf.fr}\\
\textbf{ORCID ID: 0009-0001-7618-4606}

\end{document}